\documentclass[a4paper, 10pt, DIV = 13]{scrartcl}
\usepackage[utf8]{inputenc}
\usepackage[english]{babel}
\usepackage{tikz}
\usetikzlibrary{calc,angles,quotes,arrows.meta}
\usepackage{amsmath,amsfonts,amssymb,amsthm,mathrsfs}
\usepackage{lmodern}
\usepackage{hyperref}
\usepackage{xcolor}
\usepackage{nicefrac}
\usepackage{adjustbox}
\usepackage[normalem]{ulem}
\hypersetup{breaklinks=true}
\usepackage{framed}
\usepackage{mathtools}
\mathtoolsset{showonlyrefs}
\usepackage{longtable}
\usepackage{graphicx}
\usepackage{subcaption}   
\numberwithin{equation}{section}
\newtheorem{proposition}{Proposition}[section]
\newtheorem{lemma}[proposition]{Lemma}

\newtheorem{claim}[proposition]{Claim}

\theoremstyle{definition}
\newtheorem{theorem}[proposition]{Theorem}
\newtheorem*{theorem*}{Theorem}
\newtheorem{definition}[proposition]{Definition}
\newtheorem{remark}[proposition]{Remark}
\newtheorem{corollary}[proposition]{Corollary}

\usepackage{enumitem}
\setlist{nosep}
\setlist{noitemsep}
\setlist{leftmargin=*}
\newcommand{\R}{\mathbb{R}}

\newcommand{\N}{\mathbb{N}}

\renewcommand{\epsilon}{\varepsilon}
\newcommand{\hal}{\frac{1}{2}}

\newcommand{\sineb}{\mathrm{Sine}_\beta}

\newcommand{\Ee}{\mathcal{E}}
\newcommand{\Fb}{\mathcal{F}_\beta}
\newcommand{\Fbint}{\mathcal{F}^{\mathrm{int}}_\beta}
\newcommand{\W}{\mathcal{W}}
\newcommand{\Wint}{\mathcal{W}^{\mathrm{int}}}

\newcommand{\bX}{\mathbf{X}}
\newcommand{\X}{\mathsf{X}}
\newcommand{\bXN}{\bX_N}

\newcommand{\xx}{\mathsf{x}}
\newcommand{\bXNx}{\mathbf{X}_{N, \xx}}
\newcommand{\bXNz}{\mathbf{X}_{N, 0}}

\newcommand{\Pstat}{\mathcal{P}_{\mathrm{stat}}}
\newcommand{\Pp}{\mathrm{P}}
\newcommand{\E}{\mathbb{E}}
\renewcommand{\P}{\mathbb{P}}
\newcommand{\La}{\Lambda}

\newcommand{\Dis}{\mathsf{Dis}}

\newcommand{\Dlog}{\mathcal{D}_{\mathrm{log}}}
\renewcommand{\O}{\mathcal{O}}
\newcommand{\Cc}{\mathtt{C}}
\newcommand{\cc}{\mathtt{c}}
\newcommand{\Fluct}{\mathrm{Fluct}}
\newcommand{\dd}{\mathtt{d}}
\newcommand{\Ent}{\mathsf{Ent}}

\newcommand{\ZNbC}{\mathrm{Z}_{N, \beta}^{\mathtt{c}}}
\newcommand{\ZNbGC}{\mathrm{Z}_{N, \beta}^{\mathtt{gc}}}
\newcommand{\PNbGC}{\mathbb{P}_{N, \beta}^{\mathtt{gc}}}
\newcommand{\PNbC}{\mathbb{P}_{N, \beta}^{\mathtt{c}}}
\newcommand{\Poi}{\mathsf{Poi}}
\newcommand{\Bin}{\mathsf{Bin}}
\newcommand{\Leb}{\mathsf{Leb}}
\newcommand{\Bulk}{\mathsf{Bulk}}
\newcommand{\PNb}{\mathbb{P}_{N,\beta}}
\newcommand{\ENb}{\mathbb{E}_{N,\beta}}
\newcommand{\ZNb}{\mathrm{Z}_{N,\beta}}
\newcommand{\HN}{\mathsf{H}_N}
\newcommand{\rhosc}{\rho_{\mathrm{sc}}}
\newcommand{\vrho}{\varrho}
\newcommand{\Conf}{\mathsf{Conf}}
\newcommand{\g}{\mathsf{g}}
\newcommand{\f}{\mathsf{f}}
\newcommand{\M}{\mathsf{M}}
\newcommand{\K}{\mathrm{K}}

\newcommand{\bY}{\mathbf{Y}}
\newcommand{\ind}{\mathtt{1}}
\newcommand{\sm}[2]{\delta^{(#2)}_{#1}}
\newcommand{\tr}{\tilde{\mathsf{r}}}
\newcommand{\Elec}{\mathsf{Elec}}

\newcommand{\G}{\mathsf{G}}
\newcommand{\Qq}{\mathrm{Q}}
\newcommand{\Pgen}[3]{\mathbb{P}_{#1,#2,#3}}
\newcommand{\Zgen}[3]{\mathrm{Z}_{#1,#2,#3}}
\newcommand{\Pc}[2]{\mathbb{P}^{\mathtt{c}}_{#1,#2}}
\newcommand{\Zc}[2]{\mathrm{Z}^{\mathtt{c}}_{#1,#2}}
\newcommand{\red}[1]{\textcolor{black}{#1}}

\DeclareMathOperator{\supp}{supp}

\title{\LARGE{DLR variational principle for the one-dimensional log-gas and another proof of existence of $\sineb$}}
\date{\small{\today}}
\author{\large{Sylvain Chabredier\thanks{ Université de Paris-Cité, CNRS, MAP5 UMR 8145, F-75006 Paris, France \texttt{sylvain.chabredier@ens.psl.eu}},  Thomas Leblé}\thanks{Department of Mathematics, Kyiv School of Economics, Kyiv, Ukraine \texttt{thomas.leble@math.cnrs.fr}}}
\begin{document}

\maketitle
\begin{abstract}
We show that every stationary solution of the Dobrushin--Lanford--Ruelle (DLR) equations for the one-dimensional log-gas with a quantitative control on its number variance minimizes the free energy functional. This yields a DLR variational principle for this continuous long-range system.

Moreover, we show that any limit point of the finite-volume log-gases is stationary and satisfies this condition. This allows us to give a new, physical proof of the existence of a microscopic limit for the Gaussian $\beta$-ensembles, known as the $\sineb$ point process.
\end{abstract}


\section{Introduction}
\subsection{The Gaussian \texorpdfstring{$\beta$}{beta}-ensembles}
\label{sec:intro-GbE}
The Gaussian $\beta$-ensembles form a family, indexed by $\beta \in (0, + \infty)$, of probability densities given by:
\begin{equation}
\label{def:PNb}
\dd \PNb(x_1, \dots, x_N) := \frac{1}{\ZNb} \prod_{1 \leq i < j \leq N} |x_i - x_j|^\beta \, \exp\left( - \frac{\beta}{4N} \sum_{i=1}^N x_i^2 \right) \dd x_1 \dots \dd x_N,
\end{equation}
where $\ZNb$ is a normalizing constant and $\dd x_1 \dots \dd x_N$ is the Lebesgue measure on $\R^N$. They describe the joint law of eigenvalues in important models of random matrix theory. Historically speaking, the three \emph{classical} values were $\beta \in \{1, 2, 4\}$~\cite{dyson1962threefold}, but in a celebrated work, Dumitriu and Edelman~\cite{dumitriu2002matrix} constructed, for every $\beta > 0$, tridiagonal symmetric matrices with independent entries whose eigenvalues have joint law~\eqref{def:PNb}. This constitutes the “random matrix interpretation” of the measure $\PNb$.

One major achievement of the theory of $\beta$-ensembles is the proof by Valkó--Virág~\cite{valko2009continuum} of the existence of a limiting point process for all $\beta > 0$, called the $\sineb$ process, see Section~\ref{sec:intro-sineb}.

\paragraph{Log-gas interpretation.} Dyson~\cite{dyson1962statistical} was the first to suggest rewriting~\eqref{def:PNb} as a \emph{Gibbs measure}
\begin{equation}
\label{def:HN}
\dd \PNb(\X_N) = \frac{1}{\ZNb} \exp\left( - \beta \HN(\X_N) \right) \dd \X_N, \qquad \HN(\X_N) := \hal \sum_{1 \leq i \neq j \leq N} - \log |x_i - x_j|  + \sum_{i=1}^N \frac{x_i^2}{4N},
\end{equation}
where $\X_N = (x_1, \dots, x_N) \in \R^N$. In this “statistical physics interpretation”, $\PNb$ describes $N$ particles on the line, with equal charges, repelling each other through the logarithmic interaction potential, and confined by a quadratic external field whose strength depends on $N$. The parameter $\beta$ plays the role of an \emph{inverse temperature}, $\HN$ is the \emph{energy} and $\ZNb$ the \emph{partition function}.

\medskip

Since the seminal paper~\cite{sandier2015log}, this point of view has been used to study the microscopic scale. The goal of the present paper is to give a new proof of the existence of a limiting point process, using only the statistical physics perspective on $\beta$-ensembles. This approach could be implemented in other Gibbsian models for which no random matrix theory interpretation is available.

\paragraph{Scaling.} With our scaling convention, in which the quadratic potential has strength $\frac{1}{N}$, the typical distance between neighboring eigenvalues is of order $1$, and the total size of the system is proportional to $N$. A classical consequence of Wigner's semicircle law~\cite{wigner1958distribution} (see e.g.~\cite{anderson2010introduction}) is that the rescaled empirical measure $\frac{1}{N} \sum_{i=1}^N \delta_{x_i / N}$ converges a.s. to the semicircle law, whose density is given by:
\begin{equation}
\label{def:rhosc}
\rhosc(v) := \frac{1}{2\pi} \sqrt{\left(4 - v^2 \right)_+}.
\end{equation}
 Near $N v \in \R$ with $v \in (-2, 2)$, there are approximately $\rhosc(v)$ eigenvalues per unit length. We write:
\begin{equation*}
\vrho := \rhosc(0) = \frac{1}{\pi}.
\end{equation*}
We distinguish two lengthscales: the microscopic one, of order $1$, and the macroscopic one, of order $N$.
\subsection{The \texorpdfstring{$\sineb$}{Sine-beta} point process}
\label{sec:intro-sineb}

\paragraph{Local configurations.} The microscopic behavior of the system near a given $\xx \in \R$ is encoded by the local point configuration seen from this point. For $\xx$ in the \emph{bulk} $(-2, 2)$, we let
\begin{equation}
\label{def:bXNx}
    \bXNx := \sum_{i=1}^N \delta_{\rhosc(\xx) (x_i - N \xx)},
\end{equation}
where $(x_1, \dots, x_N)$ is distributed according to $\PNb$. This is the configuration seen from $N\xx$, rescaled by $\rhosc(\xx)$ to have intensity close to one near $0$. For simplicity, we will  mostly work with $\xx = 0$ (our arguments extend to any $\xx$ in the bulk), namely with:
\begin{equation}
\label{def:bXNz}
    \bXNz = \sum_{i=1}^N \delta_{\vrho x_i}.
\end{equation}
Let $\Conf$ be the space of locally finite point configurations on $\R$, endowed with the vague topology, for which two configurations are close if they “look the same on a large interval”. A finite configuration can thus be close to an infinite one. The law of $\bXNx$ is a probability measure on $\Conf$, i.e. a point process.

\paragraph{Existence of a limit.}
For $\beta \in \{1, 2, 4\}$, proving the existence of a limit is a fairly simple analytic task due to a very specific integrable structure. For general $\beta$, there is no such structure. In fact, even the existence of \emph{limit points}, i.e. of tightness for the law of $(\bXNx)_N$, is a delicate question: the configuration $\bXNx$ could a priori have anywhere between $0$ and $N$ points in a fixed interval, and controlling this requires understanding the system \emph{down to the microscopic scale}. In a major breakthrough \cite{valko2009continuum}, Valkó and Virág proved that this limit \emph{exists and is the same, up to a scaling, at every point of the bulk}.
\begin{theorem*}
For every $\xx \in (-2, 2)$, the law of $\bXNx$ converges, as $N \to \infty$, to the law of a stationary point process of intensity one, called $\sineb$, which does not depend on $\xx$.
\end{theorem*}

 The approach of~\cite{valko2009continuum} starts from the tridiagonal matrix model of~\cite{dumitriu2002matrix}, and analyzes its eigenvalue equation through a phase function, in the spirit of Sturm--Liouville oscillation theory. The work of Killip and Stoiciu~\cite{killip2009eigenvalue}, based on a “pentadiagonal” matrix model~\cite{killip2004matrix} and on a similar analysis of phases, proved the existence of a limit for a very closely related model, the Circular $\beta$-ensemble, whose physical interpretation is a log-gas with periodic boundary conditions. Those limits are actually the same~\cite{nakano2014level}. 

 Several other descriptions of $\sineb$ have been given since: as the spectrum of a random Dirac operator~\cite{valko2017operator}, as the zero set of a random analytic function~\cite{valko2022many}, among others. They have been used to study $\sineb$ itself, beyond its mere existence. None of them uses the interpretation of $\PNb$ as the Gibbs measure of a system of particles with logarithmic interaction. 

\medskip

In recent works~\cite{dereudre2021dlr,erbar2021one}, the $\sineb$ process itself was given a statistical physics description, but the following natural question remained:
\begin{center}
\emph{Can we prove the existence of a limit for $\bXNx$, using only the Gibbsian structure of~\eqref{def:HN}?}
\end{center}

As sketched in~\cite[Section 1.4]{lebleHDR}, the following chain of results would lead to a positive answer:
\begin{enumerate}
	\item Prove \emph{local laws} for the finite-volume model, valid down to the microscopic scale. This ensures the existence of limit points for the law of $(\bXNx)_N$, with good controls on the number of points.
	\\ This was done by Peilen~\cite{peilen2022local}.
	\item Show that the limit points satisfy \emph{DLR equations}, i.e. are infinite-volume Gibbs measures. \\ This was essentially done in~\cite{dereudre2021dlr}.
	\item Show that the limit points are stationary and retain good properties from the finite-volume model.
	\item Show that every stationary solution (with good properties) of the DLR equations \emph{minimizes the free energy functional}. This is called a “DLR variational principle”.
	\item Use the fact, proven in~\cite{erbar2021one}, that this free energy functional has a \emph{unique} stationary minimizer.
	\item Conclude: the sequence $(\bXNx)_N$ has only one limit point, and thus converges!
\end{enumerate}
The goal of our paper is to complete this program by implementing steps 3 and 4.

\medskip

Compared to the “random matrix theory” approach, this has arguably two upsides:
\begin{itemize}
	\item It can be implemented for other, non-RMT-related systems, as was done recently  in \cite{boursier2026existence} for the \emph{two}-dimensional log-gas at inverse temperature $0 < \beta \leq 2$ — existence of a limiting point process was previously only known for $\beta = 2$.
	\item It “directly” yields universality of the limiting point process, retrieving some of the results of \cite{bourgade2014universality,bekerman2015transport}. Indeed, our argument does not use the specific form of the Gaussian potential: it relies on the local laws of \cite{peilen2022local}, which hold for quite general confining potentials.
\end{itemize}

\subsection{DLR equations}
\label{sec:intro-DLR}

\subsubsection*{Heuristics}
In finite volume, $\PNb$ is a Gibbs measure: its density is of the form $\exp(-\beta \times \text{energy})$ with respect to a reference “volume” measure. It is natural to ask whether the limit points of $(\bXNx)_N$ admit a similar representation. The right answer goes through the \emph{Dobrushin--Lanford--Ruelle} (DLR) formalism: instead of looking for a global density, one asks for a Gibbsian description of the \emph{local, conditional} distributions: for every bounded set $\La$, the law of the configuration in $\La$, given the configuration $\bX_{\La^c}$ outside $\La$, should be of the form
\begin{equation}
\label{eq:DLR-heuristic}
\dd \Pp\left( \bX'_\La \, | \, \bX_{\La^c} \right) \propto \exp\left( - \beta \left( \mathsf{Energy}(\bX'_\La) + \mathsf{Energy}(\bX'_\La, \bX_{\La^c}) \right) \right) \, \dd \Poi_\La(\bX'_\La),
\end{equation}
where $\Poi_\La$ is a Poisson point process in $\La$. The self-energy of $\bX'_\La$ is not mysterious, since $\La$ is bounded. The puzzling term is the interaction between the points in $\La$ and the \emph{infinitely many} points outside $\La$:
\begin{equation}
\label{eq:naive-interaction}
\mathsf{Energy}(\bX'_\La, \bX_{\La^c}) \overset{??}{:=} \iint_{\La \times \La^c} - \log |x - y| \, \dd (\bX'_\La - \Leb)(x) \, \dd (\bX_{\La^c} - \Leb)(y),
\end{equation}
where the Lebesgue measure $\Leb$ plays the role of the neutralizing background. Because the interaction is long-range, and $\La^c$ is unbounded, there is no reason for~\eqref{eq:naive-interaction} to converge.

\subsubsection*{Move functions, DLR equations}
As explained in~\cite{dereudre2021dlr}, it is enough to make sense of the \emph{difference} between the interactions of two configurations $\bX', \bX$ with the exterior, namely of the \emph{move function}
\begin{equation}
\label{eq:move-intro}
\M_\La(\bX', \bX) := \iint_{\La \times \La^c} - \log |x - y| \, \dd (\bX'_{\La} - \bX_\La)(x) \, \dd (\bX_{\La^c} - \Leb)(y),
\end{equation}
which represents the energy cost of moving the points of $\bX$ inside $\La$ to the positions given by $\bX'$, while feeling the effect of the (fixed) exterior configuration and of the background. If $\bX'$ and $\bX_\La$ have the \emph{same number of points}, the signed measure $\bX' - \bX_\La$ has total mass zero, and this neutrality provides an additional cancellation which makes~\eqref{eq:move-intro} converge. This leads to the \emph{canonical} DLR equations, which are stated as follows: conditionally on $\bX_{\La^c}$ \emph{and on the number $n$ of points in $\La$},
\begin{equation}
\label{DLR}
\tag{DLR}
\dd \Pp\left( \bX'_\La \, | \, \bX_{\La^c}, \, \bX(\La) = n \right) = \frac{1}{\K_\La(\bX)} \exp\left( - \beta \left( F_\La(\bX'_\La) + \M_\La(\bX'_\La, \bX) \right) \right) \, \dd \Bin_{\La, n}(\bX'_\La),
\end{equation}
where $F_\La$ is the self-energy in $\La$ (see~\eqref{def:FLa}), $\Bin_{\La, n}$ is the law of $n$ independent uniform points in $\La$, and $\K_\La(\bX)$ is a normalizing constant. The precise definition is given in Section~\ref{sec:DLR-def}. Fixing the number of points is natural here: the limit processes are \emph{number-rigid} in the sense of~\cite{ghosh2017rigidity}, i.e. the number of points in $\La$ is a function of the exterior configuration, so that conditioning on it adds no information~\cite{dereudre2021dlr, chhaibi2018rigidity}. 

\subsubsection*{The classical theory, and the variational principle(s)}
The DLR formalism, introduced at the end of the 1960s, is the classical framework for infinite-volume Gibbs measures. It is very well understood for lattice spin systems with absolutely summable interactions, see~\cite{georgii2011gibbs, friedli2017statistical}, and has been developed for continuous systems of particles with sufficiently short-range (e.g. superstable and regular) interactions, see e.g.~\cite{georgii1976canonical, preston1976random, dereudre2019introduction}.

One of the great classical results of the theory is the \emph{variational principle}, which connects the two possible definitions of an equilibrium state in infinite volume: \emph{For translation-invariant measures, being a solution of the DLR equations is equivalent to minimizing the \emph{free energy} functional.} This is the infinite-volume analogue of the classical fact, sometimes known as the \emph{Gibbs variational principle}, that in finite volume, the Gibbs measure is the unique minimizer of “$\beta \times$ energy $+$ relative entropy”. We refer to~\cite[Section 6.9]{friedli2017statistical} or~\cite[Chapter 8]{RAS15} for lattice systems, and to \cite{jahnel2025variational} for one of the very few results in the continuum.

The classical theory is developed for “short-range” interactions, which ensures that the interaction energy of a large box with its exterior is negligible compared to the volume of the box. This is the main ingredient of the proof. To the best of our knowledge, the DLR variational principle has \emph{not been proven for continuous systems with long-range interactions}, such as log-gases, or Coulomb and Riesz gases.

\subsection{The free energy functional}
\label{sec:intro-free-energy}
In~\cite{leble2017large}, an infinite-volume free energy functional for the one-dimensional log-gas was introduced, and was identified as the rate function of a large deviation principle for the empirical field at microscopic scale. It is defined on the set $\Pstat$ of stationary point processes on $\R$ of intensity one, and has the form
\begin{equation}
\label{def:Fb}
\Fb(\Pp) := \beta \W(\Pp) + \Ee(\Pp).
\end{equation}
\begin{itemize}
	\item $\Ee$ is the \emph{specific relative entropy} with respect to the Poisson point process, which is a fairly easy object to define. For $\Pp \in \Pstat$, let $\Pp_{\La_R}$ and $\Poi_{\La_R}$ be the restrictions of $\Pp$ and of a Poisson point process of intensity one to $\La_R := [-\frac{R}{2}, \frac{R}{2}]$. Then
	\begin{equation*}
	\Ee(\Pp) := \lim_{R \to \infty} \frac{1}{R} \Ent\left( \Pp_{\La_R} \, | \, \Poi_{\La_R} \right),
	\end{equation*}
	where $\Ent(\mu | \nu) := \int \log \frac{\dd \mu}{\dd \nu} \dd \mu$ if $\mu \ll \nu$ and $+ \infty$ otherwise. The limit exists by super-additivity, see~\cite{friedli2017statistical} for lattice systems and~\cite{georgii1993large} for point processes.
	\item $\W$ is the \emph{electric energy}, or renormalized energy, of Sandier and Serfaty~\cite{sandier2015log}, adapted to point processes in~\cite{leble2017large}. Its definition uses the electric field generated by the points on the line and a uniform background of opposite charge, and requires a “renormalization” procedure. We will not need its precise definition here, we refer to~\cite{serfaty2024lectures} for a textbook presentation.
\end{itemize}
In this paper, we will rather use an \emph{intrinsic} energy, introduced in~\cite[Def. 3.3]{leble2016logarithmic}, which has an explicit expression in terms of the pairwise interactions. For a bounded interval $\La$ and a configuration $\bX$ in $\La$, let
\begin{equation}
\label{def:FLa}
F_\La(\bX) := \hal \iint_{(x,y) \in \La^2, \, x \neq y} - \log |x - y| \, \dd (\bX - \Leb)(x) \, \dd (\bX - \Leb)(y)
\end{equation}
be the logarithmic energy of the system made of the charges of $\bX$ and of the uniform background on $\La$. For $\Pp \in \Pstat$, let
\begin{equation}
\label{def:Wint}
\Wint(\Pp) := \liminf_{L \to \infty} \frac{1}{L} \E_\Pp\left[ F_{\La_L}(\bX \cap \La_L) \right], \qquad \Fbint(\Pp) := \beta \Wint(\Pp) + \Ee(\Pp).
\end{equation}
The connection between $\W$ and $\Wint$ was studied in~\cite{leble2016logarithmic}, and we recall it in the proof of Proposition~\ref{prop:hPNbGC}. The two do not coincide in general: for the Poisson point process, $\Wint = 0$ while $\W = + \infty$. 

\paragraph{Minimizers.} The functional $\Fb$ is lower semi-continuous with compact sub-level sets~\cite{leble2017large}, hence it has minimizers, and the process $\sineb$ constructed by~\cite{valko2009continuum} is a minimizer~\cite[Corollary 1.2]{leble2017large}. Moreover:
\begin{theorem*}[\cite{erbar2021one}]
$\Fb$ has a unique minimizer in $\Pstat$.
\end{theorem*}
The proof of \cite{erbar2021one} relies on a notion of optimal transport between stationary point processes in infinite volume and on the \emph{displacement convexity} of the free energy, the strict convexity coming from the strict convexity of the logarithmic interaction. These techniques have recently been improved, and extended to one-dimensional Riesz interactions, in~\cite{huesmann2025gap}.

\subsection{Main statements}
\label{sec:intro-main}
Recall that $\Pstat$ denotes the set of stationary point processes of intensity one on $\R$. Our results hold for the following class of point processes, see Definition~\ref{def:good}: we say that $\Pp$ is \emph{good} if there is a constant $\Cc_\Pp$ such that, for every interval $I$ of length $|I| \geq 1$,
\begin{equation*}
\E_\Pp\left[ \left( \bX(I) - |I| \right)^2 \right] \leq \Cc_\Pp \log^2(2 + |I|).
\end{equation*}
This is a quantitative assumption on the fluctuations of the number of points at large scales. We will see that it is satisfied by all the limit points of the finite-volume log-gases. For our purposes, it could be relaxed, see Remark~\ref{rem:good-relax}.

\medskip

We first return to the finite-volume model, and study its limit points in the thermodynamic limit.
\begin{theorem}[DLR equations for the limit points]
\label{theo:Convergence}
Consider the sequence $(\bXNz)_N$ as in~\eqref{def:bXNz}. It is tight in the exp-local topology (Definition~\ref{def:exp-local}), and every limit point $\Pp$ satisfies:
\begin{enumerate}
	\item $\Pp$ is good (Section~\ref{sec:tightness}),
	\item $\Pp$ is a solution of~\eqref{DLR} in the sense of Definition~\ref{def:DLR} (Section~\ref{sec:DLR-limit}).
\end{enumerate}
\end{theorem}

Next, we show that the DLR equations force stationarity, by adapting the classical method of~\cite{frohlich1981absence}.
\begin{theorem}[Stationarity]
\label{theo:TI-intro}
Every good solution of~\eqref{DLR} is stationary.
\end{theorem}
This is proven in Section~\ref{sec:stationarity}. Note that the finite-volume model is not translation-invariant, so that stationarity of the limit points is not obvious a priori. This is an important difference from the circular $\beta$-ensemble, which is rotation-invariant.

\medskip

Our main new result is the following.
\begin{theorem}[DLR variational principle]
\label{theo:VarProp}
Every good stationary solution of {\eqref{DLR}} is a minimizer of $\Fb$. 
\end{theorem}
By Theorems~\ref{theo:Convergence} and~\ref{theo:TI-intro}, good solutions of~\eqref{DLR} exist, and they are stationary. Combined with Theorem~\ref{theo:VarProp} and with the uniqueness result of~\cite{erbar2021one} for minimizers of $\Fb$, we deduce:
\begin{corollary}[Uniqueness for DLR equations]
\label{cor:Uniqueness}
There is a unique good solution of~\eqref{DLR}, and it is the unique minimizer of $\Fb$.
\end{corollary}
\begin{remark}[Tail triviality]
\label{rem:tail}
Corollary~\ref{cor:Uniqueness} implies that the tail $\sigma$-algebra of the unique good solution $\Pp$ is trivial. Tail triviality of $\sineb$ for all $\beta > 0$ was recently proven in~\cite{assiotis2026tail} by a different route.
\end{remark}

Putting things together, we deduce that:
\begin{theorem}
\label{theo:SineBeta}
The sequence $(\bXNz)_N$ converges, in the exp-local topology of Definition~\ref{def:exp-local}, to a good stationary point process, which is the unique minimizer of $\Fb$.
\end{theorem}
\begin{proof}
By Theorem~\ref{theo:Convergence}, the sequence is tight in the exp-local topology, and every limit point is a good solution of~\eqref{DLR}. By Corollary~\ref{cor:Uniqueness}, all the limit points coincide with the unique minimizer of $\Fb$, hence the sequence converges.
\end{proof}
This yields a new, “physical” proof of the existence of a microscopic limit in the bulk of $\beta$-ensembles. Since the $\sineb$ process constructed in~\cite{valko2009continuum} minimizes $\Fb$~\cite[Corollary 1.2]{leble2017large}, this limit is $\sineb$.

\medskip

Our proof does not use any matrix model or integrable structure. The external inputs are the local laws of Peilen~\cite{peilen2022local}, the “intrinsic” energy of~\cite{leble2016logarithmic}, the free energy functional of~\cite{leble2017large}, and the uniqueness of its minimizer~\cite{erbar2021one}, which are all obtained by “statistical physics” methods.

\paragraph{Comments.}
\begin{enumerate}
	\item We state Theorems~\ref{theo:Convergence} and~\ref{theo:SineBeta} for the configuration seen from the origin, in order to simplify the notation. The same arguments would apply to $\bXNx$ for any $\xx$ in the bulk.
	\item The DLR variational principle usually reads as an “if and only if” statement, and here we only prove one implication. The converse also holds: by Theorem~\ref{theo:SineBeta}, the unique (\cite{erbar2021one}) minimizer of $\Fb$ is the limit of $(\bXNz)_N$, hence a good solution of~\eqref{DLR} by Theorem~\ref{theo:Convergence}.
\end{enumerate}

\subsection{{Related works}}
\label{sec:literature}
In \cite[Sec 1.3]{dereudre2021dlr}, after proving the DLR equations for $\sineb$, we asked two questions:
\begin{itemize}
	\item Can one prove a variational principle, connecting DLR solutions and free energy minimizers?
	\item Are those solutions unique?
\end{itemize}
In~\cite{erbar2021one} we answered the uniqueness question positively \emph{for free energy minimizers}. The desire to get a “fully physical” proof of existence of $\sineb$ was also natural, and became even more so once we had these two “physical” descriptions of $\sineb$ in terms of DLR equations and as unique free energy minimizer. However, the very first step, namely the proof of local laws for the finite-volume model, yielding tightness of the observable, is not easy and was only recently performed in~\cite{peilen2022local}. 

Around the same time, in~\cite{leble2024dlr2d}, we derived DLR equations for all limit points of the \emph{two-dimensional} log-gas, and proved that they are all translation-invariant by adapting techniques from~\cite{frohlich1981absence,georgii1999translation} to our singular, long-range interaction. It seemed plausible that a similar, if not simpler, argument would apply to the one-dimensional case. The only missing piece would then be the DLR variational principle, for which there were no examples with long-range interactions in the literature.

\medskip 

Let us be slightly technical here. The proof of the DLR variational principle (for the “DLR $\implies$ free energy minimizer” implication, which is the easier one) in the classical setting (lattice spin systems with short-range interactions) relies on one estimate: in a large box $B_L := [-L, L]^d$ the interaction energy between the spins in $B_L$ and the spins outside $B_L$ is {negligible compared to the volume, i.e. $o(L^d)$, uniformly in the configuration}. This is a uniform, deterministic bound, but a quick inspection of the proof reveals that it works as well with a weaker, averaged, $o(L^d)$ estimate (see Proposition~\ref{prop:petite_entropie}).

\medskip

The very recent paper~\cite{assiotis2026uniqueness} answers the question of uniqueness of \emph{stationary} DLR solutions with finite renormalized energy, by a different route. Instead of connecting DLR solutions to free energy minimizers, it shows uniqueness directly. As in the classical case, the main estimate is an $o(L^d)$ bound on an interior-exterior interaction energy. “A cyclic averaging argument and a transport inequality”, which rely on stationarity, then give the result. This is not a variational principle, although combined with the fact that $\sineb$ minimizes the free energy~\cite{leble2017large}, it implies, indirectly, that stationary DLR solutions with finite energy are indeed free energy minimizers. 

The proof of the $o(L^d)$ estimate in~\cite{assiotis2026uniqueness}, although done in a slightly different setting, was inspirational for us, as it made clear that the path we were taking for our own $o(L^d)$ estimate was unnecessarily complicated. As a result, this part of our argument bears resemblance with \cite[Sec. 7]{assiotis2026uniqueness}.

\medskip 

Unlike~\cite{assiotis2026uniqueness}, we do not assume stationarity: we first prove that good DLR solutions are stationary (Theorem~\ref{theo:TI-intro}). Without this, we could not connect to the free energy and prove existence of $\sineb$. This step follows the classical method for proving translation-invariance (or absence of symmetry breaking), which was also implemented recently in~\cite{assiotis2026tail} to prove tail triviality of $\sineb$ for all $\beta > 0$. Tail triviality can also be deduced from our uniqueness result, see Remark~\ref{rem:tail}.

\subsection{{Conventions}}
\label{sec:conventions}
In the rest of the paper, $\La_R := [- \frac{R}{2}, \frac{R}{2}]$ denotes the interval of length $R$ centered at the origin.

\paragraph{Normalization of the energy.} We define all the energies with a factor $\hal$ in front, see~\eqref{def:HN} and~\eqref{def:FLa}. With this convention, $\beta = 2$ corresponds to the determinantal case. In~\cite{sandier2015log, leble2016logarithmic, leble2017large}, the energies are defined without the factor $\hal$, and the Gibbs measures carry a factor $\frac{\beta}{2}$ instead. 

In particular, the renormalized energy $\W$ and the intrinsic energy $\Wint$ of these papers are twice the ones used here, and the free energy reads $\frac{\beta}{2} \W + \Ee$ there. This only changes some constants, and we translate the results that we use into our normalization without further comment.

\paragraph{Intensity.} All the infinite-volume point processes considered in this paper have intensity one. 

\paragraph{Constants.} Unless stated otherwise, $\Cc$ and $\cc$ denote positive constants which depend only on $\beta$ (and on fixed test functions), and may change from line to line. A subscript, as in $\Cc_r$, indicates an additional dependence.

\subsection{Plan of the paper}
\begin{itemize}
	\item In Section~\ref{sec:electric_formalism_local_laws_and_consequences}, we recall the electric formalism and the local laws of~\cite{peilen2022local}, and we deduce controls on the number of points in the finite-volume model.
	\item In Section~\ref{sec:good}, we pass these controls and the finite-volume DLR equations to the limit, and prove Theorem~\ref{theo:Convergence}.
	\item In Section~\ref{sec:stationarity}, we prove that good DLR solutions are stationary (Theorem~\ref{theo:TI-intro}).
	\item In Section~\ref{sec:varprinc}, we prove the DLR variational principle (Theorem~\ref{theo:VarProp}).
	\item Appendix~\ref{sec:auxiliary} contains the more technical proofs, in particular those of Proposition~\ref{prop:finite-truncation}, item 2 of Theorem~\ref{theo:Convergence}, and Lemmas~\ref{lem:ZNbGC_asymptotics}, \ref{lem:box-electric} and~\ref{lem:four-terms}.
\end{itemize}

\subsection*{AI disclosure}
All the mathematical ideas were designed by us or found in the relevant literature. SC did not use any large language model. TL used OpenAI and Anthropic models for writing up standard arguments, proofreading, and editing.

\section{Electric formalism, local laws, and consequences}
\label{sec:electric_formalism_local_laws_and_consequences}
Recall that $\bXNz = \sum_{i=1}^N \delta_{\vrho x_i}$ is the local configuration rescaled to have intensity one near $0$. When $I \subset \R$ is an interval, we write $\bXNz(I) := \#\{i : \vrho x_i \in I\}$ for the number of points. The goal of this section is to obtain the following estimates concerning the finite-volume measure $\PNb$.

\begin{proposition}
\label{prop:FiniteVolEstimates}
There exist $\Cc, a, c_2 > 0$ such that the following holds.
\begin{enumerate}
	\item For all $N$ large enough and every interval $I \subset [- c_2 N, c_2 N]$ of length $R \geq 1$,
	\begin{equation*}
	\log \ENb\left[ e^{a \bXNz(I)^2} \right] \leq \Cc R^2.
	\end{equation*}
	\item For every interval $I$ of length $R \geq 1$,
	\begin{equation*}
	\limsup_{N \to \infty} \ENb\left[ \left( \bXNz(I) - R \right)^2 \right] \leq \Cc \log^2(2 + R).
	\end{equation*}
\end{enumerate}
\end{proposition}
\begin{itemize}
	\item Item 1 basically says that the number of points is of the same order as the length of the interval, which is not very precise. Its real strength is the fact that we can take a super-exponential moment.
	\item Item 2 is a form of \emph{hyperuniformity}: it says that the number variance is much smaller than the expected number of points ($\log^2 R$ versus $R$). This improves on a series of “discrepancy versus energy” comparisons, which ultimately gave a $o(R)$ bound, see~\cite{leble2017large, erbar2021one}. The big improvement comes from the fact that we are not only relying on finiteness of the energy, but on the whole Gibbsian structure. 

	In fact, our estimate is off by one power of $\log R$: the number variance in an interval of length $R$ is known, by “random matrix theory” methods, to be $\O(\log R)$ \cite{kritchevski2012scaling,najnudel2021uniform}. In this sense, our result is neither new nor optimal, but we wanted a “statistical physics” proof.
\end{itemize}

Proposition~\ref{prop:FiniteVolEstimates} follows, with some work, from the local laws of~\cite{peilen2022local}, which say that the one-dimensional log-gas is “well-behaved” down to the microscopic scale. In particular, we need the electric formalism of Sandier--Serfaty, which we briefly recall next, while referring to~\cite[Section 2]{peilen2022local} and to Serfaty's lecture notes~\cite{serfaty2024lectures} for the details.

\subsection{The electric formalism}
\label{sec:electric}

We write $\g(z) := - \log |z|$ for $z \in \R^2$. In $\R^2$ (not in $\R$!), it satisfies $- \Delta \g = 2\pi \delta_0$ in the sense of distributions.

\paragraph{Next-order logarithmic interaction energy.}
For $N \geq 2$ fixed, define the (rescaled) semicircle law as
\begin{equation*}
\mathbf{m}_N(\dd x) := \rhosc\left( \frac{x}{N} \right) \dd x.
\end{equation*}
This is a measure of total mass $N$, supported on $[-2N, 2N]$, whose density is bounded by $\vrho = \frac{1}{\pi}$. For $\X_N = (x_1, \dots, x_N)$, writing $\bX_N := \sum_{i=1}^N \delta_{x_i}$, define the \emph{next-order logarithmic interaction energy} as:
\begin{equation}
\label{def:FN}
F_N(\X_N) := \hal \iint_{x \neq y} \red{- \log |x - y|} \, \dd (\bX_N - \mathbf{m}_N)(x) \, \dd (\bX_N - \mathbf{m}_N)(y).
\end{equation}

\paragraph{Electric potential and electric field.}
The electric formalism for the one-dimensional log-gas requires an \emph{extension}: we identify $\R$ with the line $\R \times \{0\} \subset \R^2$, and, if $\mu$ is a measure on $\R$, we write $\mu \delta_\R$ for its image in $\R^2$. The \emph{electric potential} and the \emph{electric field} generated by $\X_N$ and the background $\mathbf{m}_N$ are defined as
\begin{equation*}
u := \g * \left( \bX_N - \mathbf{m}_N \delta_\R \right), \qquad \nabla u := \nabla \g * \left( \bX_N - \mathbf{m}_N \delta_\R \right).
\end{equation*}
They satisfy the Poisson equation
\begin{equation}
\label{eq:poisson}
- \Delta u = 2\pi \left( \bX_N - \mathbf{m}_N \delta_\R \right) \quad \text{in } \R^2.
\end{equation}

\paragraph{Nearest-neighbor truncation and smearing out charges.}
Since $\g$ is the Coulomb kernel in $\R^2$, one would like to write $F_N(\X_N) = \frac{1}{4\pi} \int_{\R^2} |\nabla u|^2$, but the field has non-square-integrable singularities near each charge. 

One thus \emph{truncates} the field: for $\eta > 0$, let $\f_\eta(z) := \left( \log \frac{\eta}{|z|} \right)_+$, and, for a family of radii $\vec{\mathsf{r}} = (\mathsf{r}_i)_{i=1, \dots, N}$, define the \emph{truncated potential} $u_{\vec{\mathsf{r}}} := u - \sum_{i=1}^N \f_{\mathsf{r}_i}(\cdot - x_i)$. Then $\nabla u_{\vec{\mathsf{r}}}$ is the \emph{truncated electric field}. Moreover, let $\sm{x}{\eta}$ be the uniform probability measure on the circle $\partial B(x, \eta)$ (a “smeared charge”). By Newton's theorem, $u_{\vec{\mathsf{r}}} = \g * \left( \sum_i \sm{x_i}{\mathsf{r}_i} - \mathbf{m}_N \delta_\R \right)$, and thus:
\begin{equation}
\label{eq:poisson-trunc}
- \Delta u_{\vec{\mathsf{r}}} = 2\pi \left( \sum_{i=1}^N \sm{x_i}{\mathsf{r}_i} - \mathbf{m}_N \delta_\R \right) \quad \text{in } \R^2.
\end{equation}
In practice, we choose the radii $\mathsf{r}_i := \frac14 \min\left(1, \min_{j \neq i} |x_i - x_j|\right)$, so that the smeared charges do not overlap. One then has the exact identity~\cite[Section 2]{peilen2022local}
\begin{equation}
\label{eq:energy-identity}
F_N(\X_N) = \frac{1}{4\pi} \int_{\R^2} |\nabla u_{\vec{\mathsf{r}}}|^2 + \hal \sum_{i=1}^N \log \mathsf{r}_i - \sum_{i=1}^N \int \f_{\mathsf{r}_i}(x - x_i) \dd \mathbf{m}_N(x).
\end{equation}
In particular, it is classical that the energy $F_N$ is bounded below by $- \Cc N$, and that it controls the $L^2$ norm of the truncated field, see~\eqref{eq:coercivity} below and~\cite[Lemma A.4]{peilen2022local}. Because~\eqref{eq:energy-identity} expresses the energy $F_N$ in terms of electric fields, we often speak of “electric energy”.

\paragraph{Local energies.}
For an interval $\Omega \subset \R$, the \emph{local energy} $F_N^\Omega(\X_N)$ is the part of the right-hand side of~\eqref{eq:energy-identity} which “lives above $\Omega$”, computed with the truncation radii adapted to $\Omega$ of~\cite[(2.14)]{peilen2022local}:
\begin{equation*}
\tr_i = \tr_i(\Omega) := \frac14 \min\left( 1, \min_{x_j \in \Omega, \, j \neq i} |x_i - x_j| \right) \quad \text{if } \mathrm{dist}(x_i, \partial \Omega) \geq \frac12, \qquad \tr_i := \frac14 \quad \text{otherwise}.
\end{equation*}
In particular, $\tr_i \leq \frac14$ for every $i$. Here, we will simply use the following fact (\cite[Lemma A.4]{peilen2022local}):
\begin{lemma}
There exists $\Cc_0$ such that, for $\Omega \subset \R$ of length $\ell$, letting $\widetilde{\Omega} := \Omega \times [-\ell, \ell]$ and $\tr = \tr(\Omega)$,
\begin{equation}
\label{eq:coercivity}
\int_{\widetilde{\Omega}} |\nabla u_{\tr}|^2 \leq 8 \pi \left( F_N^\Omega(\X_N) + \Cc_0 \, \#(\X_N \cap \Omega) \right).
\end{equation}
The same holds for $\Omega = \R$, $\widetilde{\Omega} = \R^2$.
\end{lemma} 

\paragraph{Re-writing the Gibbs measure.}
Expanding $\bX_N = \mathbf{m}_N + (\bX_N - \mathbf{m}_N)$ in~\eqref{def:HN} and using the Euler--Lagrange equation of the semicircle law, one obtains the \emph{splitting formula} of Sandier and Serfaty~\cite{sandier2015log}: 
\begin{equation*}
\HN(\X_N) = F_N(\X_N) + \sum_{i=1}^N \zeta_N(x_i) + \mathrm{const}_N,
\end{equation*}
where the \emph{effective confinement} $\zeta_N \geq 0$ vanishes on $[-2N, 2N]$ and $\mathrm{const}_N$ does not depend on $\X_N$.

In particular, one can re-write the Gibbs measure introduced in \eqref{def:PNb} as:
\begin{equation}
\label{eq:PNb-electric}
\dd \PNb(\X_N) = \frac{1}{\K_{N,\beta}} \exp\left( - \beta \left( F_N(\X_N) + \sum_{i=1}^N \zeta_N(x_i) \right) \right) \frac{\dd \X_N}{N^N},
\end{equation}
for some (next-order) partition function $\K_{N,\beta}$ (normalized here as in~\cite[(2.2)]{peilen2022local}). We will only use the form of this density, in Section~\ref{sec:good}, to rewrite the finite-volume measure in the local coordinates.

\subsection{Finite-volume input: the local laws}
\label{sec:Peilen}
The following results are proven in~\cite{peilen2022local}, in greater generality, allowing for more general confining potentials. We state them in our coordinates. All the constants below depend on $\beta$, but not on $N$.

The (blown-up) \emph{bulk} is defined as in \cite[(1.18)]{peilen2022local} by $\Bulk_N := [-N, N]$. A quick inspection of the proof shows that anything of the form $[-(2-\epsilon)N, (2-\epsilon)N]$ would work as well if needed. On $\Bulk_N$, the density of $\mathbf{m}_N$ is bounded below by $\rhosc(1) = \frac{\sqrt{3}}{2\pi}$ and above by $\vrho$.

\begin{theorem}[Local law, {\cite[Theorem 1]{peilen2022local}}]
\label{theo:local-law}
There exist constants $\omega \geq 1$ and $\Cc, \cc > 0$ such that the following holds for $N$ large enough. For every interval $\Omega \subset \Bulk_N$ of length $\ell \geq \omega$, there exists an event $\G_\Omega$ such that
\begin{equation*}
\PNb\left( \G_\Omega^c \right) \leq \Cc e^{- \cc \ell} \qquad \text{and} \qquad F_N^\Omega(\X_N) + \Cc_0 \, \#(\X_N \cap \Omega) \leq \Cc \, \ell \quad \text{on } \G_\Omega.
\end{equation*}
In particular, by~\eqref{eq:coercivity}, on the event $\G_\Omega$ we have
\begin{equation}
\label{eq:local-law-field}
\int_{\widetilde{\Omega}} |\nabla u_{\tr}|^2 \leq 8\pi \Cc \, \ell.
\end{equation}
\end{theorem}
The main message is that, at scale $\ell$, the local energy is typically of order $\ell$.

\paragraph{Fluctuations.} For a test function $\xi$, we define the \emph{fluctuation}
\begin{equation*}
\Fluct_N(\xi) := \int \xi \, \dd \left( \bX_N - \mathbf{m}_N \right).
\end{equation*}
In Appendix~\ref{sec:how_to_read_theorem_ref_theo_bounded_fluct}, we explain how to read the following result (with a different scaling) from~\cite{peilen2022local}.
\begin{theorem}[Bounded fluctuations, {\cite[Theorem 2]{peilen2022local}}]
\label{theo:bounded-fluct}
Let $\theta \in C^3_c(\R)$ be supported in $[-4, 4]$. There exist constants $\Cc, \cc, c_1 > 0$, depending on $\beta$ and $\theta$, such that, for $N$ large enough, for every $z \in \R$ and $\ell \in [\omega, c_1 N]$ such that the function $\xi := \theta\left( \frac{\cdot - z}{\ell} \right)$ is supported in $\Bulk_N$, there exists an event $\G$ such that:
\begin{equation}
\label{eq:bounded-fluct}
\PNb(\G^c) \leq \Cc e^{- \cc \ell}, \qquad \log \ENb\left[ \exp\left( s \Fluct_N(\xi) \right) \ind_{\G} \right] \leq \Cc \quad \text{for all } |s| \leq 1.
\end{equation}
\end{theorem}

\paragraph{Comments.}
\begin{itemize}
	\item Theorem~\ref{theo:local-law} gives an exponential bound on the \emph{local electric energy}, at every scale larger than a fixed microscopic scale $\omega$, \emph{up to a rare bad event} whose probability is exponentially small in the scale.
	\item Theorem~\ref{theo:bounded-fluct} is much finer: it says that the fluctuations of smooth linear statistics are \emph{bounded}, with exponential moments, \emph{uniformly in the scale}. This is important for us.
	\item Similar controls on the eigenvalues of $\beta$-ensembles have been obtained by random matrix methods (rigidity estimates, loop equations, Dyson Brownian motion), see e.g.~\cite{bourgade2014universality, erdos2017dynamical} and the optimal local law of~\cite{bourgade2022optimal}.
\end{itemize}

\subsection{Exponential moments of the number of points}
\label{sec:number-points}
In this section, we prove item 1 of Proposition~\ref{prop:FiniteVolEstimates}: there exist $a, \Cc > 0$ such that, for $N$ large enough and every interval $I \subset [-\frac{N}{2}, \frac{N}{2}]$ of length $R \geq 1$,
\begin{equation}
\label{eq:superexp}
\log \ENb\left[ e^{a \bXN(I)^2} \right] \leq \Cc R^2.
\end{equation}

\subsubsection*{Controlling the number of points by the energy}
The following lemma expresses the fact that if the number of points is much larger than the length of the interval, this is felt by the energy.
\begin{lemma}
\label{lem:overcrowding-energy}
There exists $\Cc$ such that the following holds. Let $R > 2$ and let  $\Omega$ be an interval. Let $z$ be a point in $\Omega$ and let $I = [z-R/2, z+R/2]$, and $I' := [z - R, z + R]$. Assume that $I' \subset \Omega$. Then
\begin{equation}
\label{eq:overcrowding-test}
\bXN(I)^2 \leq  \Cc \left(R^2 + \int_{\widetilde{\Omega}} |\nabla u_{\tr}|^2\right),
\end{equation}
where $\widetilde{\Omega}$ and $\tr = \tr(\Omega)$ are as in~\eqref{eq:coercivity}.
\end{lemma}
\begin{proof}
Let $\chi : \R \to [0,1]$ be a smooth function, equal to $1$ on $[-\frac34, \frac34]$ and supported in $[-1, 1]$, and let
\begin{equation*}
\Phi(t, s) := \chi\left( \frac{t - z}{R} \right) \chi\left( \frac{s}{R} \right).
\end{equation*}
By scaling, $\|\nabla \Phi\|_{L^2(\R^2)} \leq \Cc$. The support of $\Phi$ is contained in $I' \times [-R, R] \subset \widetilde{\Omega}$. Since the radii $\tr_i$ are at most $\frac14$, the function $\Phi$ is equal to $1$ on each smeared charge centered in $I$, and $\Phi \geq 0$. Testing~\eqref{eq:poisson-trunc}, with the radii $\tr$, against $\Phi$, we obtain
\begin{equation*}
\bXN(I) \leq \int \Phi(t, 0) \, \dd \mathbf{m}_N(t) + \frac{1}{2\pi} \left| \int_{\R^2} \nabla \Phi \cdot \nabla u_{\tr} \right| \leq 2 \vrho R + \Cc \left( \int_{\widetilde{\Omega}} |\nabla u_{\tr}|^2 \right)^{1/2},
\end{equation*}
since the density of $\mathbf{m}_N$ is bounded by $\vrho$. Squaring gives~\eqref{eq:overcrowding-test}.
\end{proof}
This is a rough bound (one may want to control $(\bXN(I) -R)^2$ instead), but it is enough for us.

\subsubsection*{Strategy: controlling the deviations}
By Theorem~\ref{theo:local-law}, the right-hand side of~\eqref{eq:overcrowding-test} is controlled at every scale $\ell \geq \omega$, but \emph{only outside a bad event of probability at most $\Cc e^{- \cc \ell}$}. To get rid of this, instead of controlling the exponential moments directly, we first look for a control on the \emph{deviations} of $\bXN(I)$. Very bad deviations at a given scale, whose real unlikeliness is obfuscated by this fixed bad event of probability $e^{- \cc \ell}$, are also “felt” at a bigger scale, for which the “bad event” is more rare. 

\begin{claim}[Deviations of the number of points]
\label{claim:deviations}
There exist $\Cc, \cc > 0$ such that, for $N$ large enough, for every interval $I$ of length $R > 2$ {whose center lies in $[-\frac{N}{2}, \frac{N}{2}]$}, and for every $M \geq \Cc R$,
\begin{equation}
\label{eq:overcrowding}
\PNb\left( \bXN(I) \geq M \right) \leq \Cc e^{- \cc M^2}.
\end{equation}
\end{claim}

\begin{proof}
Let $I$ be such an interval, with center $z$, and assume that $M \geq \Cc R$, where $\Cc \geq 8$ is a large constant to be chosen later. If $R > \frac{N}{8}$, the event $\{\bXN(I) \geq M\}$ is empty, since there are only $N$ points, so we may assume that $R \leq \frac{N}{8}$. We distinguish two cases:

\paragraph{Case 1: $M^2 > KN$.} For $K$ large enough, this is covered by the large deviation estimate of~\cite[(1.38)]{sandier2015log}.

\paragraph{Case 2: $M^2 \leq KN$.} Let $\delta$ be small, to be chosen later, and let $\Omega$ be the interval of center $z$ and length
\begin{equation*}
\ell := \min\left( \delta M^2, N \right) \geq \min\left( \delta, \tfrac{1}{K} \right) M^2.
\end{equation*}
This is the larger scale at which the excess of points will be felt. Let $\Cc_1$ be the constant of Lemma~\ref{lem:overcrowding-energy}, and $\Cc_2 := 8 \pi \Cc$, with $\Cc$ as in~\eqref{eq:local-law-field}. If $\Cc$ is large enough, we have $\ell \geq 2R$ and $\ell \geq \omega$, so that $I' := [z - R, z + R] \subset \Omega$, and $\Omega \subset \Bulk_N$ since $|z| \leq \frac{N}{2}$ and $\ell \leq N$. On the event $\G_\Omega$ of Theorem~\ref{theo:local-law}, the electric energy in $\widetilde{\Omega}$ is at most $\Cc_2 \ell \leq \Cc_2 \delta M^2$. By Lemma~\ref{lem:overcrowding-energy}, on $\G_\Omega$,
\begin{equation*}
\bXN(I)^2 \leq \Cc_1 \left( R^2 + \Cc_2 \delta M^2 \right) \leq \Cc_1 \left( \frac{1}{\Cc^2} + \Cc_2 \delta \right) M^2 < M^2,
\end{equation*}
if $\delta$ is chosen small enough and $\Cc$ large enough. Hence the event $\{\bXN(I) \geq M\} \subset \G_\Omega^c$, and
\begin{equation*}
\PNb\left( \bXN(I) \geq M \right) \leq \PNb\left( \G_\Omega^c \right) \leq \Cc e^{- \cc \ell} \leq \Cc e^{- \cc \min(\delta, \frac{1}{K}) M^2}.
\end{equation*}
\end{proof}

\subsubsection*{Conclusion: proof of~\eqref{eq:superexp}}
Let $a := \frac{\cc}{2}$, with $\Cc, \cc$ as in Claim~\ref{claim:deviations}. Splitting according to the value of $\bXN(I)$, we get:
\begin{equation*}
\ENb\left[ e^{a \bXN(I)^2} \right] \leq e^{a \Cc^2 R^2} + \sum_{M \geq \Cc R} e^{a (M+1)^2} \, \PNb\left( \bXN(I) \geq M \right) \leq e^{a \Cc^2 R^2} + \Cc \sum_{M \geq 1} e^{a(M+1)^2 - \cc M^2},
\end{equation*}
and the last sum is finite. This proves~\eqref{eq:superexp}. \qed

\subsection{Number variance}
\label{sec:number-variance}
We now prove item 2 of Proposition~\ref{prop:FiniteVolEstimates}: for $N$ large enough and every interval $I \subset [-\frac{N}{2}, \frac{N}{2}]$ of length $R \in [1, c_1 N]$,
\begin{equation}
\label{eq:variance}
\ENb\left[ \left( \bXN(I) - \mathbf{m}_N(I) \right)^2 \right] \leq \Cc \log^2 (2 + R).
\end{equation}
The number variance is the second moment of the discrepancy $\bXN(I) - \mathbf{m}_N(I) = \Fluct_N(\ind_I)$, which is a fluctuation, but with a non-continuous test function. The local laws of Theorem~\ref{theo:bounded-fluct} deal with fluctuations of \emph{smooth} functions, and only outside a bad event. We thus do two things: we remove the bad event, using the exponential moments of the number of points, and we decompose $\ind_I$ into a sum of smooth functions living at all scales between $1$ and $R$, plus a small remainder.

\subsubsection*{Removing the bad event}
\begin{lemma}[Bounded fluctuations without “bad event”]
\label{lem:smooth-fluct}
Let $\Theta$ be a finite family of functions in $C^3_c(\R)$ supported in $[-4,4]$, and let $c_1 \leq \frac{1}{8}$ be such that Theorem~\ref{theo:bounded-fluct} holds for each $\theta \in \Theta$. There exist $s_0 > 0$ and $\Cc$ such that, for $N$ large enough, for every $\theta \in \Theta$, $z \in [-\frac{N}{2}, \frac{N}{2}]$ and $\ell \in [\omega, c_1 N]$, letting $\xi := \theta(\frac{\cdot - z}{\ell})$,
\begin{equation}
\label{eq:smooth-fluct}
\log \ENb\left[ e^{s \Fluct_N(\xi)} \right] \leq \Cc \quad \text{for all } |s| \leq s_0, \qquad \text{and in particular} \qquad \ENb\left[ \Fluct_N(\xi)^2 \right] \leq \Cc.
\end{equation}
\end{lemma}
\begin{proof}
The support of $\xi$ is contained in $I_\xi := [z - 4\ell, z + 4\ell] \subset \Bulk_N$, so Theorem~\ref{theo:bounded-fluct} applies, with a certain “good” event $\G$. On the complement of $\G$, we use the crude bound
\begin{equation*}
|\Fluct_N(\xi)| \leq \|\theta\|_\infty \left( \bXN(I_\xi) + \mathbf{m}_N(I_\xi) \right) \leq \Cc \left( \bXN(I_\xi) + \ell \right).
\end{equation*}
By the Gaussian tail~\eqref{eq:overcrowding} applied to $I_\xi$, which has length $8 \ell$, we have, for $t \in [0, 1]$,
\begin{equation*}
\ENb\left[ e^{t \bXN(I_\xi)} \right] \leq e^{8 \Cc t \ell} + \Cc \sum_{M \geq 1} e^{t (M+1) - \cc M^2} \leq e^{8 \Cc t \ell} + \Cc'.
\end{equation*}
By the Cauchy--Schwarz inequality, applying this with $t = 2 \Cc |s|$, we get, for $|s| \leq s_0$,
\begin{equation*}
\ENb\left[ e^{s \Fluct_N(\xi)} \ind_{\G^c} \right] \leq \ENb\left[ e^{2 \Cc |s| (\bXN(I_\xi) + \ell)} \right]^{\hal} \PNb(\G^c)^{\hal} \leq \Cc \, e^{\Cc' s_0 \ell} \, e^{- \frac{\cc}{2} \ell} \leq \Cc,
\end{equation*}
if $s_0$ is chosen small enough (depending on $\cc$ and on the family $\Theta$). Combined with~\eqref{eq:bounded-fluct}, this gives the first bound in~\eqref{eq:smooth-fluct}. The second one follows from $u^2 \leq s_0^{-2} (e^{s_0 u} + e^{-s_0 u})$.
\end{proof}

\subsubsection*{A dyadic decomposition}
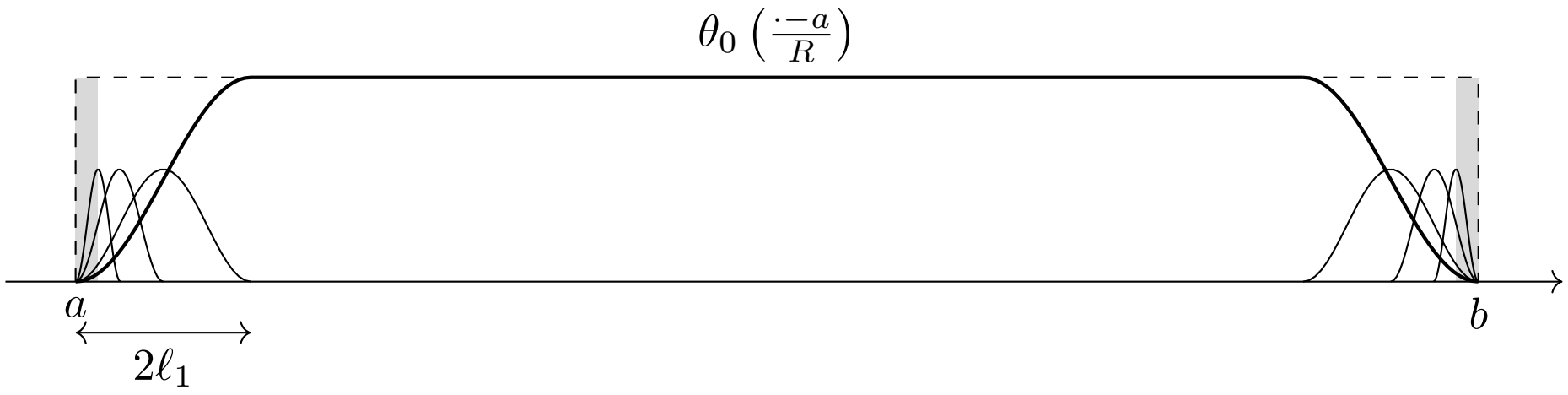
\begin{figure}[ht]
\centering
\begin{tikzpicture}[xscale=1.1, yscale=1.6]
\fill[gray!30] (0,0) rectangle (0.16,1);
\fill[gray!30] (9.84,0) rectangle (10,1);
\draw[->] (-0.5,0) -- (10.6,0);
\draw[dashed] (0,0) -- (0,1) -- (10,1) -- (10,0);
\draw[thick, domain=0:1.25, samples=30] plot (\x, {sin(deg(pi*\x/2.5))^2});
\draw[thick] (1.25,1) -- (8.75,1);
\draw[thick, domain=8.75:10, samples=30] plot (\x, {sin(deg(pi*(10-\x)/2.5))^2});
\foreach \l in {0.625, 0.3125, 0.16} {
  \draw[domain=0:2*\l, samples=30] plot (\x, {0.55*sin(deg(pi*\x/(2*\l)))^2});
  \draw[domain=10-2*\l:10, samples=30] plot (\x, {0.55*sin(deg(pi*(10-\x)/(2*\l)))^2});
}
\node[below] at (0,0) {$a$};
\node[below] at (10,0) {$b$};
\node[above] at (5,1) {$\theta_0\left(\frac{\cdot - a}{R}\right)$};
\draw[<->] (0,-0.25) -- (1.25,-0.25) node[midway, below] {$2 \ell_1$};
\end{tikzpicture}
\caption{The dyadic decomposition~\eqref{eq:dyadic}.}
\label{fig:dyadic}
\end{figure}

Let $I = [a, b] \subset [-\frac{N}{2}, \frac{N}{2}]$, with $R = b - a \leq c_1 N$. If $R \leq 16 \omega$, the bound~\eqref{eq:variance} follows from~\eqref{eq:superexp}, so we assume that $R > 16 \omega$. The indicator function of $I$ is not smooth, and we cannot apply Lemma~\ref{lem:smooth-fluct} to it directly. Instead, we use a dyadic partition of unity of $I$, adapted to its endpoints, see Figure~\ref{fig:dyadic}. Let $\ell_j := 2^{-j-3} R$ for $j \geq 0$, and let $J \geq 1$ be such that $\ell_J \in [\omega, 2\omega)$. There exist three fixed smooth functions $\theta_0, \theta, \check{\theta}$, supported in $[-4, 4]$ and independent of $N$, $I$ and $j$, such that
\begin{equation}
\label{eq:dyadic}
\ind_I = \theta_0\left( \frac{\cdot - a}{R} \right) + \sum_{j=1}^{J} \left( \theta\left( \frac{\cdot - a}{\ell_j} \right) + \check{\theta}\left( \frac{\cdot - b}{\ell_j} \right) \right) + r,
\end{equation}
where the remainder $r$ takes values in $[0,1]$ and is supported in $[a, a + 2 \ell_J] \cup [b - 2 \ell_J, b]$.

\subsubsection*{Conclusion}
By Lemma~\ref{lem:smooth-fluct}, applied to the family $\Theta := \{\theta_0, \theta, \check{\theta}\}$ at the scales $R, \ell_1, \dots, \ell_J \in [\omega, c_1 N]$, each of the $2J + 1$ smooth pieces of~\eqref{eq:dyadic} has a fluctuation with second moment bounded by $\Cc$. Since $J \leq \log_2 R$, Minkowski's inequality in $L^2(\PNb)$ gives
\begin{equation*}
\ENb\left[ \Fluct_N(\ind_I - r)^2 \right]^{\hal} \leq \Cc (2J + 1) \leq \Cc \log (2 + R).
\end{equation*}
On the other hand, we have the very rough bound
\begin{equation*}
|\Fluct_N(r)| \leq \bXN([a, a + 4 \omega]) + \bXN([b - 4 \omega, b]) + 8 \vrho \omega,
\end{equation*}
and the second moment of the right-hand side is bounded thanks to~\eqref{eq:superexp}. This proves~\eqref{eq:variance}. 
\qed

\subsubsection*{Back to the rescaled configuration: proof of Proposition~\ref{prop:FiniteVolEstimates}}
Passing from a bound on the number of points of $\bXN$ to a bound for $\bXNz$ is easy. To express the number variance, we need to replace the varying measure $\mathbf{m}_N$ appearing in \eqref{eq:variance} by the constant background. This is not difficult: let $I$ be {a fixed interval of length $R \geq 1$}, and let $\tilde{I} := {\vrho^{-1}} I${, so that $\bXNz(I) = \bXN(\tilde{I})$ and $|\tilde{I}| = \pi R$. For $N$ large enough, $\tilde{I} \subset [-\frac{N}{2}, \frac{N}{2}]$ and $\pi R \leq c_1 N$}. Since $\rhosc$ is continuous with $\rhosc(0) = \vrho$, we have
\begin{equation*}
\mathbf{m}_N({\tilde{I}}) = \int_{{\tilde{I}}} \rhosc\left( \frac{x}{N} \right) \dd x \longrightarrow \vrho \, |{\tilde{I}}| = R \qquad \text{as } N \to \infty.
\end{equation*}
Hence, by~\eqref{eq:variance},
\begin{equation*}
\limsup_{N \to \infty} \ENb\left[ \left( \bXNz(I) - R \right)^2 \right] \leq 2 \limsup_{N \to \infty} \ENb\left[ \left( \bXN({\tilde{I}}) - \mathbf{m}_N({\tilde{I}}) \right)^2 \right] \leq 2 \Cc \log^2(2 + \pi R),
\end{equation*}
which gives item 2 of Proposition~\ref{prop:FiniteVolEstimates}.

\begin{remark}
The loss of one power of $\log R$ compared to the “true” bound comes from applying a rough triangle inequality in $L^2$, which ignores the fact that fluctuations at different scales should be almost independent. Presumably this could be improved.
\end{remark}

\section{Passing good properties to the limit}
\label{sec:good}
In this section, we prove items 1 and 2 of Theorem~\ref{theo:Convergence}, by passing to the limit $N \to \infty$ some properties of the finite-volume measure: 
\begin{itemize}
 	\item The control on the number of points stated in Proposition~\ref{prop:FiniteVolEstimates}.
 	\item The finite-volume DLR equations, which are straightforward.
 \end{itemize} 
The proof of the DLR equations follows the strategy of~\cite{dereudre2021dlr}, but our finite-volume measure is not stationary. In particular, the system has edges, which we must control. {In the circular setting, stationarity provides averaging properties, which are not available here. Instead,} we rely on the local laws from~\cite{peilen2022local}. 

There is an extra difficulty: passing DLR equations to the limit requires a notion of convergence, namely local convergence of point processes, that would allow test functions which are not necessarily continuous. However, the compactness that we get “for free” from local laws holds only for the topology of \emph{weak} convergence. What is missing is either:
\begin{itemize}
	\item A technical input: knowing that $k$-point correlation functions are bounded can be enough to “regularize” any test function into a continuous one and upgrade “weak” into “local” convergence. This was done in~\cite{leble2024dlr2d} for the two-dimensional log-gas, relying on results from~\cite{thoma2023overcrowding} whose proof seems quite specific to the two-dimensional situation.
	\item A relative entropy bound: sub-level sets of the specific entropy are known to be compact for the local topology. When working with the \emph{stationary} version of the finite-volume Gibbs measure (here, it would mean working with the Circular $\beta$-ensemble instead of the Gaussian $\beta$-ensemble, or equivalently with a periodic log-gas), this bound can be obtained without proving local laws. This was the approach chosen in~\cite{dereudre2021dlr} and in related papers~\cite{dereudre2023rigidity}. 
\end{itemize}
For the \emph{one-dimensional} log-gas \emph{without} periodic boundary condition, we obtain a relative entropy bound “by hand”, using finite-volume DLR equations and local laws, see Section~\ref{sec:tightness}.

\subsection{Various topologies for convergence of point processes}
Recall that the \emph{weak topology} on point processes is the coarsest topology for which the maps $\Pp \mapsto \E_\Pp[f]$ are continuous for every bounded, continuous $f$. We give the definition of \emph{local} topology, and introduce a stronger version called \emph{exp-local} topology.
\begin{definition}[Local and exp-local topologies]
\label{def:exp-local}
Let $a > 0$ be as in Proposition~\ref{prop:FiniteVolEstimates}, and $\mathfrak{a} := \frac{a}{4}$. 
A measurable function $f : \Conf \to \R$ is said to be:
\begin{itemize}
	\item \emph{local} if there is a bounded interval $I$ such that $f(\bX) = f(\bX_I)$ for all $\bX$,
	\item \emph{exp-tame} if there are a bounded interval $I$ and a constant $\Cc_f$ such that 
	\begin{equation*}
|f(\bX)| \leq \Cc_f \exp\left( \mathfrak{a} \, \bX(I)^2 \right) \text{ for all $\bX$.}
	\end{equation*} 
\end{itemize}
Then:
\begin{itemize}
	\item The \emph{local topology}  is the coarsest topology for which the maps $\Pp \mapsto \E_\Pp[f]$ are continuous for every local, bounded, measurable $f$. 
	\item The \emph{exp-local topology} is the coarsest topology for which they are continuous for every local, exp-tame, measurable $f$.
\end{itemize}
\end{definition}
The exp-local topology allows us to pass more quantities to the limit, in particular exponential moments of the number of points, see Theorem~\ref{theo:tightness}.

\subsection{DLR equations and move functions}
\label{sec:DLR-def}
We start with some notation:
\begin{itemize}
	\item For a Borel set $A \subset \R$ and $\bX \in \Conf$, we write $\bX_A$ for the restriction of $\bX$ to $A$ and $\bX(A)$ for its number of points in $A$.
	\item  For a bounded Borel set $\La$, we write $\La^c := \R \setminus \La$, and $\bX = \bX_\La + \bX_{\La^c}$.
	\item  For $n \geq 0$, let $\Conf(\La, n)$ be the set of configurations with exactly $n$ points, all in $\La$, and $\Conf(\La) := \bigcup_{n \geq 0} \Conf(\La, n)$.
	\item We let $\Bin_{\La, n}$ be the binomial point process with $n$ points in $\La$, i.e. the law of $n$ independent points uniformly distributed in $\La$.
	\item We let  $\Poi_\La$ be the Poisson point process of intensity one in $\La$. 
	\item  For a point process $\Pp$, we write $\Pp_\La$ for the law of $\bX_\La$ under $\Pp$.
\end{itemize}
Finally, for $\rho : \R \to [0, + \infty)$ and $\bX' \in \Conf(\La, n)$, we let
\begin{equation*}
F^\rho_\La(\bX') := \hal \iint_{(x,y) \in \La^2, \, x \neq y} - \log|x - y| \, \dd (\bX' - \rho)(x) \, \dd (\bX' - \rho)(y),
\end{equation*}
so that $F^1_\La = F_\La$ is the energy defined in~\eqref{def:FLa}.

\subsubsection*{Move functions, good processes}
\begin{definition}[\cite{dereudre2021dlr}]
\label{def:move}
Let $\La$ be a bounded interval, let $\bX \in \Conf$ and $\bX' \in \Conf(\La)$ such that $\bX'$ and $\bX$ have the same number of points in $\La$. For $T > 0$, define the \emph{truncated move function} $\M^T$ as:
\begin{equation*}
\M^T_\La(\bX', \bX) := \iint_{\La \times ([-T, T] \setminus \La)} - \log |x - y| \, \dd (\bX' - \bX_\La)(x) \, \dd (\bX - \Leb)(y).
\end{equation*}
The \emph{move function} $\M_\La(\bX', \bX)$ is defined as the limit of $\M^T_\La(\bX', \bX)$ as $T \to \infty$, when it exists.
\end{definition}
The existence of this limit was proven in~\cite[Lemma 2.15 and Proposition 2.16]{dereudre2021dlr} for every \emph{stationary} point process with finite expected renormalized energy. This assumption was used to bound the number of points and the discrepancies in large intervals.

Here, we state the same result for a class of “good” processes, which are not assumed to be stationary, or to have finite energy, but where discrepancies are controlled by definition. We will apply it to the limit points of $\bXNz$ in order to derive the DLR equations, and we will only deduce stationarity later, \emph{from} the DLR equations, in Section~\ref{sec:stationarity}. 

\begin{definition}[Good point processes]
\label{def:good}
We say that a point process $\Pp$ on $\R$ is \emph{good} if there is a constant $\Cc_\Pp$ such that, for every interval $I$ of length $|I| \geq 1$,
\begin{equation*}
\E_\Pp\left[ \bX(I)^2 \right] \leq \Cc_\Pp |I|^2, \qquad \E_\Pp\left[ \left( \bX(I) - |I| \right)^2 \right] \leq \Cc_\Pp \log^2 (2 + |I|).
\end{equation*}
\end{definition}
The first bound follows from the second one, since $|I| \geq 1$. We state it for convenience.
\begin{remark}
\label{rem:good-relax}
The logarithmic bound on the number variance is what we obtain for the limit points (Theorem~\ref{theo:tightness}), but an inspection of the proofs suggests that our arguments only need some quantitative control on the fluctuations. For instance, a bound of the form $\E_\Pp[(\bX(I) - |I|)^2] \leq \Cc_\Pp |I|^{\kappa}$ with $\kappa < \frac12$ should be enough. The constraint comes from the term $\mathcal{T}_2$ in Section~\ref{sec:petite-entropie}, which would then be of order $N^{\frac34 + \frac{\kappa}{2}}$ instead of $N^{\frac34} \log^2 N$. We do not pursue this here.
\end{remark}

With this definition, we have the following existence result.
\begin{proposition}[Move functions, truncation]
\label{prop:move}
Let $\Pp$ be a good point process, and let $\La$ be a bounded interval. For $\Pp$-almost every $\bX$, $\M_\La(\bX', \bX)$ exists for every $\bX' \in \Conf(\La)$ such that $\bX'$ and $\bX$ have the same number of points in $\La$. Moreover, for all $\epsilon > 0$, we have
\begin{equation}
\label{eq:move-tail}
\lim_{T \to \infty} \Pp\left( \sup_{\bX' \in \Conf(\La, \bX(\La))} \left| \M_\La(\bX', \bX) - \M^T_\La(\bX', \bX) \right| > \epsilon \right) = 0.
\end{equation}
\end{proposition}
The proof of~\cite{dereudre2021dlr} relies on an integration by parts against the discrepancies $\bX((0, y]) - y$, that are controlled, thanks to stationarity, by the energy. Here, they are controlled directly, and the same line of computations applies. For completeness, we give a proof in Appendix~\ref{sec:existence_of_move_functions}.

\subsubsection*{Canonical DLR equations}
\begin{definition}[Canonical DLR equations]
\label{def:DLR}
We say that a point process $\Pp$ on $\R$ \emph{satisfies the canonical DLR equations}, or that it is a \emph{solution of~\eqref{DLR}}, if for every bounded interval $\La$:
\begin{enumerate}
	\item for $\Pp$-almost every $\bX$, the move functions $\M_\La(\bX', \bX)$ exist for all $\bX' \in \Conf(\La, \bX(\La))$,
	\item for every bounded measurable function $f : \Conf \to \R$, we have $\E_\Pp\left[ f \right] = \E_\Pp\left[ f_\La \right]$, where $f_\La$ is defined as:
	\begin{equation}
	\label{eq:DLR}
	f_\La(\bX) := \frac{1}{\K_\La(\bX)} \int_{\Conf(\La, \bX(\La))} f(\bX' + \bX_{\La^c}) \, e^{- \beta \left( F_\La(\bX') + \M_\La(\bX', \bX) \right)} \dd \Bin_{\La, \bX(\La)}(\bX'),
	\end{equation}
	with $\K_\La(\bX) := \int e^{- \beta \left( F_\La(\bX') + \M_\La(\bX', \bX) \right)} \dd \Bin_{\La, \bX(\La)}(\bX')$.
\end{enumerate}
\end{definition}
\paragraph{Remarks.}
\begin{itemize}
	\item The normalization $\K_\La(\bX)$ is positive and finite, because $F_\La$ and $\M_\La(\cdot, \bX)$ are bounded below on $\Conf(\La, n)$. For the move function, this follows from the uniform convergence of Proposition~\ref{prop:move}, since $\M^T_\La(\cdot, \bX)$ is bounded below.
	\item The move function $\M_\La(\bX', \bX)$ is the difference of a term depending on $(\bX', \bX_{\La^c})$ and of a term depending on $(\bX_\La, \bX_{\La^c})$. The second one cancels between the numerator and the denominator of $f_\La$, so that \emph{$f_\La(\bX)$ only depends on $\bX(\La)$ and on $\bX_{\La^c}$}. In other terms, $f_\La$ is measurable with respect to the $\sigma$-algebra $\mathcal{G}_\La$ generated by $\bX(\La)$ and $\bX_{\La^c}$, and~\eqref{eq:DLR} says that, conditionally on $\mathcal{G}_\La$, the law of $\bX_\La$ is the Gibbs measure~\eqref{DLR}. The kernels are \emph{proper}: they do not change $\mathcal{G}_\La$-measurable quantities.
	\item Resampling only changes a configuration in a bounded region, while the existence of the move functions is an asymptotic property. Thus the kernels preserve the set of configurations for which the move functions exist for every bounded interval. They are also \emph{consistent}: for bounded intervals $\La \subset \La'$, we have $(f_\La)_{\La'} = f_{\La'}$. This follows from the additivity of the energies, which is clear at the truncated level and passes to the limit. In the language of~\cite[Chap. 1]{georgii2011gibbs}, the kernels form a (canonical) specification.
	\item The kernels are \emph{covariant} under translations: $(f \circ \tau_a)_\La = f_{\La + a} \circ \tau_a$ for every $a \in \R$, where $\tau_a$ is the translation by $a$. The energy $F$ and the Lebesgue measure are translation-invariant, but the truncation window $[-T, T]$ in the definition of $\M^T_\La$ is not: translating by $a$ moves it by $a$. The two windows differ by two intervals of length $|a|$ at distance of order $T$. For $\La \subset [-r, r]$ and $T$ large, the charges in these intervals contribute at most
\begin{equation*}
2 n \, \frac{2r}{T - |a|} \left( \bX\left( [T - |a|, T + |a|] \cup [-T - |a|, -T + |a|] \right) + 4|a| \right)
\end{equation*}
to the difference of the truncated move functions, by~\eqref{eq:kx} and since $\bX'$ and $\bX_\La$ have the same number $n$ of points. This goes to zero as $T \to \infty$ as soon as $\bX([m, m+1]) = o(|m|)$, which holds almost surely for a good process. Hence the limiting move functions, and the kernels, are covariant.
	\item All of this holds as well for the truncated kernels $f^T_\La$, defined as in~\eqref{eq:DLR} with $\M^T_\La$ instead of $\M_\La$, except the covariance, since translations also move the truncation window. Moreover, if $f$ only depends on the configuration in $[-T_1, T_1]$, then, for $T \geq T_1$, $f^T_\La$ is local: it only depends on $\bX(\La)$ and on $\bX_{[-T, T] \setminus \La}$.
\end{itemize}

\subsection{Finite-volume DLR equations}
\label{sec:finite-DLR-section}
In order to derive (canonical) DLR equations in infinite volume, one must first write down a finite-volume version, starting from the Gibbs measure $\PNb$ itself. 

\subsubsection*{The finite-volume model in rescaled coordinates}
For $\X_N = (x_1, \dots, x_N)$, write $y_i := \vrho x_i$ so that the local configuration $\bXNz$ is $\bXNz = \sum_i \delta_{y_i}$. We let
\begin{equation*}
\nu_N := \bXNz - \widetilde{\mu}_N,
\end{equation*}
where $\widetilde{\mu}_N$ is the image of the background $\mathbf{m}_N$ under $x \mapsto \vrho x$:
\begin{equation}
\label{def:muNtilde}
\widetilde{\mu}_N(\dd y) := \rho_N(y) \dd y, \qquad \rho_N(y) := \frac{1}{\vrho} \rhosc\left( \frac{y}{\vrho N} \right).
\end{equation}
It has total mass $N$, $\rho_N(0) = 1$, $0 \leq \rho_N \leq 1$, and $|\rho_N(y) - 1| \leq \Cc \frac{y^2}{N^2}$. The law of $(y_1, \dots, y_N)$ has a density with respect to the Lebesgue measure on $\R^N$ proportional to
\begin{equation}
\label{eq:density-local}
\exp\left( - \beta \left( F(\bXNz, \widetilde{\mu}_N) + \Xi_N(\bXNz) \right) \right), \qquad \Xi_N(\bXNz) := \sum_{i=1}^N \zeta_N\left( \frac{y_i}{\vrho} \right),
\end{equation}
where, for a finite configuration $\bX$ and a measure $\mu$ on $\R$,
\begin{equation*}
F(\bX, \mu) := \hal \iint_{x \neq y} - \log|x - y| \, \dd (\bX - \mu)(x) \, \dd (\bX - \mu)(y).
\end{equation*}
Indeed, by~\eqref{eq:PNb-electric} and by scaling, $F(\bXNz, \widetilde{\mu}_N)$ is equal to $F_N(\X_N)$ up to an additive constant.

\subsubsection*{Finite-volume move function, kernels, DLR equations} 
For a bounded interval $\La$, a finite configuration $\bX$ and $\bX' \in \Conf(\La, \bX(\La))$, define the \emph{finite-volume move function} as
\begin{equation*}
\M^{(N)}_\La(\bX', \bX) := \iint_{\La \times \La^c} - \log|x - y| \, \dd (\bX' - \bX_\La)(x) \, \dd (\bX_{\La^c} - \widetilde{\mu}_N)(y),
\end{equation*}
and $\M^{T,(N)}_\La$ denotes the same integral restricted to $|y| \leq T$. 

We define the finite-volume kernels $f^{(N)}_\La$ (resp. $f^{T,(N)}_\La$) as in~\eqref{eq:DLR}, with $F_\La$ replaced by $F^{\rho_N}_\La$ and $\M_\La$ by $\M^{(N)}_\La$ (resp. by $\M^{T, (N)}_\La$). We write $\P_N, \E_N$ for the law of $\bXNz$ under $\PNb$ and the corresponding expectation (the inverse temperature $\beta$ is fixed, and we omit it from the notation).

\begin{proposition}[Finite-volume DLR equations]
\label{prop:finite-DLR}
Let $\La$ be a bounded interval, and $N$ large enough that $\La \subset [- 2 \vrho N, 2 \vrho N]$. For every bounded measurable $f$, we have $\E_N[f] = \E_N[f^{(N)}_\La]$. 

In other words, conditionally on $(\bXNz)_{\La^c}$ and on $\bXNz(\La) = n$, the law of $(\bXNz)_\La$ is
\begin{equation}
\label{eq:finite-kernel}
\mathcal{K}^{(N)}_{\La}(\bX, \cdot) := \frac{1}{\K^{(N)}_\La(\bX)} e^{- \beta \left( F^{\rho_N}_\La(\bX') + \M^{(N)}_\La(\bX', \bX) \right)} \dd \Bin_{\La, n}(\bX').
\end{equation}
\end{proposition}
\begin{proof}
This follows from standard manipulations, see similar computations in~\cite{dereudre2021dlr, dereudre2023rigidity}.
\end{proof}

\subsubsection*{Truncation estimates for finite-volume move functions}
\begin{proposition}[Truncation of finite-volume move functions]
\label{prop:finite-truncation}
Let $\La \subset [-r, r]$ be a bounded interval. For every $\epsilon > 0$, there exists $T_0$ such that, for every $T \geq T_0$,
\begin{equation*}
\limsup_{N \to \infty} \P_N\left( \sup_{\bX' \in \Conf(\La, \bXNz(\La))} \left| \M^{(N)}_\La(\bX', \bXNz) - \M^{T,(N)}_\La(\bX', \bXNz) \right| > \epsilon \right) \leq \epsilon.
\end{equation*}
\end{proposition}
{Let us sketch the proof. The only new ingredient is the control near the edge of the system. Let $\bX$ be a finite configuration, let $n := \bX(\La)$ and $\bX' \in \Conf(\La, n)$. By definition, for $T \geq 2r$, we need to control the interaction between $\bX' - \bX_\La$ and the charges at distance larger than $T$:
\begin{equation}
\label{eq:move-minus-truncated}
\M^{(N)}_\La(\bX', \bX) - \M^{T,(N)}_\La(\bX', \bX) = \iint_{\La \times \{|y| > T\}} - \log|x - y| \, \dd (\bX' - \bX_\La)(x) \, \dd \nu_N(y).
\end{equation}

Since $\bX'$ and $\bX_\La$ have the same number of points, we can subtract from $- \log|x - y|$ any function of $y$ alone without changing~\eqref{eq:move-minus-truncated}. We subtract $- \log |y|$, which amounts to comparing the potential created by $\nu_N$ at $x$ with the one at the origin. Let
\begin{equation*}
k_x(y) := \log \frac{|y|}{|x - y|}, \qquad V^T_N(x) := \int_{|y| > T} k_x(y) \, \dd \nu_N(y),
\end{equation*}
so that $V^T_N(x)$ is the potential created at $x$ by the charges at distance larger than $T$, minus its value at the origin. We get
\begin{equation}
\label{eq:move-potential}
\sup_{\bX' \in \Conf(\La, n)} \left| \M^{(N)}_\La(\bX', \bX) - \M^{T,(N)}_\La(\bX', \bX) \right| \leq 2 n \sup_{|x| \leq r} \left| V^T_N(x) \right|.
\end{equation}
{We thus need to control $V^T_N$ uniformly on $[-r, r]$.} Let $c_3 := \vrho \min(c_1, \frac12)$, with $c_1$ as in~\eqref{eq:variance}, and take $N$ large enough that $T < c_3 N$. We split $V^T_N$ into a bulk part and an edge part, $V^T_N = V^{T}_{\mathrm{bulk}} + V_{\mathrm{edge}}$, where
\begin{equation*}
V^{T}_{\mathrm{bulk}}(x) := \int_{T < |y| \leq c_3 N} k_x(y) \, \dd \nu_N(y), \qquad V_{\mathrm{edge}}(x) := \int_{|y| > c_3 N} k_x(y) \, \dd \nu_N(y).
\end{equation*}

\begin{claim}
\label{lem:exterior-field}
Let $r \geq 1$ and $T \geq 2r$. For $N$ large enough, we have:
\begin{itemize}
	\item for the bulk part,
	\begin{equation}
	\label{eq:exterior-field}
	\E_N\left[ \sup_{|x| \leq r} |V^{T}_{\mathrm{bulk}}(x)|^2 \right]^{\hal} \leq \Cc r \left( \frac{\log(2 + T)}{T} + \frac{\log N}{N} \right),
	\end{equation}
	\item for the edge part,
	\begin{equation}
	\label{eq:edge-potential}
	\sup_{|x| \leq r} |V_{\mathrm{edge}}(x)| \leq \Cc r^2, \qquad \sup_{|x| \leq r} |V_{\mathrm{edge}}(x)| \longrightarrow 0 \text{ in probability as } N \to \infty.
	\end{equation}
\end{itemize}
\end{claim}
\begin{itemize}
	\item The bulk part is treated as in~\cite{dereudre2021dlr}: we integrate by parts, which makes the discrepancies $\nu_N((T, y])$ appear, and we control their second moments with~\eqref{eq:variance}.
	\item The edge part is new: there is no local law near the edges of the system. However, seen from $[-r, r]$, the far charges look like a dipole, and $V_{\mathrm{edge}}$ is essentially a \emph{macroscopic} linear statistic. The law of large numbers for $\beta$-ensembles is then enough.
\end{itemize}
Proposition~\ref{prop:finite-truncation} follows from~\eqref{eq:move-potential}, Claim~\ref{lem:exterior-field} and the control on $n$ of Proposition~\ref{prop:FiniteVolEstimates}. The details are in Appendix~\ref{sec:proof-finite-truncation}.

\subsection{Tightness in the exp-local topology}
\label{sec:tightness}
We need to upgrade the tightness of $(\bXNz)_N$ from the weak topology to the local topology, in order to pass the DLR equations to the limit. This follows from a uniform bound on the local relative entropies with respect to the Poisson process, see Lemma~\ref{lem:entropy-compactness} below. We obtain this bound from the finite-volume DLR equations, together with a rough bound on the influence of the exterior.

\subsubsection*{Uniform entropy bounds imply tightness in the local topology}
\begin{lemma}
\label{lem:entropy-compactness}
Let $(\Pp^k)_{k \geq 1}$ be point processes on $\R$. Assume that, for every bounded interval $\La$
\begin{equation*}
\sup_k \Ent(\Pp^k_\La \, | \, \Poi_\La) < + \infty.
\end{equation*}
Then, up to extraction, the sequence $(\Pp^k)_{k \geq 1}$ converges to some point process $\Pp$ in the local topology.
\end{lemma}
The proof is classical, see for instance \cite[Prop. 15.14]{georgii2011gibbs} for the very similar statement that the sub-level sets of the specific entropy are compact for the local topology. For completeness, we include it in Appendix~\ref{sec:entropy-compactness}. Next, we obtain these entropy bounds.

\subsubsection*{Getting uniform entropy bounds}
\begin{proposition}
\label{prop:local-entropy}
For every bounded interval $\La$,
\begin{equation*}
\sup_{N \geq 1} \Ent\left( (\P_N)_\La \, | \, \Poi_\La \right) < + \infty.
\end{equation*}
\end{proposition}
\begin{proof}
Without loss of generality, assume that $\La = [-r, r]$ with $r \geq 1$. For each fixed $N$, the entropy is finite: $\bXNz$ has at most $N$ points, and the joint density of $\PNb$ is bounded, so that $(\P_N)_\La$ has a bounded density with respect to $\Poi_\La$. It is thus enough to consider $N$ large enough such that $[-2r, 2r] \subset [-c_2 N, c_2 N]$, and such that Claim~\ref{lem:exterior-field} applies with $T = 2r$.

\medskip 

We write $\bX := \bXNz$ and $n := \bX(\La)$. By Proposition~\ref{prop:FiniteVolEstimates}, the moments of $n$ and of $\bX([-2r, 2r])$ under $\P_N$ are bounded uniformly in $N$.

\paragraph{Step 1: conditioning on the exterior.} By Proposition~\ref{prop:finite-DLR}, conditionally on $(\bX_{\La^c}, n)$, the law of $\bX_\La$ under $\P_N$ is the kernel $\mathcal{K}^{(N)}_\La(\bX, \cdot)$. In other words, $(\P_N)_\La$ is the average under $\P_N$ of these kernels, and, by convexity of the relative entropy,
\begin{equation*}
\Ent\left( (\P_N)_\La \, | \, \Poi_\La \right) \leq \E_N\left[ \Ent\left( \mathcal{K}^{(N)}_\La(\bX, \cdot) \, | \, \Poi_\La \right) \right].
\end{equation*}
The kernel $\mathcal{K}^{(N)}_\La(\bX, \cdot)$ is concentrated on $\Conf(\La, n)$, and $\Poi_\La$ conditioned to have $n$ points is $\Bin_{\La, n}$, hence:
\begin{equation*}
\Ent\left( \mathcal{K}^{(N)}_\La(\bX, \cdot) \, | \, \Poi_\La \right) = \Ent\left( \mathcal{K}^{(N)}_\La(\bX, \cdot) \, | \, \Bin_{\La, n} \right) - \log \left( e^{-|\La|} \frac{|\La|^n}{n!} \right),
\end{equation*}
and the last term is at most $|\La| + n |\log |\La|| + n^2$.

\paragraph{Step 2: Jensen's inequality.} In view of~\eqref{eq:finite-kernel}, the kernel $\mathcal{K}^{(N)}_\La(\bX, \cdot)$ has density $\frac{1}{\K^{(N)}_\La(\bX)} e^{- \beta \mathsf{H}_\bX}$ with respect to $\Bin_{\La, n}$, where
\begin{equation*}
\mathsf{H}_\bX(\bX') := F^{\rho_N}_\La(\bX') + \M^{(N)}_\La(\bX', \bX), \qquad \K^{(N)}_\La(\bX) = \int e^{- \beta \mathsf{H}_\bX} \dd \Bin_{\La, n}.
\end{equation*}
Hence
\begin{equation*}
\Ent\left( \mathcal{K}^{(N)}_\La(\bX, \cdot) \, | \, \Bin_{\La, n} \right) = - \beta \, \E_{\mathcal{K}^{(N)}_\La(\bX, \cdot)}\left[ \mathsf{H}_\bX \right] - \log \K^{(N)}_\La(\bX),
\end{equation*}
and Jensen's inequality, $- \log \K^{(N)}_\La(\bX) \leq \beta \, \E_{\Bin_{\La, n}}[\mathsf{H}_\bX]$, gives
\begin{equation}
\label{eq:entropy-jensen}
\Ent\left( \mathcal{K}^{(N)}_\La(\bX, \cdot) \, | \, \Bin_{\La, n} \right) \leq \beta \left( \E_{\Bin_{\La, n}}\left[ \mathsf{H}_\bX \right] - \E_{\mathcal{K}^{(N)}_\La(\bX, \cdot)}\left[ \mathsf{H}_\bX \right] \right).
\end{equation}
We now study the average under $\P_N$ of both expectations in the right-hand side (which depend on $\bX$).

\paragraph{Step 3: averaging over the exterior.} The average under $\P_N$ of the second term of~\eqref{eq:entropy-jensen} can be computed with the DLR equations. For $\bX' \in \Conf(\La, n)$, let
\begin{equation*}
\mathsf{A}_\bX(\bX') := \iint_{\La \times \La^c} - \log |x - y| \, \dd \bX'(x) \, \dd (\bX - \widetilde{\mu}_N)(y)
\end{equation*}
be the interaction of $\bX'$ with the charges outside $\La$. It only depends on $\bX$ through $\bX_{\La^c}$, and, by definition of the move function,
\begin{equation*}
\M^{(N)}_\La(\bX', \bX) = \mathsf{A}_\bX(\bX') - \mathsf{A}_\bX(\bX_\La), \qquad \text{hence} \qquad \mathsf{H}_\bX(\bX') = F^{\rho_N}_\La(\bX') + \mathsf{A}_\bX(\bX') - \mathsf{A}_\bX(\bX_\La).
\end{equation*}
We apply Proposition~\ref{prop:finite-DLR} to the function $g(\bX) := F^{\rho_N}_\La(\bX_\La) + \mathsf{A}_\bX(\bX_\La)$, which is not bounded, but is integrable under $\P_N$. This gives:
\begin{equation*}
\E_N\left[ \E_{\mathcal{K}^{(N)}_\La(\bX, \cdot)}\left[ F^{\rho_N}_\La + \mathsf{A}_\bX \right] \right] = \E_N\left[ F^{\rho_N}_\La(\bX_\La) + \mathsf{A}_\bX(\bX_\La) \right],
\end{equation*}
and thus, subtracting $\E_N[\mathsf{A}_\bX(\bX_\La)]$ from both sides,
\begin{equation*}
\E_N\left[ \E_{\mathcal{K}^{(N)}_\La(\bX, \cdot)}\left[ \mathsf{H}_\bX \right] \right] = \E_N\left[ F^{\rho_N}_\La(\bX_\La) \right].
\end{equation*}
For the first term of~\eqref{eq:entropy-jensen}, let $\bY$ be a configuration of $n$ independent uniform points in $\La$, independent of $\bX$ given $n$. Then $\E_{\Bin_{\La, n}}[\mathsf{H}_\bX] = \E[ F^{\rho_N}_\La(\bY) + \M^{(N)}_\La(\bY, \bX) \, | \, \bX]$. Averaging~\eqref{eq:entropy-jensen} under $\P_N$, we obtain
\begin{equation}
\label{eq:entropy-three-terms}
\E_N\left[ \Ent\left( \mathcal{K}^{(N)}_\La(\bX, \cdot) \, | \, \Bin_{\La, n} \right) \right] \leq \beta \, \E\left[ F^{\rho_N}_\La(\bY) - F^{\rho_N}_\La(\bX_\La) + \M^{(N)}_\La(\bY, \bX) \right].
\end{equation}
So the entropy is bounded by the energy cost of replacing the points in $\La$ by independent uniform points. We bound the three terms in the next two steps.

\paragraph{Step 4: self-energies.} Since $\La$ is fixed and $0 \leq \rho_N \leq 1$, the self-energy of $n$ points in $\La$ is bounded below deterministically, and its expectation under independent uniform points is explicit:
\begin{equation}
\label{eq:crude-F}
F^{\rho_N}_\La(\bX') \geq - \Cc_r (1 + n^2) \quad \text{for all } \bX' \in \Conf(\La, n), \qquad \E\left[ F^{\rho_N}_\La(\bY) \right] \leq \Cc_r (1 + n^2).
\end{equation}
Hence
\begin{equation*}
\E\left[ F^{\rho_N}_\La(\bY) - F^{\rho_N}_\La(\bX_\La) \right] \leq \Cc_r \left( 1 + \E_N[n^2] \right).
\end{equation*}
We now focus on the contribution of the move function, i.e. of the exterior.

\paragraph{Step 5: interaction with the exterior.} By definition,
\begin{equation*}
\M^{(N)}_\La(\bY, \bX) = \iint_{\La \times \La^c} - \log |x - y| \, \dd (\bY - \bX_\La)(x) \, \dd (\bX - \widetilde{\mu}_N)(y).
\end{equation*}
We split $\La^c$ into the near region $[-2r, 2r] \setminus \La$ and the far region $\{|y| > 2r\}$.

\paragraph{Near region.} There are three contributions.
\begin{itemize}
	\item The interaction of $\bY$ with the points $z$ of $\bX_{\La^c}$ in the near region. Since the logarithm is locally integrable,
	\begin{equation*}
	\E\left[ \sum_{x \in \bY} \sum_{z \in \bX_{\La^c} \cap [-2r, 2r]} - \log |x - z| \, \middle| \, \bX \right] \leq \Cc_r \, n \, \bX([-2r, 2r]).
	\end{equation*}
	\item The interaction of $- \bX_\La$ with the same points. Here two points may be arbitrarily close to each other, on both sides of $\partial \La$, but this term comes with a minus sign, and
	\begin{equation*}
	\sum_{x \in \bX_\La} \sum_{z \in \bX_{\La^c} \cap [-2r, 2r]} \log |x - z| \leq n \, \bX([-2r, 2r]) \log (4r).
	\end{equation*}
	\item The interaction of $\bY - \bX_\La$ with the background. For $x \in \La$, since $\rho_N \leq 1$,
	\begin{equation*}
	\int_{[-2r, 2r] \setminus \La} |\log |x - y|| \, \rho_N(y) \dd y \leq \int_{-4r}^{4r} |\log |t|| \dd t \leq \Cc_r,
	\end{equation*}
	so that this contribution is at most $2 \Cc_r n$.
\end{itemize}
By the Cauchy--Schwarz inequality, the expectation of the near region contribution is bounded by
\begin{equation*}
\Cc_r \left( \E_N[n^2] + \E_N\left[\bX([-2r, 2r])^2\right] + \E_N[n] \right).
\end{equation*}

\paragraph{Far region.} As in~\eqref{eq:move-potential}, the far region contribution is bounded by $2n \sup_{|x| \leq r} |V^{2r}_N(x)|$. By Claim~\ref{lem:exterior-field} with $T = 2r$, and by the Cauchy--Schwarz inequality, its expectation is bounded by
\begin{equation*}
2 \E_N[n^2]^{\hal} \left( \E_N\left[ \sup_{|x| \leq r} |V^{2r}_{\mathrm{bulk}}(x)|^2 \right]^{\hal} + \Cc r^2 \right) \leq \Cc_r \, \E_N[n^2]^{\hal}.
\end{equation*}

\paragraph{Step 6: conclusion.} $\Ent((\P_N)_\La \, | \, \Poi_\La)$ is bounded by $\Cc_r \left(1 + \E_N[\bX([-2r, 2r])^2]\right)$, which is bounded uniformly in $N$.
\end{proof}

\subsubsection*{From local to exp-local convergence}
Combining Proposition~\ref{prop:local-entropy} with Lemma~\ref{lem:entropy-compactness} gives tightness in the local topology. Thanks to the strong controls on the number of points of Proposition~\ref{prop:FiniteVolEstimates}, we can upgrade it to the exp-local topology, and pass these controls to the limit. The following result contains the tightness part and item 1 of Theorem~\ref{theo:Convergence}.

\begin{theorem}[Tightness and properties of limit points]
\label{theo:tightness}
The laws of $(\bXNz)_N$ are tight in the exp-local topology. Every limit point $\Pp$ of the laws of $(\bXNz)_N$ (for the weak topology, or equivalently for the exp-local topology) satisfies:
\begin{enumerate}
	\item for every bounded interval $\La$, $\Ent(\Pp_\La \, | \, \Poi_\La) < + \infty$. In particular, $\Pp$-almost surely, the configuration is simple, and $\Pp$ has no fixed atom,
	\item for every interval $I$ of length $R \geq 1$, $\log \E_\Pp\left[ e^{a \bX(I)^2} \right] \leq \Cc R^2$,
	\item for every interval $I$ of length $R \geq 1$, $\E_\Pp\left[ \left( \bX(I) - R \right)^2 \right] \leq \Cc \log^2 (2 + R)$.
\end{enumerate}
In particular, every limit point is good in the sense of Definition~\ref{def:good}.
\end{theorem}
\begin{proof}
The local convergence is stronger than the weak convergence, so that the limit points for both topologies coincide. The rest follows from Lemma~\ref{lem:entropy-compactness}, Proposition~\ref{prop:local-entropy} and the controls of Proposition~\ref{prop:FiniteVolEstimates}, see~\cite[Lemmas 4.3 and 4.4]{leble2024dlr2d} for the details. For item 2, we apply the local convergence to the bounded local functions $\min(e^{a \bX(I)^2}, k)$, and let $k \to \infty$ by monotone convergence.
\end{proof}

\subsection{DLR equations for limit points}
\label{sec:DLR-limit}
We now prove item 2 of Theorem~\ref{theo:Convergence}, that we recall.
\begin{theorem*}[DLR equations for limit points, item 2 of Theorem~\ref{theo:Convergence}]
Every limit point $\Pp$ of the laws of $(\bXNz)_N$ is a solution of~\eqref{DLR} in the sense of Definition~\ref{def:DLR}.
\end{theorem*}

We pass the finite-volume DLR equations of Proposition~\ref{prop:finite-DLR} to the limit, following~\cite{dereudre2021dlr}, see also~\cite[Section 7]{leble2024dlr2d}. The finite-volume kernels $f^{(N)}_\La$ differ from the infinite-volume ones $f_\La$ in two ways: they see the background $\widetilde{\mu}_N$ instead of the Lebesgue measure, and they are not local, because the move function feels all the particles. We thus make three approximations:
\begin{enumerate}
	\item truncate the move functions at distance $T$, uniformly in $N$ (Proposition~\ref{prop:finite-truncation}) and in the limit (Proposition~\ref{prop:move}),
	\item flatten the background, which is harmless near the origin (Lemma~\ref{lem:flatten}),
	\item use the convergence in the local topology, which does not require any continuity of the truncated kernels.
\end{enumerate}
Compared with~\cite{dereudre2021dlr}, only the first step is new: since the finite-volume model is not stationary, the uniformity of the truncation in $N$ comes from the local laws and from the law of large numbers, see Claim~\ref{lem:exterior-field}. The rest of the proof is standard, and we give it for completeness in Appendix~\ref{sec:proof-DLR-limit}.

\section{Translation-invariance}
\label{sec:stationarity}
In this section, we prove Theorem~\ref{theo:TI-intro}. The finite-volume $\beta$-ensembles are not translation-invariant, so it is not obvious that the limit points of $\bXNz$ should be stationary. The following theorem shows that every {good} solution of {the} DLR {equations} has this symmetry. Note that this is false for closely related systems, for instance the one-dimensional Coulomb gas, with interaction potential $- |\cdot|$ instead of $- \log|\cdot|$, which admits periodic infinite-volume Gibbs states at all temperatures~\cite{aizenman1980structure,aizenman2001bounded}.  

\begin{theorem*}[Translation-invariance, Theorem~\ref{theo:TI-intro}]
Every good solution of~\eqref{DLR}, in the sense of Definition~\ref{def:DLR}, is stationary. In particular, every limit point of the laws of $(\bXNz)_N$ is stationary.
\end{theorem*}

\subsubsection*{Strategy of proof}
We follow the classical method of Fröhlich and Pfister~\cite{frohlich1981absence}, see also Georgii~\cite{georgii1999translation}, in the form used for the two-dimensional Coulomb gas in~\cite[Section 8]{leble2024dlr2d}. The idea is to construct \emph{localized translations}: maps of $\R$ which act as a translation in a large region, and are “dampened” so as to coincide with the identity far away. The DLR equations allow us to compute the effect of such maps on $\Pp$: if their energy cost (including the Jacobian) is bounded, uniformly with respect to the size of the region, then $\Pp$ is absolutely continuous with respect to its translates, and a classical argument on the structure of the set of Gibbs measures shows that $\Pp$ is in fact invariant. It is often convenient (and sufficient) to consider the \emph{average} effect of two localized translations in opposite directions.

The effect of a transport map on the logarithmic energy has been studied in detail in finite volume, in the context of fluctuations of linear statistics~\cite{bekerman2018clt, peilen2022local}: if $\phi = \mathrm{Id} + \psi$ is a small perturbation of the identity, then for $x \neq y$,
\begin{equation}
\label{eq:transport-energy}
\red{- \log |\phi(x) - \phi(y)| + \log |x - y|} = - \log \left( 1 + \frac{\psi(x) - \psi(y)}{x - y} \right) = - \frac{\psi(x) - \psi(y)}{x - y} + O\left( \left( \frac{\psi(x) - \psi(y)}{x - y} \right)^2 \right).
\end{equation}
Summed over all pairs, the first-order term gives the “anisotropy” term $- \hal \iint \frac{\psi(x) - \psi(y)}{x - y} \dd \nu \, \dd \nu$ of~\cite{bekerman2018clt}, where $\nu$ is the difference between the configuration and the background, whose control is the difficult part of the transport approach, while the second-order term is a sum of non-negative, explicit quantities. Two features are specific to dimension one:
\begin{itemize}
	\item A localized translation of the line cannot preserve the Lebesgue measure, so that the Jacobian must be taken into account.
	\item Symmetrizing two opposite localized translations cancels the first-order term of~\eqref{eq:transport-energy} \emph{exactly}, and the symmetrized cost is non-negative. As a consequence, a \emph{first} moment bound on the cost is enough, whereas~\cite{leble2024dlr2d} uses exponential moments. This bound only uses second moments of the number of points (Lemma~\ref{lem:cost-L1}). The control on the fluctuations given by the good hypothesis is only used to define the move functions, hence the DLR kernels.
\end{itemize}

\subsection{Localized translations and their symmetrized cost}

\subsubsection*{Localized translations}
Fix $\psi \in C^2_c(\R, [0,1])$, equal to $1$ on $[-\frac12, \frac12]$ and to $0$ outside $[-1, 1]$. For $L \geq \max(1, 2 \|\psi'\|_\infty)$, our \emph{localized translations} are the maps
\begin{equation*}
T^\pm_L(x) := x \pm \psi\left( \frac{x}{L} \right), \qquad \text{with Jacobian} \qquad (T^\pm_L)'(x) = 1 \pm \frac{1}{L} \psi'\left( \frac{x}{L} \right).
\end{equation*}
Since $|(T^\pm_L)' - 1| \leq \frac12$, they satisfy the following properties.
\begin{itemize}
	\item They act as the translations by $\pm 1$ on $[-\frac{L}{2}, \frac{L}{2}]$.
	\item They are increasing bijections of $\R$.
	\item They coincide with the identity outside $\La_{2L} = [-L, L]$, and map $\La_{2L}$ onto itself.
\end{itemize}
They are not exactly inverses of each other, but they are on a large portion of $[-\frac{L}{2}, \frac{L}{2}]$.

\subsubsection*{Cost of a pair of points}
Let us introduce the auxiliary function
\begin{equation*}
d_L(x, y) := \frac{\psi(x/L) - \psi(y/L)}{x - y} \quad (x \neq y), \qquad \text{which satisfies} \qquad |d_L(x,y)| \leq \frac{\|\psi'\|_\infty}{L} \leq \frac12.
\end{equation*}
\begin{claim}
For $x \neq y$, the symmetrized variation of the interaction of $x$ and $y$,
\begin{equation*}
\mathrm{cost}_L(x,y) := \red{- \frac{\log |T^+_L x - T^+_L y| + \log |T^-_L x - T^-_L y|}{2} + \log |x - y|},
\end{equation*}
satisfies
\begin{equation}
\label{eq:sym-cost}
0 \leq \mathrm{cost}_L(x,y) \leq d_L(x,y)^2.
\end{equation}
\end{claim}
\begin{proof}
Since $T^\pm_L x - T^\pm_L y = (x - y)\left(1 \pm d_L(x,y)\right)$, the first identity of~\eqref{eq:transport-energy} gives
\begin{equation*}
\mathrm{cost}_L(x,y) = - \hal \log\left( (1 + d_L(x,y))(1 - d_L(x,y)) \right) = - \hal \log\left( 1 - d_L(x,y)^2 \right),
\end{equation*}
and $0 \leq - \hal \log(1 - u) \leq u$ for $u \in [0, \frac14]$.
\end{proof}
The first-order terms have cancelled, and the \red{logarithmic singularity} on the diagonal has disappeared.

\subsubsection*{Energy cost}
For a configuration $\bX$, we define the \emph{energy cost}
\begin{equation}
\label{def:energy-cost}
\mathrm{Cost}^{\mathrm{en}}_L(\bX) := \hal \sum_{\substack{x \neq y \in \bX \\ x, y \in \La_{2L}}} \mathrm{cost}_L(x, y) + \sum_{\substack{x \in \bX \cap \La_{2L} \\ y \in \bX \setminus \La_{2L}}} \mathrm{cost}_L(x, y),
\end{equation}
with values in $[0, +\infty]$ (points are counted with multiplicity). It is the variation of the interaction energy between the average of the configurations $T^+_L \bX$ and $T^-_L \bX$, and the configuration $\bX$. Only the pairs with at least one point in $\La_{2L}$ contribute, since $T^\pm_L$ is the identity outside $\La_{2L}$.

\subsubsection*{Cost of the Jacobians}
The Jacobians combine in the same way: for a configuration $\bX$, let
\begin{equation*}
\mathrm{Cost}^{\mathrm{jac}}_L(\bX) := - \hal \sum_{x \in \bX \cap \La_{2L}} \log\left( (T^+_L)'(x) \, (T^-_L)'(x) \right) = - \hal \sum_{x \in \bX \cap \La_{2L}} \log\left( 1 - \frac{1}{L^2} \psi'\left( \frac{x}{L} \right)^2 \right).
\end{equation*}
For the same reason as above,
\begin{equation}
\label{eq:jac-cost}
0 \leq \mathrm{Cost}^{\mathrm{jac}}_L(\bX) \leq \frac{\|\psi'\|^2_\infty}{L^2} \, \bX(\La_{2L}).
\end{equation}

\subsubsection*{The symmetrized cost}
For a configuration $\bX$, we define the \emph{symmetrized cost}
\begin{equation}
\label{def:CL}
\mathrm{Cost}_L(\bX) := \beta \, \mathrm{Cost}^{\mathrm{en}}_L(\bX) + \mathrm{Cost}^{\mathrm{jac}}_L(\bX).
\end{equation}
In view of~\eqref{eq:sym-cost} and~\eqref{eq:jac-cost}, it is non-negative, and it can be bounded in terms of the number of points. In fact, its expectation is bounded uniformly in $L$, see Lemma~\ref{lem:cost-L1} below: this is the condition needed to make the argument of Fröhlich and Pfister work, see Proposition~\ref{prop:abs-cont}.

\begin{lemma}
\label{lem:cost-L1}
Let $\Pp$ be a point process such that, for every interval $I$, $\E_\Pp[\bX(I)^2] \leq \Cc_\Pp (|I| + 1)^2$. Then 
\begin{equation*}
\sup_{L} \E_\Pp[\mathrm{Cost}_L] \leq \Cc,
\end{equation*} 
where $\Cc$ depends only on $\beta$, $\psi$ and $\Cc_\Pp$.
\end{lemma}
\begin{proof}
By~\eqref{eq:sym-cost}, we have $\mathrm{cost}_L \leq d_L^2$. For the pairs of points in $\La_{4L} = [-2L, 2L]$ (with at least one of them in $\La_{2L}$), we have
\begin{equation*}
d_L^2 \leq \frac{\|\psi'\|^2_\infty}{L^2}, \qquad \text{so their contribution is at most} \qquad \|\psi'\|^2_\infty \frac{\bX(\La_{4L})^2}{L^2},
\end{equation*}
whose expectation is $O(1)$ by our assumption. For $x \in \La_{2L}$ and $|y| > 2L$, we have
\begin{equation*}
|d_L(x,y)| = \frac{\psi(x/L)}{|x - y|} \leq \frac{2}{|y|}.
\end{equation*} With the shells $S_k := \{ y : 2^k L < |y| \leq 2^{k+1} L \}$, $k \geq 1$, the contribution of these pairs is at most
\begin{equation*}
4 \sum_{k \geq 1} \frac{\bX(\La_{2L}) \, \bX(S_k)}{(2^k L)^2}, \quad \text{with expectation at most} \quad 4 \sum_{k \geq 1} \frac{\Cc L \cdot \Cc 2^k L}{4^k L^2} \leq \Cc,
\end{equation*}
by the Cauchy--Schwarz inequality and our assumption. Finally, $\E_\Pp[\mathrm{Cost}^{\mathrm{jac}}_L] \leq \frac{\Cc}{L}$ by~\eqref{eq:jac-cost}.
\end{proof}

\subsection{The cost in the DLR kernels}

We introduce a subset of admissible configurations.
\begin{definition}
Let $\Omega_{\mathrm{adm}}$ be the set of configurations $\bX$ such that $\bX([m, m+1]) = o(|m|)$ as $|m| \to \infty$, and such that, for every bounded interval $\La$, the move functions $\M_\La(\bX', \bX)$ exist for all $\bX' \in \Conf(\La, \bX(\La))$.  This set is invariant under translations and has full measure for every good point process.
\end{definition}

Fix $L$ as above and take $\La := \La_{2L}$. For $\bX \in \Omega_{\mathrm{adm}}$ and $\bX' \in \Conf(\La, \bX(\La))$, let
\begin{equation*}
\mathsf{H}_\bX(\bX') := F_\La(\bX') + \M_\La(\bX', \bX)
\end{equation*}
be the energy in the DLR kernel~\eqref{eq:DLR}, and define 
\begin{equation*}
\mathsf{R}^\pm(\bX') := \beta \left( \mathsf{H}_\bX(T^\pm_L \bX') - \mathsf{H}_\bX(\bX') \right) - \sum_{x \in \bX'} \log (T^\pm_L)'(x),
\end{equation*}
where $T^\pm_L \bX'$ is the image configuration, which belongs to $\Conf(\La, \bX(\La))$. By definition, these quantities only depend on $\bX'$ and on $\bX_{\La^c}$. The quantity $\mathsf{R}^\pm$ takes two things into account:
\begin{itemize}
	\item The effect of our localized translations on the energy, via the difference $\mathsf{H}_\bX(T^\pm_L \bX') - \mathsf{H}_\bX(\bX')$.
	\item Their effect on the \emph{volume}, through the Jacobian $(T^\pm_L)'$. Its logarithm goes into the exponent, next to the energy.
\end{itemize}

\begin{lemma}
\label{lem:identification}
For $\bX \in \Omega_{\mathrm{adm}}$ and $\Bin_{\La, \bX(\La)}$-a.e. $\bX' \in \Conf(\La, \bX(\La))$, we have
\begin{equation*}
\hal \left( \mathsf{R}^+ + \mathsf{R}^- \right)(\bX') = \mathrm{Cost}_L(\bX' + \bX_{\La^c}).
\end{equation*}
In particular, $\mathrm{Cost}_L(\bX) < \infty$ for every $\bX \in \Omega_{\mathrm{adm}}$ with no point on the boundary of $\La$ and no multiple point.
\end{lemma}
\begin{proof}
The Jacobian part of $\hal(\mathsf{R}^+ + \mathsf{R}^-)$ is $\mathrm{Cost}^{\mathrm{jac}}_L(\bX' + \bX_{\La^c})$, by definition. For the energy, we truncate the exterior configuration at $|y| \leq T$, $T > L$, as in Definition~\ref{def:move}, and compute the symmetrized variation of $F_\La + \M^T_\La(\cdot, \bX)$. It has three contributions.
\begin{itemize}
	\item The pairs of points of $\bX'$ contribute $\hal \sum_{x \neq y \in \bX'} \mathrm{cost}_L(x,y)$, by definition of $\mathrm{cost}_L$.
	\item The pairs made of a point of $\bX'$ and a point of $\bX$ in $[-T, T] \setminus \La$ contribute
	\begin{equation*}
	\sum_{x \in \bX', \, y \in \bX \cap ([-T,T] \setminus \La)} \mathrm{cost}_L(x,y),
	\end{equation*}
	by definition of $\mathrm{cost}_L$ again, since the exterior points do not move ($T^\pm_L$ is the identity outside $\La$).
	\item The background terms of $F_\La$ and of $\M^T_\La$ \emph{recombine}: up to terms that do not depend on $\bX'$, they are equal to $- \sum_{x \in \bX'} V_T(x)$, where $V_T(x) := \int_{-T}^T \red{- \log |x - t|} \dd t$ is the potential of the uniform background on $[-T, T]$. Their symmetrized variation is thus
	\begin{equation*}
	- \sum_{x \in \bX'} \left[ \hal \left( V_T\left(x + \psi(\tfrac{x}{L})\right) + V_T\left(x - \psi(\tfrac{x}{L})\right) \right) - V_T(x) \right].
	\end{equation*}
	Since $V_T''(x) = - \frac{2T}{T^2 - x^2}$ for $|x| < T$, this is bounded by
	\begin{equation*}
	\Cc \frac{T}{T^2 - L^2} \bX(\La),
	\end{equation*}
	and it tends to zero as $T \to \infty$.
\end{itemize}
As $T \to \infty$, the truncated move functions converge (since $\bX \in \Omega_{\mathrm{adm}}$), and the sums over pairs converge monotonically to the corresponding sums in~\eqref{def:CL}, which proves the identity. 
\end{proof}

\subsection{Mutual absolute continuity of translated DLR solutions}
We now apply the method of~\cite{frohlich1981absence}. Assume first that the symmetrized cost is bounded: $\mathrm{Cost}_L \leq K$ for all $L$ and all configurations. Let $\Qq$ be a solution of~\eqref{DLR} and $A$ be a local event. The DLR equations in $\La_{2L}$, the convexity of the exponential and the change of variables $\bX' \mapsto T^\pm_L \bX'$ in the DLR kernels give, for $L$ large enough (depending on $A$),
\begin{equation*}
\Qq(A) \leq \frac{e^K}{2} \left( \Qq^{1}(A) + \Qq^{-1}(A) \right),
\end{equation*}
where $\Qq^{\pm 1}$ are the translates of $\Qq$ by $\pm 1$. This extends to all events, so that $\Qq \ll \hal(\Qq^{1} + \Qq^{-1})$.

In fact, it is enough that the cost be bounded on average, along a subsequence: Markov's inequality then adds an error $\frac{1}{K} \E_\Qq[\mathrm{Cost}_L]$ to the previous bound, which is still enough for absolute continuity. Richthammer~\cite{richthammer2007translation} also uses costs that are only controlled in probability, for singular interactions. We need this generality in Section~\ref{sec:extremal}: the extremal components of a good solution need not be good, and they only inherit a bound on $\liminf_L \E_\Qq[\mathrm{Cost}_L]$.

\begin{proposition}
\label{prop:abs-cont}
Let $\Qq$ be a probability measure on $\Omega_{\mathrm{adm}}$ satisfying~\eqref{DLR}, and such that
\begin{equation*}
\liminf_{L \to \infty} \E_{\Qq}[\mathrm{Cost}_L] < \infty.
\end{equation*}
Let $\Qq^{\pm 1}$ be the law, under $\Qq$, of the translated configuration $\tau_{\pm 1} \bX$. Then $\Qq$ is absolutely continuous with respect to $\hal (\Qq^{1} + \Qq^{-1})$.
\end{proposition}
\begin{proof}
Let $(L_k)_k$ be a sequence tending to infinity with $\E_{\Qq}[\mathrm{Cost}_{L_k}] \leq \Cc_\Qq$ for all $k$. Let $A$ be a local event, depending only on the configuration in $[-R, R]$ for some $R \geq 1$, and choose $k$ large enough that $\frac{L_k}{2} \geq R + 1$. We write $\La := \La_{2L_k}$, $n := \bX(\La)$, and $w_\bX(\bX') := e^{-\beta \mathsf{H}_\bX(\bX')}$ for the weight of the DLR kernel in $\La$, with normalization $\K_\La(\bX)$. We use the Lebesgue measure $\dd \bX'$ on $\La^n$ instead of $\Bin_{\La, n}$, since the normalizations cancel.

\paragraph{Inserting the cost.} By the DLR equations, Lemma~\ref{lem:identification}, and the convexity of the exponential,
\begin{multline*}
\E_\Qq\left[ \ind_A e^{-\mathrm{Cost}_{L_k}} \right] = \E_\Qq\left[ \frac{1}{\K_\La(\bX)} \int \ind_A(\bX' + \bX_{\La^c}) \, e^{-\frac12 (\mathsf{R}^+ + \mathsf{R}^-)(\bX')} \, w_\bX(\bX') \dd \bX' \right] \\ \leq \hal \sum_{\pm} \E_\Qq\left[ \frac{1}{\K_\La(\bX)} \int \ind_A(\bX' + \bX_{\La^c}) \, e^{-\mathsf{R}^\pm(\bX')} \, w_\bX(\bX') \dd \bX' \right].
\end{multline*}

\paragraph{Change of variables.} Consider the “$+$” term. By definition of $\mathsf{R}^+$,
\begin{equation*}
e^{-\mathsf{R}^+(\bX')} w_\bX(\bX') = w_\bX(T^+_{L_k} \bX') \prod_{x \in \bX'} (T^+_{L_k})'(x),
\end{equation*}
so that the change of variables $\mathsf{Z} = T^+_{L_k} \bX'$ in $\La^n$ gives
\begin{equation*}
\int \ind_A(\bX' + \bX_{\La^c}) \, e^{-\mathsf{R}^+(\bX')} \, w_\bX(\bX') \dd \bX' = \int \ind_A\left( (T^+_{L_k})^{-1} \mathsf{Z} + \bX_{\La^c} \right) w_\bX(\mathsf{Z}) \dd \mathsf{Z}.
\end{equation*}
The event $A$ only sees the points of $(T^+_{L_k})^{-1} \mathsf{Z}$ in $[-R, R]$. Since $T^+_{L_k}$ acts as the translation by $+1$ on $[-\frac{L_k}{2}, \frac{L_k}{2}] \supset [-R, R]$, these are the images by $z \mapsto z - 1$ of the points of $\mathsf{Z}$ in $[-R + 1, R + 1]$, and $\bX_{\La^c}$ has no point there. Hence
\begin{equation*}
\ind_A\left( (T^+_{L_k})^{-1} \mathsf{Z} + \bX_{\La^c} \right) = f(\mathsf{Z} + \bX_{\La^c}), \qquad \text{where} \qquad f(\bX) := \ind_A(\tau_{-1} \bX).
\end{equation*}

\paragraph{Using the DLR equations again.} We have obtained that the “$+$” term is equal to
\begin{equation*}
\E_\Qq\left[ \frac{1}{\K_\La(\bX)} \int f(\mathsf{Z} + \bX_{\La^c}) \, w_\bX(\mathsf{Z}) \dd \mathsf{Z} \right] = \E_\Qq\left[ f_\La \right],
\end{equation*}
where $f_\La$ is the DLR kernel in $\La$ applied to $f$, as in~\eqref{eq:DLR}. By~\eqref{DLR}, this is equal to
\begin{equation*}
\E_\Qq[f] = \Qq(\tau_{-1} \bX \in A) = \Qq^{-1}(A).
\end{equation*} Similarly, the “$-$” term is equal to $\Qq^{1}(A)$. We have shown
\begin{equation*}
\E_\Qq\left[ \ind_A e^{-\mathrm{Cost}_{L_k}} \right] \leq \hal \left( \Qq^{1}(A) + \Qq^{-1}(A) \right).
\end{equation*}

\paragraph{Conclusion.} For $K > 0$, since $\mathrm{Cost}_{L_k} \geq 0$, Markov's inequality gives
\begin{equation}
\label{eq:abs-cont}
\Qq(A) \leq e^K \E_\Qq\left[ \ind_A e^{-\mathrm{Cost}_{L_k}} \right] + \Qq(\mathrm{Cost}_{L_k} > K) \leq \frac{e^K}{2} \left( \Qq^{1}(A) + \Qq^{-1}(A) \right) + \frac{\Cc_\Qq}{K}.
\end{equation}
For fixed $K$, the set of events $A$ satisfying
\begin{equation*}
\Qq(A) \leq \frac{e^K}{2} \left( \Qq^{1}(A) + \Qq^{-1}(A) \right) + \frac{\Cc_\Qq}{K}
\end{equation*}
is a monotone class, which contains the algebra of local events, hence it contains all measurable events. If $\Qq^{1}(A) + \Qq^{-1}(A) = 0$, letting $K \to \infty$ gives $\Qq(A) = 0$.
\end{proof}

\subsection{The extremal decomposition: proof of Theorem~\ref{theo:TI-intro}}
\label{sec:extremal}
Absolute continuity with respect to the translates does not, by itself, imply stationarity. To conclude, we use the extremal decomposition of DLR solutions, as in~\cite{bricmont1979translation}. Let $\mathscr{G}$ be the (convex) set of probability measures on $\Omega_{\mathrm{adm}}$ satisfying~\eqref{DLR}. We use the following classical facts, see~\cite[Chapter 7]{georgii2011gibbs}, and~\cite{preston1976random} for specifications of canonical type such as ours. The abstract theory applies to any consistent family of proper kernels on a standard Borel space, indexed by a countable directed set. Here, by consistency, it is enough to consider the intervals $[-k, k]$, $k \in \N$.
\begin{itemize}
	\item The extremal elements of $\mathscr{G}$ are exactly those which are trivial on the tail $\sigma$-algebra $\mathcal{G}_\infty := \bigcap_\La \mathcal{G}_\La$, and two distinct extremal elements are mutually singular.
	\item Every $\Pp \in \mathscr{G}$ has an extremal decomposition: $\Pp = \int \Qq \, \kappa(\dd \Qq)$, for a certain probability measure $\kappa$ on the set of extremal elements of $\mathscr{G}$. 
	\item By covariance of the interaction energy, $\Qq \mapsto \Qq^{\pm 1}$ maps $\mathscr{G}$ onto itself, and preserves extremality.
\end{itemize}
Let us comment on the first item, since we are working with a restricted set of configurations, which is not quite the usual setting. If $\Pp \in \mathscr{G}$ and $E \in \mathcal{G}_\infty$ with $0 < \Pp(E) < 1$, then $\Pp(\cdot \, | \, E)$ and $\Pp(\cdot \, | \, E^c)$ still satisfy~\eqref{DLR} (this only uses the structure of the DLR equations: the kernel in $\La$ does not modify $\mathcal{G}_\La$-measurable quantities, see~\cite[Prop. 6.61]{friedli2017statistical}), and they are still concentrated on $\Omega_{\mathrm{adm}}$. Hence an extremal element is tail-trivial. The converse, and the mutual singularity, follow as in~\cite[Lemma 6.62, Theorem 6.58]{friedli2017statistical}. Note that we do not need the components $\Qq$ of the extremal decomposition to be good with uniform constants: only the bound on $\liminf_L \E_\Qq[\mathrm{Cost}_L]$ will be used.

\begin{proof}[Proof of Theorem~\ref{theo:TI-intro}]
Let $\Pp$ be a good point process satisfying the DLR equations. Since $\Pp(\Omega_{\mathrm{adm}}) = 1$, we have $\Pp \in \mathscr{G}$. Let $\Pp = \int \Qq \, \kappa(\dd \Qq)$ be its extremal decomposition. By Fatou's lemma and Lemma~\ref{lem:cost-L1},
\begin{equation*}
\int \liminf_{L \to \infty} \E_\Qq[\mathrm{Cost}_L] \, \kappa(\dd \Qq) \leq \liminf_{L \to \infty} \E_\Pp[\mathrm{Cost}_L] < \infty,
\end{equation*}
so that $\liminf_L \E_\Qq[\mathrm{Cost}_L] < \infty$ for $\kappa$-almost every $\Qq$. By Proposition~\ref{prop:abs-cont}, $\Qq \ll \hal(\Qq^{1} + \Qq^{-1})$ for $\kappa$-almost every $\Qq$. If $\Qq \neq \Qq^{1}$, then also $\Qq \neq \Qq^{-1}$ (otherwise, translating by $1$ would give $\Qq^{1} = \Qq$), and $\Qq$ would be singular with respect to both $\Qq^{1}$ and $\Qq^{-1}$, which are extremal: this contradicts the absolute continuity. Hence $\Qq = \Qq^{1}$ for $\kappa$-almost every $\Qq$, and $\Pp = \Pp^{1}$ after integration. We wrote the proof for the translation by $1$. For the translation by $a \in \R$, we use the maps $x \mapsto x \pm a \psi(\frac{x}{L})$, with $L \geq \max(1, 2|a| \|\psi'\|_\infty)$, and the same argument applies with constants depending on $a$. Hence $\Pp$ is stationary.
\end{proof}

\section{The DLR variational principle}
\label{sec:varprinc}
In this section, we prove Theorem~\ref{theo:VarProp}.

\paragraph{Notation.} Throughout this section, $\Pp$ is a good stationary point process. In Section~\ref{sec:petite-entropie}, we will also assume that it satisfies~\eqref{DLR}, in the sense of Definition~\ref{def:DLR}. For $N \geq 1$, we let $\La_N := [-\frac{N}{2}, \frac{N}{2}]$ be the interval of length $N$ centered at the origin, as in Section~\ref{sec:intro-free-energy}. From now on, $N$ denotes the length of a box, and is unrelated to the number of particles of the Gaussian $\beta$-ensembles of Sections~\ref{sec:electric_formalism_local_laws_and_consequences} and~\ref{sec:good}. We write $\Pp_{\La_N}$ for the law of $\bX_{\La_N}$ under $\Pp$, and $n := \bX(\La_N)$ for the number of points in $\La_N$. Since $\Pp$ is good and stationary, it has intensity one, and
\begin{equation}
\label{eq:n-moments}
\E_\Pp[n] = N, \qquad \E_\Pp\left[ (n - N)^2 \right] \leq \Cc \log^2 N, \qquad \E_\Pp\left[ n^2 \right] \leq \Cc N^2.
\end{equation}

\subsubsection*{Strategy of the proof}
We compare $\Pp_{\La_N}$ with a finite-volume Gibbs measure $\PNbGC$ in $\La_N$, whose free energy per unit volume converges to $\inf \Fb$, using the relative entropy as a distance.
\begin{enumerate}
	\item The relative entropy controls the difference of free energies: for every good stationary $\Pp$, the growth rate of $\Ent(\Pp_{\La_N} \, | \, \PNbGC)$ is at least $\Fbint(\Pp) - \inf \Fb$ (Proposition~\ref{prop:hPNbGC}). This is the finite-volume Gibbs variational principle, and it does not use the DLR equations.
	\item If $\Pp$ satisfies the DLR equations, then $\Ent(\Pp_{\La_N} \, | \, \PNbGC) = o(N)$ (Proposition~\ref{prop:petite_entropie}). Conditionally on the exterior, the configuration in $\La_N$ is a Gibbs measure, which differs from the reference measure through the move function, i.e. through the effect of the exterior on the interior. We show that this effect is $o(N)$, by writing the move function in terms of the discrepancies of the interior and exterior configurations. This is the long-range analogue of the classical boundary estimate discussed in Section~\ref{sec:literature}.
\end{enumerate}
Combining both statements gives $\Fbint(\Pp) \leq \inf \Fb$. The reference measure $\PNbGC$ is grand-canonical, with a penalization of the number of points, because $\Pp$ does not have a fixed number of points in $\La_N$.

\subsection{The reference measures}
\label{sec:reference-measures}

\subsubsection*{Non-neutral measures in a box}
For $N \geq 1$ and $M \geq 0$, let $F_{\La_N}$ be the energy~\eqref{def:FLa} in the box $\La_N$, and let $\Pgen{N}{M}{\beta}$ be the probability measure on $\Conf(\La_N, M)$ given by
\begin{equation}
\label{def:PNMb}
\dd \Pgen{N}{M}{\beta}(\bX) := \frac{1}{\Zgen{N}{M}{\beta}} \exp\left( - \beta F_{\La_N}(\bX) \right) \, \dd \Bin_{\La_N, M}(\bX),
\end{equation}
where $\Zgen{N}{M}{\beta}$ is the normalizing constant. It describes $M$ charges in $\La_N$, together with the uniform neutralizing background on $\La_N$. 

The \emph{canonical} log-gas is the neutral case $M = N$, and we write
\begin{equation}
\label{def:PNbC}
\PNbC := \Pgen{N}{N}{\beta}, \qquad \ZNbC := \Zgen{N}{N}{\beta}.
\end{equation}

The free energy enters through the following lemma.
\begin{lemma}[Thermodynamic limit]
\label{lem:thermoLimit}
For every $\beta > 0$,
\begin{equation}
\label{eq:ZNbC_asymptotics}
\lim_{N \to \infty} \frac{1}{N} \log \ZNbC = - \inf \Fb.
\end{equation}
\end{lemma} 
\begin{proof}
This is a consequence of~\cite[Corollary 1.1]{leble2017large}, see Section~\ref{sec:conventions} for the normalization. After rescaling $\La_N$ to $[-\frac12, \frac12]$, the measure $\PNbC$ is a one-dimensional log-gas whose confining potential is the logarithmic potential generated by the uniform measure on the segment, and $+ \infty$ outside the segment. Its equilibrium measure is the uniform measure on the segment.
\end{proof}

\subsubsection*{The grand-canonical reference measure}
In a \emph{grand-canonical} framework, the number of points is not fixed. Here, we define
\begin{equation}
\label{def:grand-canonical}
\dd \PNbGC(\bX) := \frac{1}{\ZNbGC}  \sum_{M=0}^\infty \exp\left( - \beta F_{\La_N}(\bX) \right) e^{- 2 \beta (M - N)^2 \log N} \, \dd \Bin_{\La_N, M}(\bX),
\end{equation}
with a grand-canonical partition function given by:
\begin{equation}
\label{def:ZNbGC}
\ZNbGC := \sum_{M=0}^\infty e^{- 2 \beta (M - N)^2 \log N} \, \Zgen{N}{M}{\beta}.
\end{equation}
One should think of the factor $e^{- 2 \beta (M - N)^2 \log N}$ as penalizing large discrepancies. It forces $M$ to be close to $N$, without changing the partition function at exponential order (Lemma~\ref{lem:ZNbGC_asymptotics}). It is of the optimal order: a non-neutral system with $M$ points in $\La_N$ gains a factor of order $e^{\frac{\beta}{2} (M - N)^2 \log N}$ in its partition function, see~\eqref{borne phi-minus}, and any constant larger than $\frac12$ in front of $\beta (M - N)^2 \log N$ would do. The series converges, since $F_{\La_N} \geq - (M + N)^2 \log N$ on $\Conf(\La_N, M)$ for $N \geq 3$, so that $\Zgen{N}{M}{\beta} \leq e^{\beta (M + N)^2 \log N}$, and $2 (M - N)^2 - (M + N)^2 \to + \infty$ as $M \to \infty$.

\begin{lemma}[Grand-canonical asymptotics]
\label{lem:ZNbGC_asymptotics}
We have
\begin{equation*}
\lim_{N \to \infty} \frac{1}{N} \log \ZNbGC = - \inf \Fb.
\end{equation*}
\end{lemma}
The lower bound is clear: keep only the term $M = N$. For the upper bound, we compare $\Zgen{N}{M}{\beta}$ (a non-neutral system in a box of size $N$) to $\Zc{M}{\beta}$ (a neutral system in a box of size $M$), by rescaling a configuration of $M$ points in $\La_N$ into a neutral configuration in $\La_M$: when $\bX \in \Conf(\La_N, M)$, we define
\begin{equation*}
\bX^{(M)} := \sum_{p \in \bX} \delta_{\frac{M}{N} p} \in \Conf(\La_M, M).
\end{equation*}
The energy cost of this rescaling is explicit, up to the fluctuation of the potential created by the background, i.e. of
\begin{equation*}
V_M(t) := \int_{\La_M} - \log|t - u| \, \dd u, \qquad \Fluct_\bY(V_M) := \sum_{p \in \bY} V_M(p) - \int_{\La_M} V_M(t) \, \dd t, \qquad \bY \in \Conf(\La_M, M),
\end{equation*}
which we control using the electric formalism, see Lemma~\ref{lem:box-electric} below. The proof is given in Appendix~\ref{sec:proof-GC}.

\subsection{Proof of the variational principle}
\label{sec:proof-VarProp}
We consider the growth rate of the relative entropy of $\Pp_{\La_N}$ with respect to the reference measure:
\begin{equation}
\label{def:hgammapi}
h(\Pp \, | \, \PNbGC) := \limsup_{N \to \infty} \frac{1}{N} \Ent\left( \Pp_{\La_N} \, | \, \PNbGC \right) \geq 0.
\end{equation}
The next proposition shows that $h(\Pp \, | \, \PNbGC)$ controls the free energy of $\Pp$.
\begin{proposition}
\label{prop:hPNbGC}
Let $\Pp$ be a good stationary point process. Then
\begin{equation}
\label{hPNbGC}
h(\Pp \, | \, \PNbGC) \geq \Fbint(\Pp) - \inf \Fb \geq \Fb(\Pp) - \inf \Fb.
\end{equation}
\end{proposition}
\begin{proof}
We prove the two inequalities separately.

\paragraph{Expansion of the relative entropy.} Let $\pi_N(M) := e^{-N} \frac{N^M}{M!}$ be the probability that $\Poi_{\La_N}$ has $M$ points. Since $\Poi_{\La_N}$ conditioned to have $M$ points is $\Bin_{\La_N, M}$, the definition~\eqref{def:grand-canonical} gives, for $\bX \in \Conf(\La_N, M)$,
\begin{equation*}
\log \frac{\dd \PNbGC}{\dd \Poi_{\La_N}}(\bX) = - \log \ZNbGC - \beta F_{\La_N}(\bX) - 2 \beta (M - N)^2 \log N - \log \pi_N(M).
\end{equation*}
Since $\PNbGC$ has a positive density with respect to $\Poi_{\La_N}$, the relative entropy with respect to $\PNbGC$ can be computed through $\Poi_{\La_N}$:
\begin{multline*}
\Ent\left( \Pp_{\La_N} \, | \, \PNbGC \right) = \E_\Pp\left[ \log \frac{\dd \Pp_{\La_N}}{\dd \Poi_{\La_N}}(\bX_{\La_N}) \right] - \E_\Pp\left[ \log \frac{\dd \PNbGC}{\dd \Poi_{\La_N}}(\bX_{\La_N}) \right] \\
= \Ent\left( \Pp_{\La_N} \, | \, \Poi_{\La_N} \right) - \E_\Pp\left[ \log \frac{\dd \PNbGC}{\dd \Poi_{\La_N}}(\bX_{\La_N}) \right].
\end{multline*}
Inserting the previous expression of the density, we get
\begin{multline}
\label{EntgPGC1}
\Ent\left( \Pp_{\La_N} \, | \, \PNbGC \right) = \Ent\left( \Pp_{\La_N} \, | \, \Poi_{\La_N} \right) + \beta \E_{\Pp}\left[ F_{\La_N}(\bX_{\La_N}) \right] + \log \ZNbGC \\
+ 2 \beta \log N \, \E_\Pp\left[ (n - N)^2 \right] + \E_\Pp\left[ \log \pi_N(n) \right].
\end{multline}
\paragraph{The Poisson terms.} The fourth term of~\eqref{EntgPGC1} is non-negative. For the last one, we use Stirling's bound and the inequality $\log u \leq u - 1$, in the form
\begin{equation*}
\log n! \leq n \log n - n + \log(n+1) + 1, \qquad n \log \frac{n}{N} - (n - N) \leq \frac{(n - N)^2}{N}.
\end{equation*}
They give
\begin{equation*}
\log \pi_N(n) \geq - \frac{(n - N)^2}{N} - \log(n + 1) - 1.
\end{equation*}
By~\eqref{eq:n-moments} and Jensen's inequality, the last term of~\eqref{EntgPGC1} is thus $o(N)$.

\paragraph{First inequality.} By stationarity and by the definition~\eqref{def:Wint},
\begin{equation*}
\lim_{N \to \infty} \frac{1}{N} \Ent\left( \Pp_{\La_N} \, | \, \Poi_{\La_N} \right) = \Ee(\Pp), \qquad \liminf_{N \to \infty} \frac{1}{N} \E_{\Pp}\left[ F_{\La_N}(\bX_{\La_N}) \right] \geq \Wint(\Pp),
\end{equation*}
and $\frac{1}{N} \log \ZNbGC \to - \inf \Fb$ by Lemma~\ref{lem:ZNbGC_asymptotics}.
Dividing~\eqref{EntgPGC1} by $N$ and taking the $\limsup$ gives the first inequality of~\eqref{hPNbGC}.

\paragraph{From the intrinsic to the electric energy.} By~\cite[Prop.~5.1]{leble2016logarithmic}, for every stationary point process $\Pp$,
\begin{equation*}
\W(\Pp) \leq \Wint(\Pp) + \Dlog(\Pp), \qquad \text{where} \qquad \Dlog(\Pp) := \Cc \limsup_{L \to \infty} \frac{\log L}{L} \E_\Pp\left[ \left( \bX(\La_L) - L \right)^2 \right]
\end{equation*}
penalizes the fluctuations of the number of points in large boxes. Here $\Pp$ is good, so $\Dlog(\Pp) = 0$. Hence $\W(\Pp) \leq \Wint(\Pp)$, and thus $\Fb(\Pp) \leq \Fbint(\Pp)$, which gives the second inequality of~\eqref{hPNbGC}.
\end{proof}

The main step is the following.
\begin{proposition}
\label{prop:petite_entropie}
Let $\Pp$ be a good stationary solution of the DLR equations. Then $h(\Pp \, | \, \PNbGC) = 0$.
\end{proposition}
The proof is given in the next section. We now explain how to deduce the DLR variational principle.

\begin{proof}[Proof of Theorem~\ref{theo:VarProp}]
Let $\Pp$ be a good stationary DLR solution. By Propositions~\ref{prop:hPNbGC} and~\ref{prop:petite_entropie}, we have:
\begin{equation*}
0 = h(\Pp \, | \, \PNbGC) \geq  \Fb(\Pp) - \inf \Fb \geq 0.
\end{equation*}
Hence $\Fb(\Pp) = \inf \Fb$.
\end{proof}

\subsection{Smallness of the relative entropy: proof of Proposition~\ref{prop:petite_entropie}}
\label{sec:petite-entropie}
Let $\Pp$ be a good stationary solution of the DLR equations. We want to bound $\Ent\left( \Pp_{\La_N} \, | \, \PNbGC \right)$ from above. We follow the proof of Proposition~\ref{prop:local-entropy}, but the reference measure is now $\PNbGC$, and we need a bound in $o(N)$ instead of $O(1)$. As in Step~3 of that proof, the kernels differ from the reference measure through the move function. We write it in terms of the discrepancies, which are controlled since $\Pp$ is good. In Claim~\ref{lem:exterior-field}, we did so for the exterior configuration only, the interior one being in a fixed window. Here, both live at the scale $N$, and we integrate by parts in both variables.

\subsubsection*{Discrepancies}
For an interval $I$ and a configuration $\bX$, let $\Dis(\bX, I) := \bX(I) - |I|$ be the discrepancy of $\bX$ in $I$, and let
\begin{equation*}
\sigma^2(x) := \E_\Pp\left[ \Dis(\bX, [0, x])^2 \right], \qquad x \geq 0.
\end{equation*}
By stationarity, $\E_\Pp[\Dis(\bX, I)^2] = \sigma^2(|I|)$ for every interval $I$. Since $\Pp$ is good, $\sigma^2(x) \leq \Cc \log^2(2 + x)$ for $x \geq 1$. We will also need a bound on $\sigma^2(x)$ which vanishes as $x \to 0$. For an integer $k \geq 1$, split $[0, 1]$ into $k$ intervals $J_1, \dots, J_k$ of length $\frac{1}{k}$. Since $\sum_i \bX(J_i)^2 \leq \left( \sum_i \bX(J_i) \right)^2 = \bX([0,1])^2$, and since all the $\bX(J_i)$ have the same law by stationarity,
\begin{equation*}
k \, \E_\Pp\left[ \bX\left(\left[0, \tfrac{1}{k}\right]\right)^2 \right] = \sum_{i=1}^k \E_\Pp\left[ \bX(J_i)^2 \right] \leq \E_\Pp\left[ \bX([0,1])^2 \right] \leq \Cc.
\end{equation*}
Monotonicity then gives $\E_\Pp[\bX([0, x])^2] \leq \Cc x$ for $x \in [0, 1]$. Thus, for all $x > 0$ we can write:
\begin{equation}
\label{eq:omega}
\sigma^2(x) \leq \Cc \min\left( x, \log^2(2 + x) \right).
\end{equation}

\paragraph{Discrepancies of the gas in a box.} We will also need to control the discrepancies of the canonical gas in a box, which is not stationary. As in Section~\ref{sec:electric}, they are controlled by the electric energy. For $M \geq 1$ and $\bY \in \Conf(\La_M, M)$, let $u_\bY := \g * \left( \sum_{p \in \bY} \delta_p - \ind_{\La_M} \delta_\R \right)$ be the electric potential generated by $\bY$ and by the uniform background on $\La_M$, as in Section~\ref{sec:electric}. Let $u_{\bY, \vec{\mathsf{r}}}$ be its truncation at the nearest-neighbor radii $\mathsf{r}_p := \frac14 \min\left(1, \min_{q \neq p} |p - q|\right)$, and
\begin{equation*}
\Elec(\bY) := \int_{\R^2} |\nabla u_{\bY, \vec{\mathsf{r}}}|^2,
\end{equation*}
which is finite since the system is neutral. Testing the equation satisfied by $u_{\bY, \vec{\mathsf{r}}}$, which is~\eqref{eq:poisson-trunc} with $\ind_{\La_M}$ instead of $\mathbf{m}_N$, integrating by parts and using the Cauchy--Schwarz inequality, we get, for every continuous $\Phi$ with compact support and $\nabla \Phi \in L^2(\R^2)$,
\begin{equation}
\label{eq:duality-box}
\left| \int \Phi \, \dd \left( \sum_{p \in \bY} \sm{p}{\mathsf{r}_p} - \ind_{\La_M} \delta_\R \right) \right| \leq \frac{1}{2\pi} \|\nabla \Phi\|_{L^2} \sqrt{\Elec(\bY)}.
\end{equation}

\begin{lemma}[Electric controls in a box]
\label{lem:box-electric}
Up to increasing the constant $\Cc_0$ of~\eqref{eq:coercivity}, there exists $\Cc$ such that, for every $M \geq 1$ and $\bY \in \Conf(\La_M, M)$:
\begin{enumerate}
	\item $F_{\La_M}(\bY) + \Cc_0 M \geq \frac{1}{8\pi} \Elec(\bY) \geq 0$,
	\item for every interval $I \subset \La_M$, $|\Dis(\bY, I)| \leq \Cc \left( 1 + (|I| + 3)^{1/4} \sqrt{\Elec(\bY)} \right)$,
	\item $|\Fluct_\bY(V_M)| \leq \frac{M}{2} + \Cc M \sqrt{\Elec(\bY)}$.
\end{enumerate}
\end{lemma}
The proof is given in Appendix~\ref{sec:proof-box-electric}.

\begin{remark}
Item 2 controls the discrepancy in $I$ by the \emph{total} electric energy, which is far from optimal: in Lemma~\ref{lem:overcrowding-energy}, one could control it by the \emph{local} energy in a neighborhood of $I$, see~\cite[Lemma A.2]{peilen2022local}. The point, however, is that we have no information on the local energies of the canonical gas in a box. Indeed, the local laws  are proven in great generality, but not for a hard-wall potential. They presumably extend to this setting, but checking it would require an inspection of the long proof of~\cite{peilen2022local}, which seems unreasonable here. Thankfully, it turns out that:
\begin{itemize}
	\item The \emph{global} energy of the canonical gas is controlled by much simpler arguments, which go back to Sandier and Serfaty~\cite{sandier2015log}.
	\item This crude bound is enough for our purposes.
\end{itemize}
Global controls for the one-dimensional log-gas are notably easier than local laws: about ten years separate the global laws of~\cite{sandier2015log} from the local laws of~\cite{peilen2022local}.
\end{remark}

\subsubsection*{The interaction with the exterior in terms of discrepancies}
As in Step~3 of the proof of Proposition~\ref{prop:local-entropy}, we write the move function as a difference of interactions with the exterior, the background being now the Lebesgue measure. Here, these interactions do not converge separately, so we truncate them. For $T > 0$, $n := \bX(\La_N)$ and $\bY \in \Conf(\La_N, n)$, let
\begin{equation*}
\mathsf{A}^T_\bX(\bY) := \iint_{\La_N \times E_T} - \log|u - v| \, \dd (\bY - \Leb)(u) \, \dd (\bX - \Leb)(v), \qquad E_T := \left[ - \tfrac{N}{2} - T, \tfrac{N}{2} + T \right] \setminus \La_N.
\end{equation*}
Since $\bY$ and $\bX_{\La_N}$ have the same number of points, subtracting the Lebesgue measure on $\La_N$ does not change the difference, and, by Definition~\ref{def:move},
\begin{equation}
\label{eq:move-as-limit}
\M_{\La_N}(\bY, \bX) = \lim_{T \to \infty} \left( \mathsf{A}^T_\bX(\bY) - \mathsf{A}^T_\bX(\bX_{\La_N}) \right).
\end{equation}
We now write $\mathsf{A}^T_\bX$ in terms of discrepancies. The boundary terms only depend on $n$ and $\bX$, so they cancel in the difference. For $x \in [0, N]$ and $y \geq 0$, consider the following intervals inside $\La_N$
\begin{equation*}
I^+_x := \left[ \tfrac{N}{2} - x, \tfrac{N}{2} \right], \qquad I^-_x := \left[ - \tfrac{N}{2}, - \tfrac{N}{2} + x \right],
\end{equation*}
and the intervals outside $\La_N$
\begin{equation*}
E^+_y := \left( \tfrac{N}{2}, \tfrac{N}{2} + y \right], \qquad E^-_y := \left[ - \tfrac{N}{2} - y, - \tfrac{N}{2} \right).
\end{equation*}
For configurations $\bY$ in $\La_N$ and $\bX$ in $\R$, let
\begin{equation}
\label{def:JN}
\mathsf{J}_N(\bY, \bX) := \sum_{\pm} \int_0^\infty \int_0^N \frac{\Dis(\bY, I^\pm_x) \, \Dis(\bX, E^\pm_y)}{(x + y)^2} \, \dd x \, \dd y.
\end{equation}

\begin{lemma}[The move function in terms of discrepancies]
\label{lem:move-discrepancies}
For $\Pp$-almost every $\bX$ and $\Bin_{\La_N, \bX(\La_N)}$-almost every $\bY$, the integrals in~\eqref{def:JN} converge absolutely, for $\bY$ and for $\bX_{\La_N}$, and we have:
\begin{equation}
\label{eq: reecriture mlambda}
\M_{\La_N}(\bY, \bX) = \mathsf{J}_N(\bY, \bX) - \mathsf{J}_N(\bX_{\La_N}, \bX).
\end{equation}
\end{lemma}
\begin{proof}
We treat the right side of $\La_N$, the left side being symmetric. The discrepancies $\Dis(\bY, I^+_x)$ and $\Dis(\bX, E^+_y)$ vanish at $x = 0$ and $y = 0$, and $\Dis(\bY, I^+_N) = n - N$. After the two changes of variables $u = \frac{N}{2} - x$ and $v = \frac{N}{2} + y$, and two integrations by parts (in $y$, then in $x$), the right part of $\mathsf{A}^T_\bX(\bY)$ can be written as
\begin{equation*}
\int_0^N \int_0^{T} - \log(x + y) \, \dd_y \Dis(\bX, E^+_y) \, \dd_x \Dis(\bY, I^+_x) = \mathrm{(i)} + \mathrm{(ii)} + \mathrm{(iii)},
\end{equation*}
where
\begin{equation*}
\mathrm{(i)} := \int_0^{T} \int_0^N \frac{\Dis(\bY, I^+_x) \, \Dis(\bX, E^+_y)}{(x + y)^2} \, \dd x \, \dd y,
\end{equation*}
\begin{equation*}
\mathrm{(ii)} := (n - N) \left( \int_0^{T} \frac{\Dis(\bX, E^+_y)}{N + y} \, \dd y - \Dis(\bX, E^+_T) \log(N + T) \right),
\end{equation*}
\begin{equation*}
\mathrm{(iii)} := \Dis(\bX, E^+_T) \int_0^N \frac{\Dis(\bY, I^+_x)}{x + T} \, \dd x.
\end{equation*}
\begin{itemize}
	\item The term (i) represents the interaction in terms of discrepancies. It converges absolutely as $T \to \infty$. Indeed, for $y \geq 1$, we have $|\Dis(\bY, I^+_x)| \leq n + N$ and, almost surely,
	\begin{equation*}
	\int_1^\infty \frac{|\Dis(\bX, E^+_y)|}{y^2} \, \dd y < \infty,
	\end{equation*}
	see the proof of Proposition~\ref{prop:move}.
	\item The term (ii) only depends on $n$ and $\bX$. It is the same for $\bY$ and for $\bX_{\La_N}$, and cancels in~\eqref{eq:move-as-limit}.
	\item The term (iii) is bounded by
	\begin{equation*}
	(n + N) N \frac{|\Dis(\bX, E^+_T)|}{T},
	\end{equation*}
	which tends to zero as $T \to \infty$, almost surely, again by the proof of Proposition~\ref{prop:move}.
\end{itemize}
Letting $T \to \infty$ in~\eqref{eq:move-as-limit} gives~\eqref{eq: reecriture mlambda}.
\end{proof}

By Lemma~\ref{lem:move-discrepancies}, the move function $\M_{\La_N}(\bY, \bX)$ and $\mathsf{J}_N(\bY, \bX)$ differ by a quantity which does not depend on $\bY$, which does not change the DLR kernel. The DLR kernel of $\Pp$ in $\La_N$ is thus, for $\Pp$-almost every $\bX$,
\begin{equation}
\label{def:kernel-N}
\mathcal{K}_N(\bX, \dd \bY) := \frac{1}{\K_N(\bX)} e^{- \beta \left( F_{\La_N}(\bY) + \mathsf{J}_N(\bY, \bX) \right)} \, \dd \Bin_{\La_N, n}(\bY),
\end{equation}
where $\K_N(\bX)$ is the normalizing constant.

\subsubsection*{Convexity and the four terms}
By~\eqref{DLR} in $\La_N$, applied to functions of $\bX_{\La_N}$, we can write $\Pp_{\La_N}$ as the mixture
\begin{equation}
\label{formule proba totale}
\Pp_{\La_N} = \int \mathcal{K}_N(\bX, \cdot) \, \Pp(\dd \bX).
\end{equation}
In particular, if $\bX$ has law $\Pp$ and, given $\bX$, $\bY$ has law $\mathcal{K}_N(\bX, \cdot)$, then $\bY$ has law $\Pp_{\La_N}$. By stationarity, we thus have, for every $x \in [0, N]$ and $y \geq 0$,
\begin{equation}
\label{eq:marginals}
\E\left[ \Dis(\bY, I^\pm_x)^2 \right] = \sigma^2(x), \qquad \E\left[ \Dis(\bX, E^\pm_y)^2 \right] = \sigma^2(y).
\end{equation}
By convexity of the relative entropy, we get
\begin{equation*}
\Ent\left( \Pp_{\La_N} \, | \, \PNbGC \right) \leq \int \Ent\left( \mathcal{K}_N(\bX, \cdot) \, | \, \PNbGC \right) \Pp(\dd \bX).
\end{equation*}
Comparing~\eqref{def:kernel-N} with~\eqref{def:grand-canonical}, we obtain, for $n = \bX(\La_N)$,
\begin{equation*}
\frac{\dd \mathcal{K}_N(\bX, \cdot)}{\dd \PNbGC}(\bY) = \frac{\ZNbGC}{\K_N(\bX)} \, e^{- \beta \mathsf{J}_N(\bY, \bX)} \, e^{2 \beta (n - N)^2 \log N} \, \ind_{\bY(\La_N) = n}.
\end{equation*}
We split the constant as
\begin{equation*}
\frac{\ZNbGC}{\K_N(\bX)} = \frac{\Zgen{N}{n}{\beta}}{\K_N(\bX)} \cdot \frac{\ZNbGC}{\Zgen{N}{n}{\beta}},
\end{equation*}
and we obtain
\begin{equation}
\label{eq:four-terms}
\Ent\left( \Pp_{\La_N} \, | \, \PNbGC \right) \leq \mathcal{T}_1 + \mathcal{T}_2 + \mathcal{T}_3 + \mathcal{T}_4,
\end{equation}
where
\begin{equation*}
\mathcal{T}_1 := - \beta \, \E_\Pp\left[ \E_{\mathcal{K}_N(\bX, \cdot)}\left[ \mathsf{J}_N(\bY, \bX) \right] \right], \qquad \mathcal{T}_2 := \E_\Pp\left[ \log \frac{\Zgen{N}{n}{\beta}}{\K_N(\bX)} \right],
\end{equation*}
\begin{equation*}
\mathcal{T}_3 := \E_\Pp\left[ \log \frac{\ZNbGC}{\Zgen{N}{n}{\beta}} \right], \qquad \mathcal{T}_4 := 2 \beta \log N \, \E_\Pp\left[ (n - N)^2 \right].
\end{equation*}

Recall that $h(\Pp \, | \, \PNbGC) = \limsup_N \frac{1}{N} \Ent(\Pp_{\La_N} \, | \, \PNbGC)$ is non-negative, see~\eqref{def:hgammapi}. In order to show that it vanishes, it is thus enough to bound each of the four terms from above by $o(N)$.
\begin{lemma}[The four terms are small]
\label{lem:four-terms}
For $N$ large enough,
\begin{equation*}
\mathcal{T}_1 \leq \Cc \sqrt{N} \log(2 + N), \qquad \mathcal{T}_2 \leq \Cc N^{3/4} \log^2(2 + N), \qquad \limsup_{N \to \infty} \frac{\mathcal{T}_3}{N} \leq 0, \qquad \mathcal{T}_4 \leq \Cc \log^3 N.
\end{equation*}
\end{lemma}
The proof is a computation, given in Appendix~\ref{sec:proof-four-terms}.

\paragraph{How the four terms are bounded.}
\begin{itemize}
	\item The terms $\mathcal{T}_1$ and $\mathcal{T}_2$ are of the same nature. Indeed, by Jensen's inequality,
	\begin{equation*}
	\mathcal{T}_2 = - \E_\Pp\left[ \log \E_{\Pgen{N}{n}{\beta}}\left[ e^{- \beta \mathsf{J}_N(\bY, \bX)} \right] \right] \leq \beta \, \E_\Pp\left[ \E_{\Pgen{N}{n}{\beta}}\left[ \mathsf{J}_N(\bY, \bX) \right] \right].
	\end{equation*}
	Both are boundary terms: we bound the interaction $\mathsf{J}_N$ between the interior and the exterior, with its representation~\eqref{def:JN} in terms of discrepancies and the Cauchy--Schwarz inequality. In $\mathcal{T}_1$, the interior configuration $\bY$ follows the DLR kernel, so that its discrepancies are those of $\Pp$, controlled by~\eqref{eq:omega}. In $\mathcal{T}_2$, it follows the canonical gas $\Pgen{N}{n}{\beta}$, whose discrepancies are controlled by Lemma~\ref{lem:box-electric}.
	\item The term $\mathcal{T}_3$ compares the grand-canonical and the canonical partition functions. This is not completely obvious, since non-neutral configurations could change the partition function a lot, but the penalization in~\eqref{def:grand-canonical} makes the relevant numbers of points close to $N$, and we can apply Lemma~\ref{lem:ZNbGC_asymptotics}.
	\item The term $\mathcal{T}_4$ is directly controlled by~\eqref{eq:n-moments}.
\end{itemize}

\begin{proof}[Proof of Proposition~\ref{prop:petite_entropie}]
By~\eqref{eq:four-terms} and Lemma~\ref{lem:four-terms}, we have
\begin{equation*}
h(\Pp \, | \, \PNbGC) = \limsup_{N \to \infty} \frac{1}{N} \Ent\left( \Pp_{\La_N} \, | \, \PNbGC \right) \leq 0,
\end{equation*}
hence $h(\Pp \, | \, \PNbGC) = 0$.
\end{proof}

\appendix

\section{Auxiliary proofs}
\label{sec:auxiliary}

\subsection{How to read Theorem~\ref{theo:bounded-fluct} from~\cite{peilen2022local}}
\label{sec:how_to_read_theorem_ref_theo_bounded_fluct}

By~\cite[Theorem 2, (1.28)]{peilen2022local}, the log-Laplace transform of $\Fluct_N(\xi)$ on the good event is equal to $\frac{s^2}{\beta} \|\theta\|_{H^{1/2}}^2$ up to an error of order
\begin{equation*}
|s| \left( 1 + \left|1 - \tfrac{1}{\beta}\right| \right) \|\theta\|_{C^3} + \frac{s^2}{\beta} \max\left( \|\theta\|^2_{C^1}, \|\theta\|^3_{C^3} \right) + \frac{|s|^3}{\beta^2 \ell} \max\left( \|\theta\|^2_{C^2}, \|\theta\|^3_{C^3} \right),
\end{equation*}
which is bounded for $|s| \leq 1$.

The implicit constants, as well as the constants in the bound $\PNb(\G^c) \leq \Cc e^{- \cc \beta \ell}$, depend neither on $z$ nor on $\ell$.

The theorem holds for $N \geq N_c$, where the threshold $N_c$ of~\cite[Definition 4.7]{peilen2022local} requires that $\supp \xi$ is contained in the bulk, that $L$ is less than half the length of the support of the equilibrium measure (here $L < 2$), and that $\|t \psi'\|_\infty < \frac12$, where $\psi$ is the transport map associated to $\xi$ and $|t| \leq \Cc \frac{|s|}{N}$. By the transport estimates of~\cite[Lemma 4.8]{peilen2022local}, $|\psi'| \leq \frac{\Cc}{L}$, so that the last condition holds for $|s| \leq 1$ as soon as $\ell$ is larger than a constant, which we include in $\omega$.

The implicit constants of~\cite[Lemma 4.8]{peilen2022local} may depend on the size of the support of $\theta$, which is why we state Theorem~\ref{theo:bounded-fluct} for a fixed $\theta$: in Section~\ref{sec:number-variance}, we only apply it to finitely many fixed test functions.

\subsection{Truncation of finite-volume move functions: proof of Proposition~\ref{prop:finite-truncation}}
\label{sec:proof-finite-truncation}
We use the notation of Section~\ref{sec:finite-DLR-section}. For $|x| \leq r$ and $|y| \geq 2r$, a Taylor expansion gives
\begin{equation}
\label{eq:kx}
|k_x(y)| \leq \frac{2r}{|y|}, \qquad |\partial_y k_x(y)| \leq \frac{4r}{y^2}, \qquad \left| k_x(y) - \frac{x}{y} \right| \leq \frac{2r^2}{y^2}.
\end{equation}

\begin{proof}[Proof of Claim~\ref{lem:exterior-field}]
We treat the bulk part and the edge part separately.

\paragraph{Bulk.} We can focus on $T < y \leq c_3 N$ and treat the other side similarly. As in~\cite{dereudre2021dlr}, we integrate by parts, which makes the discrepancies $\nu_N((T, y])$ appear:
\begin{equation*}
\int_{(T, c_3 N]} k_x(y) \, \dd \nu_N(y) = k_x(c_3 N) \, \nu_N((T, c_3 N]) - \int_T^{c_3 N} \partial_y k_x(y) \, \nu_N((T, y]) \, \dd y.
\end{equation*}
In the unscaled coordinates, $\nu_N((T, y]) = \bXN(J) - \mathbf{m}_N(J)$ with $J := \vrho^{-1} (T, y] \subset [0, \frac{N}{2}]$, of length at most $c_1 N$, by the choice of $c_3$. Hence, by the control~\eqref{eq:variance} on the variance, and by~\eqref{eq:superexp} if $|J| < 1$, we have $\E_N[\nu_N((T, y])^2]^{\hal} \leq \Cc \log(2 + y)$. Inserting~\eqref{eq:kx} to bound $\partial_y k_x$, we obtain
\begin{equation*}
\E_N\left[ \sup_{|x| \leq r} \left| \int_{(T, c_3 N]} k_x \, \dd \nu_N \right|^2 \right]^{\hal} \leq \frac{2r}{c_3 N} \Cc \log N + \int_T^\infty \frac{4r}{y^2} \Cc \log(2 + y) \, \dd y \leq \Cc r \left( \frac{\log N}{N} + \frac{\log(2 + T)}{T} \right).
\end{equation*}

\paragraph{Edge.} By the last bound of~\eqref{eq:kx}, and since the total variation of $\nu_N$ is at most $2N$, the edge part is a dipole term up to a small error:
\begin{equation*}
\left| V_{\mathrm{edge}}(x) - x \int_{|y| > c_3 N} \frac{\dd \nu_N(y)}{y} \right| \leq \frac{2 r^2}{(c_3 N)^2} \cdot 2N \leq \frac{\Cc r^2}{N}.
\end{equation*}
The remaining integral is a \emph{macroscopic} linear statistic: the test function $y \mapsto \frac{\ind_{|y| > c_3 N}}{y}$ varies on the scale $N$. We go back to the macroscopic coordinates $v_i := \frac{x_i}{N} = \frac{y_i}{\vrho N}$. By the law of large numbers for $\beta$-ensembles, see e.g.~\cite{benarous1997large, johansson1998fluctuations}, the empirical measure $\frac{1}{N} \sum_{i=1}^N \delta_{v_i}$ converges almost surely to the semicircle law $\rhosc(v) \dd v$ of~\eqref{def:rhosc}. Since the image of $\widetilde{\mu}_N$ is $N \rhosc(v) \dd v$, we get
\begin{equation*}
\int_{|y| > c_3 N} \frac{\dd \nu_N(y)}{y} = \frac{1}{\vrho} \int_{|v| > c_3/\vrho} \frac{1}{v} \, \dd \left( \frac{1}{N} \sum_{i=1}^N \delta_{v_i} - \rhosc(v) \dd v \right).
\end{equation*}
The function $v \mapsto \frac{1}{v} \ind_{|v| > c_3/\vrho}$ is bounded by $\frac{\vrho}{c_3}$, so the integral is bounded by $\frac{2}{c_3}$. This function is continuous except at $\pm \frac{c_3}{\vrho}$, which carry no mass for the semicircle law, hence the integral goes to zero almost surely. Multiplying by $|x| \leq r$ and taking the supremum gives the claims on $V_{\mathrm{edge}}$.
\end{proof}

\begin{proof}[Proof of Proposition~\ref{prop:finite-truncation}]
Let $n := \bXNz(\La)$ and $\epsilon > 0$. By Proposition~\ref{prop:FiniteVolEstimates}, we can choose $n_0$ such that $\P_N(n > n_0) \leq \frac{\epsilon}{2}$ for all $N$ large enough. Let $\eta := \frac{\epsilon}{4 n_0}$. By~\eqref{eq:move-potential}, on the event $\{n \leq n_0\}$, the supremum in the proposition is larger than $\epsilon$ only if $\sup_{|x| \leq r} |V^T_{\mathrm{bulk}}(x)| > \eta$ or $\sup_{|x| \leq r} |V_{\mathrm{edge}}(x)| > \eta$. Hence
\begin{equation*}
\P_N\left( \sup_{\bX'} \left| \M^{(N)}_\La - \M^{T,(N)}_\La \right| > \epsilon \right) \leq \frac{\epsilon}{2} + \frac{1}{\eta^2} \E_N\left[ \sup_{|x| \leq r} |V^T_{\mathrm{bulk}}(x)|^2 \right] + \P_N\left( \sup_{|x| \leq r} |V_{\mathrm{edge}}(x)| > \eta \right).
\end{equation*}
By~\eqref{eq:edge-potential}, the last term tends to zero as $N \to \infty$. By~\eqref{eq:exterior-field}, the limsup of the second one is at most
\begin{equation*}
\frac{\Cc^2 r^2}{\eta^2} \frac{\log^2(2 + T)}{T^2},
\end{equation*}
which is smaller than $\frac{\epsilon}{2}$ for $T \geq T_0$, with $T_0$ depending on $\epsilon$, $r$ and $n_0$.
\end{proof}

\subsection{Existence of move functions}
\label{sec:existence_of_move_functions}
This is the same argument as for the bulk part of Claim~\ref{lem:exterior-field}, with $\Leb$ instead of $\widetilde{\mu}_N$. Since there is no edge, this is all we need.

\begin{proof}[Proof of Proposition~\ref{prop:move}]
Let $r \geq 1$ be such that $\La \subset [-r, r]$, let $n := \bX(\La)$, and let $\bX' \in \Conf(\La, n)$.

\paragraph{Reduction to the origin.} The part of the integral defining $\M^T_\La$ where $y \in [-2r, 2r] \setminus \La$ does not depend on $T \geq 2r$. For the rest, as in~\eqref{eq:move-potential}, we may replace $- \log |x - y|$ by $k_x(y)$ since $\bX'$ and $\bX_\La$ have the same number of points. It is thus enough to show that, almost surely, the potential
\begin{equation*}
W_T(x) := \int_{2r < |y| \leq T} k_x(y) \, \dd (\bX - \Leb)(y)
\end{equation*}
converges as $T \to \infty$, uniformly in $|x| \leq r$.

\paragraph{Integration by parts.} We consider the positive side, the negative one being symmetric. Let
\begin{equation*}
D(y) := \Dis(\bX, (2r, y]) = \bX((2r, y]) - (y - 2r)
\end{equation*}
be the discrepancy, with the notation $\Dis$ of Section~\ref{sec:petite-entropie}. Then
\begin{equation*}
\int_{(2r, T]} k_x(y) \, \dd (\bX - \Leb)(y) = k_x(T) \, D(T) - \int_{2r}^T \partial_y k_x(y) \, D(y) \, \dd y.
\end{equation*}
Since $\Pp$ is good, $\E_\Pp[|D(y)|] \leq \Cc \log(2 + y)$. By~\eqref{eq:kx}, we get
\begin{equation*}
\E_\Pp\left[ \int_{2r}^\infty \sup_{|x| \leq r} |\partial_y k_x(y)| \, |D(y)| \, \dd y \right] < \infty,
\end{equation*}
so the last integral converges absolutely as $T \to \infty$, almost surely and uniformly in $|x| \leq r$.

\paragraph{Boundary term.} For $T \in [m, m+1]$ with $m$ an integer, $|D(T)| \leq |D(m)| + \bX([m, m+1]) + 1$. Since
\begin{equation*}
\sum_m \frac{\E_\Pp[D(m)^2] + \E_\Pp[\bX([m, m+1])^2]}{m^2} < \infty,
\end{equation*}
the Borel--Cantelli lemma shows that $\frac{D(T)}{T} \to 0$ almost surely, hence
\begin{equation*}
\sup_{|x| \leq r} |k_x(T) D(T)| \leq \frac{2r |D(T)|}{T} \longrightarrow 0.
\end{equation*} This proves the existence of $\M_\La(\bX', \bX)$.

\paragraph{Tail estimate.} Letting $T' \to \infty$ in the integration by parts between $T$ and $T'$, we get, as in~\eqref{eq:move-potential},
\begin{equation*}
\sup_{\bX' \in \Conf(\La, n)} \left| \M_\La(\bX', \bX) - \M^T_\La(\bX', \bX) \right| \leq 2n \sum_\pm \left( \frac{2r |D_\pm(T)|}{T} + \int_T^\infty \frac{4r}{y^2} |D_\pm(y)| \, \dd y \right),
\end{equation*}
where $D_\pm$ are the discrepancies on the positive and negative sides. The expectation of the bracket is at most $\Cc r \frac{\log(2 + T)}{T}$. Since $n$ is finite almost surely, this gives~\eqref{eq:move-tail}.
\end{proof}

\subsection{Compactness in the local topology: proof of Lemma~\ref{lem:entropy-compactness}}
\label{sec:entropy-compactness}
For $m \geq 1$, recall that $\La_m = [-\frac{m}{2}, \frac{m}{2}]$, and let
\begin{equation*}
\varphi^k_m := \frac{\dd \Pp^k_{\La_m}}{\dd \Poi_{\La_m}}
\end{equation*}
be the local densities.

\paragraph{Compactness of the densities.} For each $m$, since $t \mapsto t \log t$ is superlinear, the entropy bound implies that $(\varphi^k_m)_k$ is uniformly integrable in $L^1(\Poi_{\La_m})$, hence relatively sequentially compact for the weak topology $\sigma(L^1, L^\infty)$. By a diagonal extraction, there is a subsequence $(k_j)_j$ such that, for every $m$, $\varphi^{k_j}_m$ converges weakly in $L^1(\Poi_{\La_m})$ to some $\varphi_m$, which is non-negative with integral one.

\paragraph{The limit process.} The measures $\varphi_m \Poi_{\La_m}$ are consistent: the restriction from $\La_{m+1}$ to $\La_m$ corresponds, at the level of densities, to a conditional expectation with respect to $\Poi_{\La_{m+1}}$, which is continuous for the weak $L^1$ topology. By Kolmogorov's extension theorem, there is a point process $\Pp$ with $\Pp_{\La_m} = \varphi_m \Poi_{\La_m}$ for all $m$. If $f$ is local, bounded and measurable, say $f(\bX) = f(\bX_{\La_m})$, then
\begin{equation*}
\E_{\Pp^{k_j}}[f] = \int f \varphi^{k_j}_m \, \dd \Poi_{\La_m} \longrightarrow \int f \varphi_m \, \dd \Poi_{\La_m} = \E_\Pp[f].
\end{equation*}

\subsection{DLR equations for limit points: proof of item 2 of Theorem~\ref{theo:Convergence}}
\label{sec:proof-DLR-limit}
We compare the kernels with the following elementary lemma.

\begin{lemma}[Comparison of kernels]
\label{lem:comparison-kernels}
Let $w_1, w_2 > 0$ be integrable weights on a probability space, and assume that $e^{- \beta \delta} \leq \frac{w_1}{w_2} \leq e^{\beta \delta}$ for some $\delta \geq 0$. Then, for every bounded measurable $f$,
\begin{equation*}
\left| \frac{\int f w_1}{\int w_1} - \frac{\int f w_2}{\int w_2} \right| \leq \left( e^{2\beta\delta} - 1 \right) \|f\|_\infty.
\end{equation*}
\end{lemma}
\begin{proof}
Write $w_1 = \lambda w_2$, and let
\begin{equation*}
\bar{\lambda} := \frac{\int \lambda w_2}{\int w_2}.
\end{equation*}
Both $\lambda$ and $\bar{\lambda}$ lie in $[e^{-\beta\delta}, e^{\beta\delta}]$, and the difference is equal to
\begin{equation*}
\frac{\int f (\lambda - \bar{\lambda}) w_2}{\bar{\lambda} \int w_2},
\end{equation*}
which is bounded by $(e^{\beta \delta} - e^{- \beta \delta}) e^{\beta \delta} \|f\|_\infty$.
\end{proof}
In particular, multiplying the weights by a constant does not change the kernel, so the terms of the energy which do not depend on the configuration $\bX'$ in $\La$ play no role below.

\begin{lemma}[Flattening the background]
\label{lem:flatten}
Let $\La \subset [-r, r]$ be a bounded interval and $T \geq r$. For every $n_0$ and every bounded measurable $f$,
\begin{equation*}
\lim_{N \to \infty} \sup\left\{ \left| f^{T, (N)}_\La(\bX) - f^T_\La(\bX) \right| \, : \, \bX \in \Conf, \ \bX(\La) \leq n_0 \right\} = 0.
\end{equation*}
\end{lemma}
\begin{proof}
Let $\bX' \in \Conf(\La, n)$ with $n \leq n_0$. The difference $F^{\rho_N}_\La(\bX') - F_\La(\bX')$ is the sum of a term which does not depend on $\bX'$ and of
\begin{equation*}
\sum_{x \in \bX'} \int_\La \log|x - t| \, (\rho_N(t) - 1) \, \dd t, \qquad \text{which is bounded by } \Cc_r \, n \sup_{[-r,r]} |\rho_N - 1| \leq \Cc_r \, n \frac{r^2}{N^2}.
\end{equation*}
Similarly,
\begin{equation*}
\left| \M^{T,(N)}_\La(\bX', \bX) - \M^T_\La(\bX', \bX) \right| \leq \Cc_T \, n \sup_{[-T,T]} |\rho_N - 1| \leq \Cc_T \, n \frac{T^2}{N^2}.
\end{equation*}
We conclude by Lemma~\ref{lem:comparison-kernels}.
\end{proof}

\begin{proof}[Proof of item 2 of Theorem~\ref{theo:Convergence}]
Let $\Pp$ be the limit of the laws of $\bXNz$ along a subsequence $(N_k)$. By Theorem~\ref{theo:tightness}, the convergence holds in the exp-local topology, and $\Pp$ is good, so that the move functions exist $\Pp$-almost surely (Proposition~\ref{prop:move}). Let $\La \subset [-r, r]$ be a bounded interval. By a monotone class argument, it is enough to prove that $\E_\Pp[f] = \E_\Pp[f_\La]$ for $f$ local, bounded and measurable, say depending only on the configuration in $[-T_1, T_1]$, with $\|f\|_\infty \leq 1$.

\paragraph{Decomposition.} Let $T \geq \max(T_1, r)$. With $N = N_k$, we write
\begin{multline*}
\E_\Pp[f] - \E_\Pp[f_\La] = \underbrace{\left( \E_\Pp[f] - \E_N[f] \right)}_{(1)} + \underbrace{\left( \E_N[f] - \E_N[f^{(N)}_\La] \right)}_{(2)} + \underbrace{\left( \E_N[f^{(N)}_\La] - \E_N[f^{T,(N)}_\La] \right)}_{(3)} \\
+ \underbrace{\left( \E_N[f^{T,(N)}_\La] - \E_N[f^T_\La] \right)}_{(4)} + \underbrace{\left( \E_N[f^T_\La] - \E_\Pp[f^T_\La] \right)}_{(5)} + \underbrace{\left( \E_\Pp[f^T_\La] - \E_\Pp[f_\La] \right)}_{(6)}.
\end{multline*}

\paragraph{The six terms.}
\begin{itemize}
	\item \emph{Terms (1) and (5): local convergence.} Since $f$ and $f^T_\La$ are local, bounded and measurable, they tend to zero as $k \to \infty$.
	\item \emph{Term (2): finite-volume DLR equations.} It vanishes, by Proposition~\ref{prop:finite-DLR}.
	\item \emph{Terms (3) and (6): truncation.} Let $\delta > 0$. On the event where the move function and its truncation differ by at most $\delta$, uniformly in $\bX'$, Lemma~\ref{lem:comparison-kernels} bounds the difference of the kernels by $e^{2\beta\delta} - 1$, and elsewhere it is at most $2$. By Proposition~\ref{prop:finite-truncation} (with $\epsilon = \delta$) for (3), and by~\eqref{eq:move-tail} for (6), both terms are thus at most $e^{2\beta \delta} - 1 + 2\delta$ for $T$ large enough, and then $k$ large enough.
	\item \emph{Term (4): flattening.} It is at most
	\begin{equation*}
	\sup_{\bX(\La) \leq n_0} \left| f^{T,(N)}_\La - f^T_\La \right| + 2 \sup_N \P_N\left( \bXNz(\La) > n_0 \right),
	\end{equation*}
	which is small for $n_0$ large, by Proposition~\ref{prop:FiniteVolEstimates}, and then $k$ large, by Lemma~\ref{lem:flatten}.
\end{itemize}
Choosing $\delta$ small, then $T$ large, then $n_0$ large, then $k$ large, we obtain $\E_\Pp[f] = \E_\Pp[f_\La]$.
\end{proof}

\subsection{Proof of Lemma~\ref{lem:box-electric}: electric controls in a box}
\label{sec:proof-box-electric}
\begin{proof}[Proof of Lemma~\ref{lem:box-electric}]
Items 1 and 2 are essentially contained in~\cite[Appendix A]{peilen2022local}. Item 1 is~\cite[Lemma A.4]{peilen2022local}, with $\Omega = \R$ and the background $\ind_{\La_M}$, whose density is at most $1$. Item 2 follows from~\cite[Lemma A.2]{peilen2022local}, whose proof only uses the Poisson equation, the fact that the truncation radii are at most $\frac14$, and an upper bound on the density of the background, and thus applies here, up to the edges of $\La_M$.

\paragraph{Item 2.} For an interval $J$ of length $L \geq 3$,~\cite[Lemma A.2]{peilen2022local} gives: either $|\Dis(\bY, J)| \leq 4$, or
\begin{equation*}
\Elec(\bY) \geq \cc \, \Dis(\bY, J)^2 \min\left( 1, \sqrt{\frac{|\Dis(\bY, J)|}{L}} \right).
\end{equation*}
In the second case, we distinguish two sub-cases.
\begin{itemize}
	\item If $|\Dis(\bY, J)| \geq L$, then $|\Dis(\bY, J)| \leq \Cc \sqrt{\Elec(\bY)}$.
	\item If $|\Dis(\bY, J)| < L$, then $|\Dis(\bY, J)|^{5/2} \leq \Cc L^{1/2} \Elec(\bY)$, i.e. $|\Dis(\bY, J)| \leq \Cc L^{1/5} \Elec(\bY)^{2/5}$.
\end{itemize} Since $L^{1/5} \Elec^{2/5} = (L^{1/4} \Elec^{1/2})^{4/5} \leq \max(1, L^{1/4} \Elec^{1/2})$, we get in all cases
\begin{equation*}
|\Dis(\bY, J)| \leq \Cc \left( 1 + L^{1/4} \sqrt{\Elec(\bY)} \right).
\end{equation*}
This gives item 2 if $|I| \geq 3$. If $|I| < 3$, then $\Dis(\bY, I) \geq -3$, and, applying the previous bound to an interval $J \supset I$ of length $3$, we get $\Dis(\bY, I) \leq \bY(J) \leq 3 + \Cc (1 + 3^{1/4} \sqrt{\Elec(\bY)})$.

\paragraph{Item 3.} We use the duality inequality~\eqref{eq:duality-box}. Let $\chi(s) := \min(1, (2 - |s|)_+)$. The fluctuation does not change if a constant is added to $V_M$, since $\bY$ has $M$ points and $|\La_M| = M$. Let
\begin{equation*}
\widetilde{V}(z) := \int_{\La_M} - \log|z - t| \, \dd t + M \log M = M V_1\left( \frac{z}{M} \right), \qquad z \in \R^2,
\end{equation*}
where $V_1(w) := \int_{-1/2}^{1/2} - \log|w - s| \, \dd s$ for $w \in \R^2$. The function $V_1$ is continuous, bounded on $\{1 \leq |w| \leq 2\}$, and $\nabla V_1 \in L^2_{\mathrm{loc}}(\R^2)$: it is bounded, except for logarithmic singularities at $(\pm \frac12, 0)$. Let $\Psi(z) := \widetilde{V}(z) \chi(\frac{|z|}{M})$. By scaling,
\begin{equation*}
\|\nabla \Psi\|^2_{L^2} \leq 2 M^2 \int_{\{|w| \leq 2\}} |\nabla V_1|^2 + 8\pi M^2 \sup_{1 \leq |w| \leq 2} |V_1(w)|^2 \leq \Cc M^2.
\end{equation*}
The smeared charges lie in $\{|z| \leq M\}$, where $\Psi = \widetilde{V}$. By the mean value property of the logarithm, the average of $\widetilde{V}$ on a circle of radius $\mathsf{r} \leq \frac14$ centered at $p \in \La_M$ is
\begin{equation*}
\widetilde{V}(p) + \int_{\La_M \cap \{|t - p| < \mathsf{r}\}} \log \frac{|p - t|}{\mathsf{r}} \, \dd t \in \left[ \widetilde{V}(p) - 2 \mathsf{r}, \widetilde{V}(p) \right].
\end{equation*}
Hence the left-hand side of~\eqref{eq:duality-box} with $\Phi = \Psi$ is equal to $\Fluct_\bY(V_M)$, up to an error of at most
\begin{equation*}
\sum_{p \in \bY} 2 \mathsf{r}_p \leq \frac{M}{2},
\end{equation*}
which proves item 3.
\end{proof}

\subsection{Proof of Lemma~\ref{lem:ZNbGC_asymptotics}: grand-canonical asymptotics}
\label{sec:proof-GC}
Recall the rescaling $\bX \mapsto \bX^{(M)}$, the potential $V_M$ of the background on $\La_M$, and the fluctuation $\Fluct_{\bY}(V_M)$, from Section~\ref{sec:reference-measures}. By scaling, since
\begin{equation*}
\iint_{[0,1]^2} - \log|a - b| \, \dd a \, \dd b = \frac{3}{2},
\end{equation*}
we get
\begin{equation}
\label{eq:EM}
\int_{\La_M} V_M(t) \, \dd t = \iint_{\La_M^2} - \log|t - u| \, \dd t \, \dd u = - M^2 \log M + \frac{3}{2} M^2.
\end{equation}

\begin{proof}[Proof of Lemma~\ref{lem:ZNbGC_asymptotics}]
We proceed in four steps.

\paragraph{Step 1: the rescaling formula.} We show that, for $M \geq 1$ and $\bX \in \Conf(\La_N, M)$,
\begin{equation}
\label{eq:rescaling}
F_{\La_N}(\bX) = F_{\La_M}\left(\bX^{(M)}\right) + \varphi_N(M) - \left( \frac{N}{M} - 1 \right) \Fluct_{\bX^{(M)}}(V_M),
\end{equation}
where
\begin{equation*}
\varphi_N(M) := \frac{M}{2} \log \frac{N}{M} - \frac{(N - M)^2}{2} \left( \log N - \frac32 \right).
\end{equation*}
Let $s := \frac{N}{M}$. Expanding the square in~\eqref{def:FLa},
\begin{equation*}
2 F_{\La_N}(\bX) = \sum_{p \neq q \in \bX} - \log|p - q| - 2 \sum_{p \in \bX} \int_{\La_N} - \log|p - t| \, \dd t + \iint_{\La_N^2} - \log|t - u| \, \dd t \, \dd u.
\end{equation*}
We change variables $p = s p'$, $t = s t'$, $u = s u'$ in each of the three terms. Since $- \log|s a| = - \log|a| - \log s$,
\begin{equation*}
\sum_{p \neq q \in \bX} - \log|p - q| = \sum_{p \neq q \in \bX^{(M)}} - \log|p - q| - M(M-1) \log s,
\end{equation*}
\begin{equation*}
\sum_{p \in \bX} \int_{\La_N} - \log|p - t| \, \dd t = s \sum_{p \in \bX^{(M)}} V_M(p) - s M^2 \log s,
\end{equation*}
\begin{equation*}
\iint_{\La_N^2} - \log|t - u| \, \dd t \, \dd u = s^2 \iint_{\La_M^2} - \log|t - u| \, \dd t \, \dd u - s^2 M^2 \log s.
\end{equation*}
We subtract the same expansion of $2 F_{\La_M}(\bX^{(M)})$, and write $\sum_{p \in \bX^{(M)}} V_M(p) = \Fluct_{\bX^{(M)}}(V_M) + \int_{\La_M} V_M$. By~\eqref{eq:EM}, this gives
\begin{multline*}
2 F_{\La_N}(\bX) - 2 F_{\La_M}\left(\bX^{(M)}\right) \\
= \left( M - (s-1)^2 M^2 \right) \log s + (s - 1)^2 \left( \frac32 M^2 - M^2 \log M \right) - 2 (s - 1) \Fluct_{\bX^{(M)}}(V_M),
\end{multline*}
which is~\eqref{eq:rescaling}, since $(s-1) M = N - M$.

\paragraph{Bound on $\varphi_N$.} The inequality $\log u \leq u - 1$ gives
\begin{equation*}
0 \leq M \log \frac{N}{M} \leq N - M \quad \text{if } M \leq N, \qquad 0 \leq M \log \frac{M}{N} \leq \frac{M (M - N)}{N} \leq (M - N) + \frac{(M-N)^2}{N} \quad \text{if } M \geq N.
\end{equation*}
Since $|M - N| \leq (M - N)^2$ for integers, there is $\Cc$ such that, for every $N \geq 3$ and $M \geq 1$,
\begin{equation}
\label{borne phi}
|\varphi_N(M)| \leq \Cc (M - N)^2 \log N.
\end{equation}
More precisely, since $\frac{M}{2} \log \frac{M}{N} \leq (M - N)^2$ in both cases, the definition of $\varphi_N$ gives the one-sided bound
\begin{equation}
\label{borne phi-minus}
- \varphi_N(M) \leq \frac{1}{2} (M - N)^2 \log N + (M - N)^2.
\end{equation}

\paragraph{Step 2: partition functions.} We show that, for $M \geq 1$,
\begin{equation}
\label{eq: N to M}
\Zgen{N}{M}{\beta} = e^{- \beta \varphi_N(M)} \, \Zc{M}{\beta} \, \E_{\Pc{M}{\beta}}\left[ \exp\left( \beta \left( \frac{N}{M} - 1 \right) \Fluct_{\bY}(V_M) \right) \right].
\end{equation}
The image of $\Bin_{\La_N, M}$ by $\bX \mapsto \bX^{(M)}$ is $\Bin_{\La_M, M}$, since both are laws of $M$ independent uniform points. Hence, by~\eqref{eq:rescaling},
\begin{equation*}
\Zgen{N}{M}{\beta} = e^{- \beta \varphi_N(M)} \int \exp\left( - \beta F_{\La_M}(\bY) + \beta \left( \frac{N}{M} - 1 \right) \Fluct_{\bY}(V_M) \right) \, \dd \Bin_{\La_M, M}(\bY),
\end{equation*}
and the right-hand side is equal to that of~\eqref{eq: N to M} by definition of $\Pc{M}{\beta} = \Pgen{M}{M}{\beta}$.

\paragraph{Step 3: the fluctuation term.} We bound the last factor of~\eqref{eq: N to M}: there is $\Cc$ such that, for every $M \geq 1$ and $N \geq 1$,
\begin{equation}
\label{borne1}
\left| \log \E_{\Pc{M}{\beta}}\left[ \exp\left( \beta \left( \frac{N}{M} - 1 \right) \Fluct_{\bY}(V_M) \right) \right] \right| \leq \Cc \left( |N - M| \sqrt{M} + (N - M)^2 \right).
\end{equation}
The same bound holds, with another constant, if $\beta$ is replaced by a fixed multiple of $\beta$ in front of $\Fluct_{\bY}(V_M)$.

By item 3 of Lemma~\ref{lem:box-electric},
\begin{equation}
\label{borne fluct VN}
\left| \Fluct_{\bY}(V_M) \right| \leq \Cc M + \Cc M \sqrt{\Elec(\bY)}.
\end{equation}

We now take exponential moments. By~\eqref{eq:ZNbC_asymptotics} at $\beta$ and $\frac{\beta}{2}$,
\begin{equation*}
\E_{\Pc{M}{\beta}}\left[ e^{\frac{\beta}{2} (F_{\La_M} + \Cc_0 M)} \right] = e^{\frac{\beta}{2} \Cc_0 M} \, \frac{\Zc{M}{\frac{\beta}{2}}}{\Zc{M}{\beta}} \leq e^{\Cc M},
\end{equation*}
hence, by item 1 of Lemma~\ref{lem:box-electric}, with $\alpha := \frac{\beta}{16 \pi}$,
\begin{equation}
\label{eq:exp-moment-energy}
\E_{\Pc{M}{\beta}}\left[ e^{\alpha \Elec} \right] \leq e^{\Cc M}.
\end{equation}
Let $\lambda > 0$ and $K \geq 1$. Young's inequality and Jensen's inequality for the concave function $u \mapsto u^{1/K}$ give
\begin{equation*}
\lambda \sqrt{\Elec} \leq \frac{\alpha \Elec}{K} + \frac{K \lambda^2}{4 \alpha}, \qquad \E\left[ e^{\alpha \Elec / K} \right] \leq \E\left[ e^{\alpha \Elec} \right]^{1/K}.
\end{equation*}
By~\eqref{eq:exp-moment-energy}, we obtain
\begin{equation*}
\log \E_{\Pc{M}{\beta}}\left[ e^{\lambda \sqrt{\Elec}} \right] \leq \frac{K \lambda^2}{4 \alpha} + \frac{\Cc M}{K},
\end{equation*}
and the choice $K := \max\left( 1, \frac{2 \sqrt{\alpha \Cc M}}{\lambda} \right)$ gives
\begin{equation*}
\log \E_{\Pc{M}{\beta}}\left[ e^{\lambda \sqrt{\Elec}} \right] \leq \Cc \left( \lambda \sqrt{M} + \lambda^2 \right).
\end{equation*}
By~\eqref{borne fluct VN}, the expectation in~\eqref{borne1} is at most
\begin{equation*}
e^{\Cc \beta |N - M|} \, \E_{\Pc{M}{\beta}}\left[ e^{\lambda \sqrt{\Elec}} \right], \qquad \text{with} \qquad \lambda := \Cc \beta |N - M|,
\end{equation*}
hence its logarithm is at most $\Cc \left( |N - M| \sqrt{M} + (N - M)^2 \right)$, which gives the upper bound in~\eqref{borne1}. The lower bound in~\eqref{borne1} follows from the upper bound applied to $- \Fluct_{\bY}(V_M)$, since $\E[e^{Z}] \, \E[e^{-Z}] \geq 1$ by the Cauchy--Schwarz inequality.

\paragraph{Comparison of the partition functions.} Let
\begin{equation*}
\mathrm{err}_N(M) := (M - N)^2 \log N + |M - N| \sqrt{M}.
\end{equation*}
By~\eqref{eq: N to M},~\eqref{borne phi} and~\eqref{borne1}, there is $\Cc$ such that, for every $N \geq 3$ and $M \geq 1$,
\begin{equation}
\label{eq:Z-comparison}
\left| \log \Zgen{N}{M}{\beta} - \log \Zc{M}{\beta} \right| \leq \Cc \, \mathrm{err}_N(M).
\end{equation}

\paragraph{Step 4: conclusion.} The term $M = N$ in~\eqref{def:ZNbGC} gives $\ZNbGC \geq \ZNbC$, hence, by~\eqref{eq:ZNbC_asymptotics},
\begin{equation*}
\liminf_{N \to \infty} \frac{1}{N} \log \ZNbGC \geq - \inf \Fb.
\end{equation*}
For the upper bound, let $\epsilon > 0$. By~\eqref{eq:ZNbC_asymptotics}, there is $\Cc_\epsilon$ such that
\begin{equation*}
\log \Zc{M}{\beta} \leq - M \inf \Fb + \epsilon M + \Cc_\epsilon \qquad \text{for every } M \geq 1.
\end{equation*}
Let $M \geq 1$ and $t := |M - N|$. By~\eqref{eq: N to M},~\eqref{borne phi-minus} and~\eqref{borne1}, and since $- M \inf \Fb + \epsilon M \leq - N \inf \Fb + \epsilon N + \Cc t$,
\begin{equation*}
\log \left( e^{- 2 \beta t^2 \log N} \, \Zgen{N}{M}{\beta} \right) \leq - N \inf \Fb + \epsilon N + \Cc_\epsilon - \frac{3 \beta}{2} t^2 \log N + \Cc \left( t^2 + t \sqrt{M} + t \right).
\end{equation*}
If $M > 2N$, then $t \sqrt{M} \leq \sqrt{2} \, t^{3/2} \leq 2 t^2$. Hence, for $N$ large enough, the right-hand side is at most
\begin{equation*}
- N \inf \Fb + \epsilon N + \Cc_\epsilon - \beta t^2 \log N + \Cc t \sqrt{N} \leq - N \inf \Fb + \epsilon N + \Cc_\epsilon - \frac{\beta t^2 \log N}{2} + \frac{\Cc N}{\log N},
\end{equation*}
where we maximized $- \frac{\beta t^2 \log N}{2} + \Cc t \sqrt{N}$ over $t \geq 0$. The term $M = 0$ is bounded directly by $e^{- 2 \beta N^2 \log N + \beta N^2 \log N} \leq 1$. Summing over $M$, we obtain, for every $\epsilon > 0$,
\begin{equation*}
\limsup_{N \to \infty} \frac{1}{N} \log \ZNbGC \leq - \inf \Fb + \epsilon.
\end{equation*}
\end{proof}

\subsection{Proof of Lemma~\ref{lem:four-terms}: the four terms}
\label{sec:proof-four-terms}
We bound the four terms of~\eqref{eq:four-terms} one by one. Recall that $n := \bX(\La_N)$, and that the moments of $n$ are controlled by~\eqref{eq:n-moments}.

\subsubsection*{The term $\mathcal{T}_1$}
We show that $\mathcal{T}_1 \leq \Cc \sqrt{N} \log(2 + N)$. Let $\bX$ have law $\Pp$ and, given $\bX$, let $\bY$ have law $\mathcal{K}_N(\bX, \cdot)$. By~\eqref{def:JN}, the Cauchy--Schwarz inequality and~\eqref{eq:marginals},
\begin{equation*}
\mathcal{T}_1 \leq 2 \beta \int_0^N \int_0^\infty \frac{\sigma(x) \, \sigma(y)}{(x + y)^2} \, \dd y \, \dd x.
\end{equation*}
By~\eqref{eq:omega}, $\sigma^2(y) \leq \Cc y$, and the change of variables $y = x u$ gives
\begin{equation*}
\int_0^\infty \frac{\sqrt{y}}{(x + y)^2} \, \dd y = \frac{1}{\sqrt{x}} \int_0^\infty \frac{\sqrt{u}}{(1 + u)^2} \, \dd u = \frac{\pi}{2 \sqrt{x}}.
\end{equation*}
Hence, by~\eqref{eq:omega} again,
\begin{equation*}
\mathcal{T}_1 \leq \Cc \int_0^N \frac{\sigma(x)}{\sqrt{x}} \, \dd x \leq \Cc \left( 1 + \int_1^N \frac{\log(2 + x)}{\sqrt{x}} \, \dd x \right) \leq \Cc \sqrt{N} \log(2 + N).
\end{equation*}

\subsubsection*{The term $\mathcal{T}_2$}
We show that $\mathcal{T}_2 \leq \Cc N^{3/4} \log^2(2 + N)$. By~\eqref{def:PNMb} and~\eqref{def:kernel-N}, we have
\begin{equation*}
\frac{\K_N(\bX)}{\Zgen{N}{n}{\beta}} = \E_{\Pgen{N}{n}{\beta}}\left[ e^{- \beta \mathsf{J}_N(\bY, \bX)} \right],
\end{equation*}
and Jensen's inequality gives
\begin{equation}
\label{eq:T2-jensen}
\log \frac{\Zgen{N}{n}{\beta}}{\K_N(\bX)} \leq \beta \, \E_{\Pgen{N}{n}{\beta}}\left[ \mathsf{J}_N(\bY, \bX) \right] \leq \beta \sum_\pm \int_0^\infty \int_0^N \frac{\E_{\Pgen{N}{n}{\beta}}\left[ |\Dis(\bY, I^\pm_x)| \right] |\Dis(\bX, E^\pm_y)|}{(x + y)^2} \, \dd x \, \dd y.
\end{equation}
Unlike for $\mathcal{T}_1$, $\bY$ now has law $\Pgen{N}{n}{\beta}$ instead of the DLR kernel. The discrepancies $\Dis(\bX, E^\pm_y)$ are still controlled by $\sigma$, but those of $\bY$ are not, and we control them with the electric energy.

\paragraph{Interior discrepancies.} We show that, for every $x \in [0, N]$,
\begin{equation}
\label{eq:discrepancy-energy}
|\Dis(\bY, I^\pm_x)| \leq \Cc + 3 |n - N| + \Cc N^{1/4} \sqrt{\Elec(\bY^{(n)})} \, \ind_{2 \leq n \leq 2N}.
\end{equation}
If $n \leq 1$ or $n > 2N$, we use the trivial bound $|\Dis(\bY, I^\pm_x)| \leq n + N$, which is at most $3 + 3|n - N|$ in both cases. If $2 \leq n \leq 2N$, recall from Section~\ref{sec:reference-measures} that
\begin{equation*}
\bY^{(n)} := \sum_{p \in \bY} \delta_{\frac{n}{N} p} \in \Conf(\La_n, n)
\end{equation*}
is the configuration obtained by multiplying the points of $\bY$ by $\frac{n}{N}$. This rescaling maps $I^\pm_x$ onto an interval $I' \subset \La_n$ of length $\frac{x n}{N} \leq 2N$, so that
\begin{equation*}
\Dis(\bY, I^\pm_x) = \Dis(\bY^{(n)}, I') + x \left( \frac{n}{N} - 1 \right), \qquad \left| x \left( \frac{n}{N} - 1 \right) \right| \leq |n - N|.
\end{equation*}
Since $\bY^{(n)}$ is neutral, Lemma~\ref{lem:box-electric} gives
\begin{equation*}
|\Dis(\bY^{(n)}, I')| \leq \Cc \left( 1 + (2N + 3)^{1/4} \sqrt{\Elec(\bY^{(n)})} \right),
\end{equation*}
which proves~\eqref{eq:discrepancy-energy}. Here $I'$ touches the edge of $\La_n$, which is allowed in Lemma~\ref{lem:box-electric}.

\paragraph{Moments of the electric energy.} We show that, for $n \geq 2$, with $\mathrm{err}_N(n)$ as in~\eqref{eq:Z-comparison}
\begin{equation}
\label{eq:energy-moment-N}
\E_{\Pgen{N}{n}{\beta}}\left[ \sqrt{\Elec(\bY^{(n)})} \right] \leq \Cc \sqrt{\mathrm{err}_N(n) + n}.
\end{equation}
The law of $\bY^{(n)}$ under $\Pgen{N}{n}{\beta}$ has a density with respect to the canonical measure $\Pc{n}{\beta}$, given by
\begin{equation}
\label{eq:density-rescaled}
\frac{e^{Z}}{\E_{\Pc{n}{\beta}}\left[ e^{Z} \right]}, \qquad \text{where} \qquad Z(\bY') := \beta \left( \frac{N}{n} - 1 \right) \Fluct_{\bY'}(V_n), \qquad \bY' \in \Conf(\La_n, n).
\end{equation}
Indeed, for a bounded measurable $f$, by the rescaling formula~\eqref{eq:rescaling} and since the image of $\Bin_{\La_N, n}$ by $\bY \mapsto \bY^{(n)}$ is $\Bin_{\La_n, n}$,
\begin{equation*}
\E_{\Pgen{N}{n}{\beta}}\left[ f(\bY^{(n)}) \right] = \frac{e^{- \beta \varphi_N(n)}}{\Zgen{N}{n}{\beta}} \int f(\bY') \, e^{- \beta F_{\La_n}(\bY') + Z(\bY')} \, \dd \Bin_{\La_n, n}(\bY') = \frac{e^{- \beta \varphi_N(n)} \, \Zc{n}{\beta}}{\Zgen{N}{n}{\beta}} \, \E_{\Pc{n}{\beta}}\left[ f e^{Z} \right].
\end{equation*}
 Let $\alpha := \frac{\beta}{16 \pi}$ be as in~\eqref{eq:exp-moment-energy}. By Jensen's inequality, then the Cauchy--Schwarz inequality,
\begin{equation*}
\E_{\Pgen{N}{n}{\beta}}\left[ \sqrt{\Elec(\bY^{(n)})} \right]^2 \leq \frac{2}{\alpha} \log \frac{\E_{\Pc{n}{\beta}}\left[ e^{\frac{\alpha}{2} \Elec + Z} \right]}{\E_{\Pc{n}{\beta}}\left[ e^{Z} \right]} \leq \frac{2}{\alpha} \log \frac{\E_{\Pc{n}{\beta}}\left[ e^{\alpha \Elec} \right]^{1/2} \E_{\Pc{n}{\beta}}\left[ e^{2 Z} \right]^{1/2}}{\E_{\Pc{n}{\beta}}\left[ e^{Z} \right]}.
\end{equation*}
By~\eqref{eq:exp-moment-energy}, $\E_{\Pc{n}{\beta}}[e^{\alpha \Elec}] \leq e^{\Cc n}$. The two other factors are bounded by~\eqref{borne1}, with $\beta$ and $2 \beta$ in front of $\Fluct$. This proves~\eqref{eq:energy-moment-N}.

\paragraph{Conclusion.} Since
\begin{equation*}
\int_0^N \frac{\dd x}{(x + y)^2} = \frac{N}{y (N + y)} \qquad \text{and} \qquad |n - N| \leq \sqrt{\mathrm{err}_N(n)},
\end{equation*}
inserting~\eqref{eq:discrepancy-energy} and~\eqref{eq:energy-moment-N} into~\eqref{eq:T2-jensen} gives
\begin{equation*}
\log \frac{\Zgen{N}{n}{\beta}}{\K_N(\bX)} \leq \Cc N^{5/4} \left( \sqrt{\mathrm{err}_N(n) + n} + 1 \right) \sum_\pm \int_0^\infty \frac{|\Dis(\bX, E^\pm_y)|}{y (N + y)} \, \dd y.
\end{equation*}
We take the expectation under $\Pp$ and apply the Cauchy--Schwarz inequality and~\eqref{eq:marginals}:
\begin{equation*}
\mathcal{T}_2 \leq \Cc N^{5/4} \left( \E_\Pp\left[ \mathrm{err}_N(n) + n \right]^{1/2} + 1 \right) \int_0^\infty \frac{\sigma(y)}{y (N + y)} \, \dd y.
\end{equation*}
By~\eqref{eq:omega}, the last integral is at most
\begin{equation*}
\frac{\Cc}{N} \left( \int_0^1 \frac{\dd y}{\sqrt{y}} + \int_1^N \frac{\log(2 + y)}{y} \, \dd y \right) + \Cc \int_N^\infty \frac{\log(2 + y)}{y^2} \, \dd y \leq \Cc \frac{\log^2(2 + N)}{N}.
\end{equation*}
By~\eqref{eq:n-moments} and the Cauchy--Schwarz inequality, $\E_\Pp[|n - N| \sqrt{n}] \leq \E_\Pp[(n-N)^2]^{1/2} \E_\Pp[n]^{1/2} \leq \Cc \sqrt{N} \log N$, so that
\begin{equation}
\label{eq:psi-moment}
\E_\Pp[\mathrm{err}_N(n)] \leq \Cc \left( \log^3 N + \sqrt{N} \log N \right), \qquad \E_\Pp[\mathrm{err}_N(n) + n] \leq \Cc N.
\end{equation}
Hence $\mathcal{T}_2 \leq \Cc N^{3/4} \log^2(2 + N)$.

\subsubsection*{The term $\mathcal{T}_3$}
We show that $\limsup_N \frac{\mathcal{T}_3}{N} \leq 0$. For $n \geq 1$, by~\eqref{eq:Z-comparison},
\begin{equation*}
\log \Zgen{N}{n}{\beta} \geq \log \Zc{n}{\beta} - \Cc \mathrm{err}_N(n).
\end{equation*}
On the event $\{n = 0\}$, whose probability is at most $\Cc \frac{\log^2 N}{N^2}$ by~\eqref{eq:n-moments}, we use $\Zgen{N}{0}{\beta} \geq 1$ (the self-energy of the background on $\La_N$ is negative for $N$ large) and $\log \ZNbGC \leq \Cc N$, so that this event contributes $O(\frac{\log^2 N}{N})$ to $\mathcal{T}_3$. By~\eqref{eq:ZNbC_asymptotics} and Lemma~\ref{lem:ZNbGC_asymptotics}, for every $\epsilon > 0$, there is $\Cc_\epsilon$ such that
\begin{equation*}
\log \Zc{n}{\beta} \geq - n \inf \Fb - \epsilon n - \Cc_\epsilon, \qquad \log \ZNbGC \leq - N \inf \Fb + \epsilon N + \Cc_\epsilon.
\end{equation*}
Since $\E_\Pp[n] = N$, using~\eqref{eq:psi-moment}, we get
\begin{equation*}
\mathcal{T}_3 \leq \E_\Pp\left[ (n - N) \inf \Fb + \epsilon (n + N) + \Cc \mathrm{err}_N(n) \right] + 2 \Cc_\epsilon + o(1) \leq 2 \epsilon N + \Cc \sqrt{N} \log N + 2 \Cc_\epsilon + o(1).
\end{equation*}
Since $\epsilon$ is arbitrary, $\limsup_N \frac{\mathcal{T}_3}{N} \leq 0$.

\subsubsection*{The term $\mathcal{T}_4$}
By~\eqref{eq:n-moments}, $\mathcal{T}_4 \leq \Cc \log^3 N$.

\bibliographystyle{alpha}
\bibliography{ExistenceSinebeta}

\end{document}